\documentclass{amsart}
\usepackage{amssymb,amscd,verbatim}
\usepackage{epsfig}
\usepackage[all]{xy}
\newcommand\datver[1]{\def\datverp%
 {\par\boxed{\boxed{\text{#1; Run: \today}}}}}
\datver{1.1A; Revised: 23-5-2007}

\newcommand{\ie}{i\@.e\@.\ }

\newcommand\Cinf{\CI}
\newcommand\dCinf{\dCI}
\newcommand\dCI{\dot{\mathcal{C}}^\infty}
\newcommand\CI{\mathcal{C}^\infty}
\newcommand\CmI{\mathcal{C}^{-\infty}}
\newcommand\Lap{\Delta}

\newcommand\pa{\partial}
\newcommand\Mand{\text{ and }}

\newcommand\Mwith{\text{ with }}

\newcommand\Mwhere{\text{ where }}

\newcommand\ep{\epsilon}
\newcommand\ev{\sigma}

\newcommand\Sp{\operatorname{Sp}}
\newcommand\fraks{\mathfrak{s}}
\newcommand\frakh{\mathfrak{h}}
\newcommand\frakg{\mathfrak{g}}

\newcommand\Ran{\operatorname{Ran}}
\newcommand\cR{\mathcal{R}}
\newcommand\cL{\mathcal{L}}
\newcommand\cO{\mathcal{O}}
\newcommand\cP{\mathcal{P}}
\newcommand\cJS{\mathcal{J}_S}
\newcommand\cI{\mathcal{I}}

\newcommand\cS{\mathcal{S}}

\newcommand\cM{\mathcal{M}}
\newcommand\tcM{\widetilde{\mathcal{M}}}

\renewcommand\inf{\operatorname{inf}}
\renewcommand\sup{\operatorname{sup}}
\newcommand\supp{\operatorname{supp}}
\newcommand\pp{\operatorname{pp}}
\newcommand\ess{\operatorname{ess}}
\newcommand\mic{\operatorname{mic}}

\newcommand\WF{\operatorname{WF}}
\newcommand\WFsc{\WF_{\text{sc}}}
\newcommand\WFscp{\WF'_{\text{sc}}}
\newcommand\Hsc{H_{\text{sc}}}
\newcommand\Nat{\mathbb{N}}
\newcommand\Real{\mathbb{R}}
\newcommand\Cx{\mathbb{C}}
\newcommand\im{\operatorname{Im}}
\newcommand\re{\operatorname{Re}}
\newcommand\bbR{\mathbb{R}}

\newcommand\RR{\mathbb{R}}
\newcommand\Rp{q}
\newcommand\RP{\operatorname{RP}}

\newcommand\Cp{q}

\newcommand\nin{\not\in}

\newcommand\Span{\operatorname{span}}

\newcommand\Cv{\operatorname{Cv}}

\newcommand\Min{\operatorname{Min}}

\newcommand\Vb{{\mathcal V}_{\bl}}

\newcommand\bl{{\operatorname{b}}}
\newcommand\scT{{}^{\scl} T}

\newcommand\scl{{\operatorname{sc}}}

\newcommand\Psisc{\Psi_{\scl}}

\newcommand\Dist{C^{-\infty}}
\newcommand\dist{C^{-\infty}}
\newcommand\CIdot{\dot C^{\infty}}

\newcommand\scH{\kern3pt{}^{\scl} H}
\newcommand\Isc{I_{\scl}}

\renewcommand\Im{\operatorname{Im}}

\def\Id{\operatorname{Id}}

\newtheorem{lemma}{Lemma}[section]
\newtheorem{prop}[lemma]{Proposition}
\newtheorem{thm}[lemma]{Theorem}
\newtheorem{cor}[lemma]{Corollary}

\newtheorem*{thm*}{Main Results}
\newtheorem*{theorem*}{Theorem}
\newtheorem*{prop*}{Proposition}
\numberwithin{equation}{section}
\theoremstyle{remark}
\newtheorem{rem}[lemma]{Remark}
\theoremstyle{definition}

\theoremstyle{definition}
\newtheorem{Def}[lemma]{Definition}

\newcommand\Xsch{X_{\operatorname{Sch}}}

\newcommand\half{\frac{1}{2}}

\newcommand\NN{\mathbb{N}}

\newcommand\Pt{\tilde P}

\newcommand\inc{-}
\newcommand\out{+}
\newcommand\Vy{{V_0}}

\newcommand\Four{\mathcal{F}}

\newcommand\symp{\mathfrak{sp}}

\newcommand\effnr{\operatorname{enr}}
\newcommand\effr{\operatorname{er}}
\newcommand\hesst{\operatorname{Ht}}

\newcommand\norm{\operatorname{norm}}

\newcommand\Rbar{\overline{\mathbb{R}}}
\newcommand\Schwartz{\mathcal{S}}
\newcommand\ac{\operatorname{ac}}

\datver{0.9Z; Revised: 6-10-2003}

\begin{document}

\author[Andrew Hassell]{Andrew Hassell$^*$}
\address{Department of Mathematics, Australian National
   University, Canberra ACT 0200 Australia}
\email{hassell@maths.anu.edu.au}
\author[Richard Melrose]{Richard Melrose$^\dag$}
\address{Department of Mathematics, Massachusetts Institute of Technology,
Cambridge MA 02139}
\email{rbm@math.mit.edu}
\author[Andr\'as Vasy]{Andr\'as Vasy$^{\ddag}$}
\address{Department of Mathematics, Massachusetts Institute of Technology,
Cambridge MA 02139 and Northwestern University, Evanston IL 60208}
\curraddr{Department of Mathematics, Stanford University, Stanford, CA
94305-2125, U.S.A.}
\email{andras@math.stanford.edu}

\title[Microlocal propagation near radial points]
{Microlocal propagation near radial points and scattering for symbolic potentials of order zero}


\thanks{$^*$ Supported in part by an Australian Research Council Fellowship,
$^\dag$ Supported in part by the National Science Foundation under grant
\#DMS-0408993,
$^\ddag$ Supported in part by the National Science Foundation under grant
\#DMS-0201092, a Clay 
Research Fellowship and a Fellowship from the Alfred P.\,Sloan Foundation.}

\begin{abstract} In this paper, the scattering and spectral theory of
  $H=\Lap_g+V$ is developed, where $\Lap_g$ is the Laplacian with respect
  to a scattering metric $g$ on a compact manifold $X$ with boundary and
  $V\in\CI(X)$ is real; this extends our earlier results in the
  two-dimensional case. Included in this class of operators are
  perturbations of the Laplacian on Euclidean space by potentials
  homogeneous of degree zero near infinity. Much of the particular
  structure of geometric scattering theory can be traced to the occurrence
  of radial points for the underlying classical system. In this case the
  radial points correspond precisely to critical points of the restriction,
  $V_0,$ of $V$ to $\partial X$ and under the additional assumption that
  $V_0$ is Morse a functional parameterization of the generalized
  eigenfunctions is obtained.
  
  The main subtlety of the higher dimensional case arises from additional
  complexity of the radial points. A normal form near such points obtained
  by Guillemin and Schaeffer is extended and refined, allowing a microlocal
  description of the null space of $H-\ev$ to be given for all but a
  finite set of `threshold' values of the energy; additional complications arise at the discrete
  set of `effectively resonant' energies. It is shown that each critical
  point at which the value of $V_0$ is less than $\ev$ is the source of
  solutions of $Hu=\ev u.$ The resulting description of the generalized
  eigenspaces is a rather precise, distributional, formulation of
  asymptotic completeness. We also derive the closely related $L^2$ and
  time-dependent forms of asymptotic completeness, including the absence of
  $L^2$ channels associated with the non-minimal critical points.  This
  phenomenon, observed by Herbst and Skibsted, can be attributed to the
  fact that the eigenfunctions associated to the non-minimal critical
  points are `large' at infinity; in particular they are too large to lie
  in the range of the resolvent $R(\ev \pm i0)$ applied to compactly
  supported functions.

\end{abstract}

\maketitle

\tableofcontents

\section{Introduction}

In this paper, which is a continuation of \cite{HMV1}, scattering
theory is developed for symbolic potentials of order zero. The general
setting is the same as in \cite{HMV1}, consisting of a compact
manifold with boundary, $X,$ equipped with a scattering metric, $g,$
and a real potential, $V\in\CI(X).$ Recall that such a scattering
metric on $X$ is a smooth metric in the interior of $X$ taking the
form
\begin{equation}
g = \frac{dx^2}{x^4} + \frac{h}{x^2}
\label{sc-metric}\end{equation}
near the boundary, where $x$ is a boundary defining function and $h$ is a
smooth cotensor which restricts to a metric on $\{x=0\}=\partial X.$ This
makes the interior, $X^\circ,$ of $X$ a complete manifold which is
asymptotically flat and is metrically asymptotic to the large end of a
cone, since in terms of the singular normal coordinate $r =x^{-1},$ the
leading part of the metric at the boundary takes the form $dr^2 + r^2
h(y,dy).$ In the compactification of $X^\circ$ to $X,$ $\pa X$ corresponds
to the set of asymptotic directions of geodesics. In particular, this setting
subsumes the case of the standard metric on Euclidean space, or a compactly
supported perturbation of it, with a potential which is a classical symbol
of order zero, hence not decaying at infinity but rather with leading term
which is asymptotically homogeneous of degree zero. The study of the
scattering theory for such potentials was initiated by Herbst
\cite{Herbst1}.

Let $V_0 \in \CI(\partial X)$ be the restriction of $V$ to $\partial X,$
and denote by $\Cv(V)$ the set of critical values of $V_0.$ It is shown in
\cite{HMV1} that the operator $H=\Lap_g+V$ (where the Laplacian is
normalized to be positive) is essentially self-adjoint with continuous
spectrum occupying $[\min V_0,\infty)$. There may be discrete spectrum of
finite multiplicity in $(-\min_XV, \max V_0]$ with possible accumulation
points only at $\Cv(V).$ To obtain finer results, it is natural to assume,
as we do throughout this paper unless otherwise noted, that $V_0$ is a
Morse function, \ie has only nondegenerate critical points; in particular
$\Cv(V)$ is a then finite set; by definition this is the set of
\emph{threshold energies,} or \emph{thresholds}.

Fr{o}m the microlocal point of view scattering theory is largely about the
study of radial points, i.e. the points in the cotangent bundle where the
Hamilton vector field is a multiple of the radial vector field (i.e. the
vector field $A = \sum_i z_i \partial_{z_i}$ on Euclidean space, where
$(z_1, \dots, z_n) \in \RR^n$). These correspond in the classical dynamical
system to the places where the particle is moving either in purely incoming
or outgoing sense.  In scattering theory for potentials decaying at
infinity, there is a radial point for each point on the sphere at infinity;
thus there is a manifold of radial points and the behaviour of the flow in
a neighbourhood of these points is rather simple, either attracting (at the
outgoing radial surface) or repelling (at the incoming radial surface) in
the transverse direction.  Estimates involving commutation with the radial
vector field $A$ multiplied by suitable powers of $|z|$ and perhaps
additional microlocalizing operators, are usually sufficient to control the
behaviour of generalized eigenfunctions. These are known as Mourre-type
estimates and play a fundamental role in conventional scattering theory.
In the present case, assuming $V_0$ is a Morse function, the radial points
are isolated and occur in pairs, one pair (incoming/outgoing) for each
critical point of $V_0$. The linearized Hamiltonian flow at the radial
points is rather more complicated since it depends on the Hessian of $V_0$
at the critical point, which is arbitrary apart from being
nondegenerate. This makes the higher dimensional case more intricate than
the case $\dim X = 2$ which we treated in \cite{HMV1}. Correspondingly one
needs more elaborate commutator estimates in order to control the behaviour
of generalized eigenfunctions. We give a rather general and complete
analysis of the regularity of solutions of $Pu = 0$ in a microlocal
neighbourhood of a radial point of $P$, using the concept of a test module
of operators. This is a family of pseudodifferential operators which is a
module over the zero-order operators, contains $P$, and is closed under
commutation. By choosing a test module closely tailored to the Hamilton
flow of $P$ near the radial point we are able to produce enough
positive-commutator estimates to parametrize the microlocal solutions of
$Pu = 0$. The construction of appropriate test modules (which can be
thought of as simply an effective bookkeeping device for keeping track of a
rather intricate set of commutator estimates) to analyze general radial
points is the main technical innovation of this paper.

The general study of radial points was initiated by Guillemin and Schaeffer
\cite{MR55:3504}. This was done in a slightly different context, where $P$
is a standard pseudodifferential operator with homogeneous principal symbol
and a radial point is one where the Hamilton vector field is a multiple of
the vector field $\sum_i \xi_i \partial_{\xi_i}$ generating dilations in
the cotangent space. This setting is completely equivalent to ours, via
conjugation by a `local Fourier Transform' (see Section~\ref{lft}). They
analyzed the situation in the nonresonant case. We refine their analysis by
treating the resonant case, which is crucial in our application since we
have a family of operators parametrized by the energy level, and the
closure of the set of energies which give rise to resonant radial points
may have nonempty interior. Moreover, we show that our parametrization of
microlocal solutions is smooth except at a set of `effectively resonant'
energies which is always discrete.

Recently, Bony, Fujiie, Ramond and Zerzeri have studied the microlocal kernel of 
pseudodifferential operators at a hyperbolic fixed point \cite{BFRZ}, 
corresponding, in our setting, to a radial point associated to a 
local maximum of $V_0.$ Their results partially overlap ours, being
most closely related to Section 10 of \cite{HMV1} and \cite{HMV1EDP}. 

\subsection{Previous results} The Euclidean setting described above was
first studied by Herbst \cite{Herbst1}, who showed that any finite energy
solution of the time dependent Schr\"odinger equation, so $u = e^{-itH} f$
with $f\in L^2(\bbR^n),$ can concentrate, in an $L^2$ sense, asymptotically
as $t\to\infty$ only in directions which are critical points of $V_0.$ This
was subsequently refined by Herbst and Skibsted \cite{Herbst-Skibsted2}, 
who showed that such
concentration can only occur near local minima of $V_0.$ In contrast,
solutions of the classical flow can concentrate near any critical point of
$V_0.$

Asymptotic completeness has been studied by Agmon, Cruz and Herbst
\cite{Agmon-Cruz-Herbst1}, by Herbst and Skibsted \cite{Herbst-Skibsted1},
\cite{Herbst-Skibsted2}, \cite{Herbst-Skibsted3}
and the present authors in \cite{HMV1}. Agmon, Cruz
and Herbst showed asymptotic completeness for sufficiently high energies,
while Herbst and Skibsted extended this to all energies except for an
explicitly given union of bounded intervals; in the two dimensional case,
they showed asymptotic completeness for all energies. These results were
obtained by time-dependent methods. On the other hand the principal result
of \cite{HMV1} involves a precise description of the generalized
eigenspaces of $H$ 
\begin{equation}
E^{-\infty}(\ev)=\{u\in\CmI(X);(H-\ev)u=0\};
\label{HMV2r.195}\end{equation}
note that the space of `extendible distributions' $\CmI(X)$ is the analogue
of tempered distributions and reduces to it in case $X$ is the radial
compactification of $\bbR^n.$ Thus we are studying all \emph{tempered}
eigenfunctions of $H.$ Let us recall these results in more detail.

For any $\ev\notin\Cv(V)$ the space $E_{\pp}(\ev)$ of $L^2$ eigenfunctions
is finite dimensional, and reduces to zero except for $\ev$ in a discrete
(possibly empty) subset of $[\min_XV,\max V_0]\setminus\Cv(V).$ It is
always the case that $E_{\pp}(\ev)\subset\dCI(X)$ consists of rapidly
decreasing functions. Hence $E^{-\infty}_{\ess}(\ev)\subset
E^{-\infty}(\ev),$ the orthocomplement of $E_{\pp}(\ev),$ is well defined for
$\ev\notin\Cv(V).$ Furthermore, as shown in the Euclidean case by Herbst in
\cite{Herbst1}, the resolvent, $R(\ev)$ of $H,$ acting on this
orthocomplement, has a limit, $R(\ev\pm i0),$ on $[\min
  V_0,\infty)\setminus\Cv(V)$ from above and below. The subspace of
  `smooth' eigenfunctions is then defined as
\begin{equation}\begin{gathered}
E^{\infty}_{\ess}(\ev)=\Sp(\ev)\left(\CIdot(X) \ominus E_{\pp}(\ev)\right)
\subset E^{-\infty}(\ev)\\
\Sp(\ev) \equiv \frac1{2\pi i} \big( R(\ev + i0) - R(\ev - i0) \big) .
\end{gathered}\label{spsp}\end{equation}
In fact 
\begin{equation*}
E^{\infty}_{\ess}(\ev)\subset \bigcap_{\epsilon >0}x^{-1/2 - \epsilon} L^2(X).
\label{HMV2r.196}\end{equation*}
An alternative characterization of $E^{\infty}_{\ess}(\ev)$ can be given in
terms of the \emph{scattering wavefront set} at the boundary of $X$.

The scattering cotangent bundle, $\scT^*X,$ of $X$ is naturally isomorphic
to the cotangent bundle over the interior of $X,$ and indeed globally
isomorphic to $T^*X$ by a non-natural isomorphism; the natural identification
exhibits both `compression' and `rescaling' at the boundary. If $(x,y)$
are local coordinates near a boundary point of $X$, with $x$ a boundary
defining function, then linear coordinates $(\nu, \mu)$ are defined on the
scattering cotangent bundle by requiring that $\Rp\in\scT^*X$ be written as
\begin{equation}
\Rp=-\nu\,\frac{dx}{x^2}+\sum_i\mu_i\,\frac{dy_i}{x}, \quad \nu \in \RR,
\quad \mu \in \RR^{n-1}.
\label{numu}\end{equation}
This makes $(\nu, \mu)$ dual to the basis $(-x^2 \partial_x, x
\partial_{y_i})$ of vector fields which form an approximately unit length
basis, uniformly up to the boundary, for any scattering metric. In
Euclidean space, $\nu$ is dual to $\partial_r$ and $\mu_i$ is dual to the
constant-length angular derivative $r^{-1} \partial_{y_i}.$ In the analysis
of the microlocal aspects of $H-\ev,$ in part for compatibility with
\cite{MR55:3504}, it is convenient pass to an operator `of first order' by
multiplying $H-\ev$ by $x^{-1}$, \ie to replace it by
\begin{equation*}
P=P(\ev)=x^{-1}(H-\ev).
\end{equation*}

The classical dynamical system giving the behaviour of particles,
asymptotically near $\pa X,$ moving under the influence of the potential
corresponds to `the bicharacteristic vector field,' see
\eqref{eq:scHp-local}, determined by the \emph{boundary symbol,} $p,$ of
$P.$ This vector field is defined on $\scT^*_{\partial X} X$, which is to
say on $\scT^* X$ at, and tangent to, the boundary $\scT^*_{\partial X} X =
\scT^*X \cap \{x=0\}.$ It has the property that $\nu$ is nondecreasing
under the flow; we refer to points $(y, \nu, \mu)$ where $\mu = 0$ as
\emph{incoming} if $\nu < 0$ and \emph{outgoing} if $\nu > 0$. What is
important in understanding the behaviour of the null space of $P,$ \ie
tempered distributions, $u,$ satisfying $Pu=0,$ is bicharacteristic flow
inside $\{ p = 0, \, x = 0 \},$ a submanifold to which it is tangent. The
only critical points of the flow are at points $(y,\nu,0)$ where $y$ is a
critical point of $P$ and $\nu = \pm \sqrt{ \ev - V(y)}.$ Thus, the only
possible asymptotic escape directions of classical particles under the
influence of the potential $V$ are the finite number of critical points of
$V_0.$ Moreover, only the local minima are stable; the others have unstable
directions according to the number of unstable directions as a critical
point of $V_0:\partial X\longrightarrow\RR.$

The classical dynamics of $p$ and the quantum dynamics of $P$ are linked
via the scattering wavefront set. Let $u \in C^{-\infty}(X)$ be a tempered
distribution on $X$ (\ie in the dual space of $\CIdot(X;\Omega)).$ The part
of the scattering wavefront set, $\WFsc(u),$ of $u$ lying over the boundary
$ \{ x = 0 \},$ which is all that is of interest here, is a closed subset
of $\scT^*_{\partial X} X$ which measures the linear oscillations (Fourier
modes, in the case of Euclidean space) present in $u$ asymptotically near
boundary points; see \cite{Melrose43} for the precise definition. We shall
also need to use the scattering wavefront set $\WFsc^s(u)$ with respect to
the space $x^s L^2(X)$ which measures the microlocal regions where $u$
fails to be in $x^s L^2(X).$ There is a propagation theorem for the
scattering wavefront set in the style of the theorem of H\"ormander in the
standard setting; if $P u \in \CIdot(X),$ then the scattering wavefront set
of $u$ is contained in $\{ p = 0 \}$ and is invariant under the
bicharacteristic flow of $P,$ see \cite{Melrose43}. In particular,
generalized eigenfunctions of $u$ have scattering wavefront set invariant
under the bicharacteristic flow of $P.$ Note that the elliptic part
of this statement is already a uniform version of the smoothness of solutions.

In view of this propagation theorem, it is possible to consider where
generalized eigenfunctions `originate', although the direction of
propagation is fixed by convention. Let us say that a generalized
eigenfunction \emph{originates} at a radial point $q,$ if $q \in \WFsc(u)$
and if $\WFsc(u)$ is contained in the forward flowout $\Phi_+(q)$ of $q;$
thus each point in $\WFsc(u)$ can be reached from $q$ by travelling along
curves that are everywhere tangent to the flow and with $\nu$ nondecreasing
along the curve, so allowing the possibility of passing through radial
points, where the flow vanishes, on the way. In Part I of this paper we
showed, in the two-dimensional case and provided the eigenvalue $\ev$ is a
non-threshold value,

\begin{itemize}

\item Every $L^2$ eigenfunction is in $\CIdot(X).$ 

\item Every nontrivial generalized eigenfunction pairing to zero with the
$L^2$ eigenspace fails to be in $x^{-1/2} L^2(X).$

\item There are generalized eigenfunctions originating at each of the
incoming radial points in $\{ p = 0\},$ \ie at each critical point of $V_0$
with value less than $\ev.$

\item There are fundamental differences between the behaviour of
eigenfunctions near a local minimum and at other critical points. The
radial point corresponding to a local minimum is always an isolated point
of the scattering wavefront set for some non-trivial eigenfunction. For
other critical points, the scattering wavefront set necessarily propagates
and in generic situations each nontrivial generalized eigenfunction is
singular at some minimal radial point.

\item A generalized eigenfunction, $u,$ with an isolated point in its
scattering wavefront set, necessarily a radial point corresponding to a
local minimum of $V_0,$ has a complete asymptotic expansion there. The
expansion is determined by its leading term, which is a Schwartz function
of $n-1$ variables. The resulting map extends by continuity to an injective
map from $E^{\infty}_{\ess}(\ev)$ into $\oplus_q L^2(\RR^{n-1}),$ where the
direct sum is over local minima of $V_0$ with value less than the energy
$\ev.$

\item The space $E^0_{\ess}(\ev),$ consisting of those generalized
eigenfunctions which are in $x^{-1/2} L^2$ microlocally near $\{ \nu = 0
\},$ is a Hilbert space and the map above extends to a unitary isomorphism,
$M_+(\ev),$ from $E^0_{\ess}(\ev)$ to $\oplus_q L^2(\RR^{n-1}).$ A similar
map $M_-(\ev)$ can be defined by reversal of sign or complex conjugation
and the the scattering matrix for $P=P(\ev)$ at energy $\ev$ may be written
$$
S(\ev) =M_+(\ev) M^{-1}_-(\ev).
$$
\end{itemize}

In this paper we extend these results to higher dimensions.

\subsection{Results and structure of the paper}

We treat this problem by microlocal methods. Thus, the `classical' system,
consisting of the bicharacteristic vector field, plays a dominant role. The
main step involves reducing this vector field to an appropriate normal form
in a neighbourhood of each of its zeroes, which are just the radial
points. Nondegeneracy of the critical points of $V_0$ implies nondegeneracy
of the linearization of the bicharacteristic vector field at the
corresponding radial points. If there are no resonances, Sternberg's
Linearization Theorem, following an argument of Guillemin and Schaeffer,
allows the bicharacteristic vector field to be reduced to its linearization
by a contact transformation of $\scT^*_{\partial X} X.$ At the quantum
level this means that conjugation by a (scattering) Fourier integral
operator, associated to this contact transformation, microlocally replaces
$P$ by an operator with principal symbol in normal form. For this normal
form we construct `test modules' of pseudodifferential operators and
analyze the commutators with the transformed operator. Modulo lower order
terms, the operator itself becomes a quadratic combination of elements of
the test module. Just as in Part I, we use the resulting system of
regularity constraints to determine the microlocal structure of the
eigenfunctions and ultimately show the existence of asymptotic expansions
for eigenfunctions with some additional regularity.

However, the problem of resonances cannot be avoided. Even for a fixed
operator and fixed critical point, the closure of the set of values of
$\ev$ for which resonances occur may have non-empty interior. Such
resonances prevent the reduction of the bicharacteristic vector field to
its linearization, and hence of the symbol of $P$ to an associated model,
although partial reductions are still possible. In general it is necessary
to allow many more terms in the model. Fortunately most of these terms are
not relevant to the construction of the test modules and to the derivation
of the asymptotic expansions. We distinguish between `effectively
nonresonant' energies, where the additional resonant terms are such that
the definition of the test modules, now only to finite order, proceeds much
as before and the `effectively resonant' energies, where this is not the
case. Ultimately, we analyze the regularity of solutions at all
(non-threshold) energies. Near effectively nonresonant energies, smoothness of
families of eigenfunctions may still be readily shown. Effectively resonant
energies are harder to analyze, but the set of these is shown to be
\emph{discrete}. In any case, the space of microlocal eigenfunctions is
parameterized at all non-threshold energies. At effectively resonant
energies the problems arising from the failure of the direct analogue of
Sternberg's linearization are overcome by showing that, to an appropriate
finite order, the operator may be reduced to a non-quadratic function of
the test module.

In outline, the discussion proceeds as follows. In
sections~\ref{sec:radial} -- \ref{sec:micro-sol} we study radial
points. This is a general microlocal study except that we work under the
assumption that the symplectic map associated to the linearization of the
flow at each radial point (see Lemma~\ref{lemma:nondeg}) has no
4-dimensional irreducible invariant subspaces; this assumption is always
fulfilled in the case of our operator $\Lap + V-\ev.$ The main result is
Theorem~\ref{thm:model} in which the operator is microlocally conjugated to
a linear vector field plus certain `error terms'. In the nonresonant case
the error terms can be made to vanish identically, while in the effectively
nonresonant case the error terms have a good property with respect to a
test module of pseudodifferential operators, namely they can be expressed
as a positive power $x^\epsilon$, $\epsilon>0,$ times a power of the
module. In the effectively resonant case this is no longer possible and we
must allow `genuinely' resonant terms, but the set of effectively resonant
energies is discrete in the parameter $\ev$ in all dimensions.

We then turn in sections~\ref{sec:modules} -- \ref{sec:er} to studying
microlocal eigenfunctions which are microlocally outgoing at a given radial
point $q.$ The main result here is Theorem~\ref{thm:enr-structure} (or
Theorem~\ref{thm:er-structure} in the effectively resonant case) which
gives a parameterization of such microlocal eigenfunctions. For a minimal
radial point, they are parameterized by $\Schwartz(\RR^{n-1}),$ Schwartz
functions of $n-1$ variables, for a maximal radial point they are
parameterized by formal power series in $n-1$ variables, and in the
intermediate case of a saddle point with $k$ positive directions, they are
parameterized by formal power series in $n-1-k$ variables with values in
$\Schwartz(\RR^k).$ In all cases, the parameterizing data appear explicitly
in the asymptotic expansion of the eigenfunction at the critical point.

We next investigate in sections~\ref{sec:micro} and \ref{sec:mmd} the
manner in which the various radial points interact, and prove, in
Theorem~\ref{thm:mmd}, a `microlocal Morse decomposition.' This shows that
for each non-threshold energy $\ev$ there are genuine eigenfunctions (as
opposed to microlocal eigenfunctions) in $E^\infty_{\ess}(\ev)$ associated
to each energy-permissible critical point.

Then we turn in sections~\ref{sec:para} and \ref{sec:time} to the spectral
decomposition of $P$ and prove several versions of asymptotic
completeness. First this is established at a fixed, non-threshold energy;
see Theorem~\ref{thm:iso} which shows that the natural map from
$E^0_{\ess}(\ev)$ to the leading term in its asymptotic expansion (\ie to
its parameterizing data) is unitary. Next we prove a form valid uniformly
over an interval of the spectrum, Theorem~\ref{thm:AC}. In
section~\ref{sec:time} a time-dependent formulation is derived, as
Theorem~\ref{thm:time-AC}. This is based on the behaviour at large times of
solutions of the time-dependent Schr\"odinger equation $D_t u = Pu$ and is
subsequently used to derive a result of Herbst and Skibsted's on the
absence of $L^2$-channels corresponding to non-minimal critical points
(Corollary~\ref{cor:channels}).

\subsection{Results used from \cite{HMV1}} Throughout this paper we state
the specific location of results used from \cite{HMV1}. For the convenience
of the reader we summarize here the relevant locations. Sections~1-3 of
\cite{HMV1} are used as the basic background (and \cite[Section~3]{HMV1}
relies on Section~4 there).  The present Section~\ref{sec:micro-sol} is the
analogue of \cite[Section~5]{HMV1}, although we restate many of the
arguments due to the slightly different (more general) setting. The basic
analytic technique using test modules in Section~\ref{sec:modules} comes
from \cite[Section~6]{HMV1}. Certain results and methods from Sections~11
and 12 of \cite{HMV1} are used here in Sections~\ref{sec:mmd} and
\ref{sec:para}. However, the results of the intermediate sections 7-10 of
\cite{HMV1}, while certainly of interest when comparing to the results of
Section~\ref{sec:enr} and \ref{sec:er} here, are never used in the present
work directly or indirectly.

In addition, there was an error in the proof of Proposition~6.7 of
\cite{HMV1}.
While this error is minor and is easily remedied,
we present the modified proof,
together with some of the context, here in the Appendix since
this proposition lies at the heart of the analysis in both papers.

\subsection{Notation} 
The items listed below without a reference whose definition is not
immediate from the stated brief description are defined in \cite{Melrose43}.

\begin{tabbing}
\bf{Notation} \hskip 25pt \= \bf{Description/definition of notation} \hskip 60pt \= \bf{Reference} \\
$V_0$ \> restriction of $V$ to $\partial X$ \> \\
$\Cv(V)$ \> set of critical values of $V_0$ \> \\
$\scT^* X$ \> scattering cotangent bundle over $X$ \> \eqref{numu} \\
$\scT^*_{\partial X} X$ \> restriction of $\scT^* X$ to $\partial X$ \> \eqref{numu} \\
$x$ \> boundary defining function of $X$ s.t. \eqref{sc-metric} holds \> \\
$y$ \> coordinates on $\partial X$ \> \\
$(\nu, \mu)$ \> fibre coordinates on $\scT^* X$ \> \eqref{numu} \\
$y = (y', y'', y''')$ \> decomposition of $y$ variable \> \eqref{y-decomp} \\
$\mu = (\mu', \mu'', \mu''')$ \> dual decomposition of $\mu$ variable \> \eqref{y-decomp} \\
$r'_i, r''_j, r'''_k$ \> eigenvalues of the contact map $A$ \> \eqref{y-decomp} \\ 
$Y''_j$ \> $y''_j / x^{r''_j}$ \> \eqref{eq:def-Phi} \\
$Y'''_k$ \> $y'''_k / x^{1/2}$ \> \eqref{eq:def-Phi} \\
$\Delta$ \> (positive) Laplacian with respect to $g$ \> \\
$P$ \> $x^{-1} (\Delta + V - \ev)$ \> Sec. ~\ref{sec:radial} \\
$H$ \> $\Delta + V$ \> \\
$R(\ev)$ \> resolvent of $H$,  $(H - \ev)^{-1}$ \> \\
$R(\ev \pm i0)$ \> limit of resolvent on real axis from above/below \> \\
$\tilde V$ \> modified potential \> Lem. \ref{lemma:Vt} \\
$\Sp(\ev)$ \> (generalized) spectral projection of $H$ at energy $\ev$ \>  \eqref{spsp} \\
$\tilde R(\ev)$ \> resolvent of modified potential $(\Delta + \tilde V - \ev)^{-1}$ \> \\
$L^2_{\text{sc}}(X)$ \> $L^2$ space with respect to Riemannian density of $g$ \> \\ 
$\Hsc^{m,0}(X)$ \> Sobolev space; image of $L^2_{\text{sc}}(X)$ under $(1 + \Lap)^{-m/2}$ \> \\
$\Hsc^{m,l}(X)$ \> $x^l \Hsc^{m,0}(X)$ \> \\
$\Psisc^{m,0}(X)$ \> scattering pseudodiff. ops. of differential order $m$ \> \\
$\Psisc^{m,l}(X)$ \> $x^l \Psisc^{m,0}(X)$; maps $\Hsc^{m',l'}(X)$  to $\Hsc^{m' - m, l'+l}(X)$ \> \\
$\sigma_{\partial, l}(A)$ \> boundary symbol of $A \in \Psisc^{m,l}(X)$; $\CI$ fn. on $\scT^*_{\partial X} X$ \> \\
$\sigma_{\partial}(A)$ \> $\sigma_{\partial, 0}(A)$ \> \\
$\WFsc(u)$ \> scattering wavefront set of $u$; closed subset of $\scT^*_{\partial X} X$ \> \\
$\WFsc^{m,l}(u)$ \> scattering wavefront set with respect to $\Hsc^{m,l}$ \> \\
$\WFsc'(A)$ \>operator scattering wave front set; in its complement \> \\
\>$A$ is microlocally in $\Psisc^{*,\infty}(X),$ i.e.\ is trivial \> \\
$\scH_p$ \> scattering Hamilton vector field \> Sec. \ref{sec:radial} \\
$\Phi_+(q)$ \> forward flowout from $q \in \scT^*_{\partial X} X$ \> Sec. 1.1  \\
radial point  \> point in $\scT^*_{\partial X} X$ where $p$ and $\scH_p$  vanish \> Sec. \ref{sec:radial} \\
$\RP_{\pm}(\ev)$ \> set of radial points of $H - \ev$ where $\pm \nu > 0$ \> \\
$\Min_+(\ev)$ \> subset of $\RP_+(\ev)$ associated to local minima of $V_0$ \> \\
$\leq$ \> partial order on $\RP_+(\ev)$ compatible with $\Phi_+$ \> Def. \ref{Def:partial-order} \\
$\tilde E_{\mic, +}(O,P)$ \> microlocal solutions of $Pu = 0$ in the set $O$ \> \eqref{HMV.42} \\
$E_{\mic, +}(q, \ev)$ \> microlocal solutions of $(H - \ev) u = 0$ near q \> \eqref{Emic-defn}\\
$E^s_{\ess}(\ev)$ \> space of generalized $\ev$-eigenfunctions of $H$ \> \eqref{Edefn} \\
$E^s(\Gamma, \ev)$ \> subset of $u \in E^s_{\ess}(\ev)$ with $\WFsc(u) \cap \RP_+(\ev) \subset \Gamma$ \> \eqref{EGamma-defn} \\
$E^s_{\Min, +}(\ev)$ \> $E^s(\Gamma, \ev)$, with $\Gamma = \Min_+(\ev)$ \> \\
$\mathcal{M}$ \> test module \> Sec. \ref{sec:modules} \\
$\Isc^{(s)}(O,\cM)$ \> space of iteratively-regular functions w.r.t.  $\cM$ \> \eqref{Isc-defn} \\
$\tau$ \> rescaled time variable; $\tau = xt$ \>  Sec. \ref{sec:time} \\
$\Xsch$ \> $X \times \Rbar_{\tau}$ \> \eqref{Xsch} \\
\end{tabbing}

\emph{Acknowledgement.} We thank an anonymous referee for helpful comments
and for pointing out some errors in the first version of this manuscript. 


\section{Radial points}\label{sec:radial}

Let $X$ be a compact $n$-dimensional manifold with smooth boundary. Recall
that if $(x,y)$ are local coordinates on $X,$ with $x$ a boundary defining
function, then dual scattering coordinates $(\nu,\mu)$ on the scattering
cotangent bundle are determined. The restriction of the scattering
cotangent bundle to $\partial X$ is denoted $\scT^*_{\pa X}X$ and has a
natural contact structure, the contact form at the boundary being
\begin{equation}
\alpha = -d\nu+ \sum_i \mu_i dy_i
\label{contact-form}\end{equation}
in local coordinates. Recall that a contact structure on a
$2n-1$-dimensional manifold, here $\scT^*_{\pa X}X$, is given by a
nondegenerate one-form, i.e.\ a one-form $\alpha$ with
$\alpha\wedge(d\alpha)^{n-1}$ everywhere non-zero; correspondingly its
kernel is a maximally non-integrable hyperplane field on $\scT^*_{\pa
X}X$. One refers to either the line bundle given by the span of $\alpha$,
or the hyperplane field given by its kernel, as the contact structure.

Suppose that $P\in\Psisc^{*,-1}(X)$ is a scattering pseudodifferential
operator of order $-1$ at the boundary; for example, $P = x^{-1} (\Delta +
V - \ev).$ Then the boundary part of its principal symbol,
$p=\sigma_{\pa}(xP),$ is a $\Cinf$ function on $\scT^*_{\pa X}X.$ In this,
and the next, section we consider radial points of a general real-valued
function, $p\in\CI(\scT^*_{\pa X}X),$ with only occasional references to
the particular case, $p=|\zeta|^2+V_0-\ev,$ of direct interest in this
paper.  Although we discuss radial points in the context of boundary points
in the scattering calculus this analysis applies directly (and could
alternatively be done for) radial points in the usual microlocal picture,
as described in the Introduction.  Our objective in this section is to find
a change of coordinates, preserving the contact structure, in which the
form of $p$ is simplified. In this section we consider the simplification
of $p$ up to second order, in a sense made precise below.

The basic non-degeneracy assumption we make is that 
\begin{equation}
p=0\text{ implies }dp\neq 0;
\label{HMV.177}\end{equation}
this excludes true `thresholds' which however do occur for our problem,
when $\ev$ is a critical value of $V_0.$ It follows directly from
\eqref{HMV.177} that the boundary part of the characteristic variety
\begin{equation*}
\Sigma=\{q\in\scT^*_{\pa X}X;p(q)=0\}\text{ is smooth;}
\end{equation*}
we shall assume that $\Sigma$ is compact, corresponding to the ellipticity
of $P.$

\begin{Def} A \emph{radial point} for a function $p$ satisfying
\eqref{HMV.177} is a point $q \in \Sigma$ such that $dp(q)$ is a
(necessarily nonzero) multiple of the contact form $\alpha$ given by
\eqref{contact-form}. Conversely, if $q \in \Sigma$ and $dp$ and $\alpha$
are linearly independent at $q$ then we say that $p$ is of \emph{principal
type} at $q.$
\end{Def}

We may extend $p$ to a $\Cinf$ function on $\scT^*X$, still denoted by $p.$
Over the interior $\scT^*_{X^\circ}X$ is naturally identified with $T^*X^\circ,$
which is a symplectic manifold with canonical symplectic form $\omega.$ Near
the boundary, expressed in terms of scattering-dual coordinates, 
\begin{equation}
\omega=d\left(-\nu\,\frac{dx}{x^2}+\sum_i\mu_i\,\frac{dy_i}{x}\right)
=(-d\nu+\sum_i\mu_i\,dy_i)\wedge\frac{dx}{x^2}+
\sum_id\mu_i\wedge\frac{dy_i}{x}.
\label{HMV.178}\end{equation}
Consider the Hamilton vector field, $H_{x^{-1}p},$ of $x^{-1}p$, which we
shall denote $\scH_p$, fixed by the identity $\omega(\cdot,\scH_p)=dp.$
Then $\scH_p$ extends to a vector field on $\scT^*X$ tangent to its
boundary, so $\scH_p\in\Vb(\scT^*X).$\footnote{Here $\Vb(M)$ denotes the
space of smooth vector fields on the manifold with boundary $M$ that are
tangent to $\partial M$.} At the boundary $\scH_p$, as an element of
$\Vb(\scT^*X)$, is independent of the extension of $p.$ We denote the
restriction of $\scH_p$ (as a vector field) to $\scT^*_{\pa X}X$ by $W,$ so
$W$ is a vector field on $\scT^*_{\pa X}X.$ Explicitly in local
coordinates
\begin{equation}\begin{split}\label{eq:scHp-local}
\scH_p=&-(\pa_\nu p) (x\pa_x+\mu\cdot\pa_\mu)+(x\pa_x p - p +\mu
\cdot\pa_\mu p)\pa_\nu\\
&+\sum_j \left(\pa_{\mu_j} p\,\pa_{y_j}-\pa_{y_j} p\,\pa_{\mu_j}\right)
+x\Vb(\scT^* X);
\end{split}\end{equation}
since $p$ is smooth up to the boundary, $x\pa_x p=0$ at $\scT^*_{\pa X}X.$ Thus,
\begin{equation}
W=-(\pa_\nu p) \mu\cdot\pa_\mu+(\mu\cdot\pa_\mu p - p)\pa_\nu
+\sum_j \left(\pa_{\mu_j} p\,\pa_{y_j}-\pa_{y_j} p\,\pa_{\mu_j}\right).
\label{W-coords}\end{equation}

Alternatively $W$ may be described in terms of the contact structure on
$\scT^*_{\pa X}X$. Namely $W$ is the Legendre vector field of $p,$ determined by
\begin{equation} 
d\alpha(.,W)+\gamma\alpha=dp,\ \alpha(W)=p 
\label{HMV.179}\end{equation}
for some function $\gamma.$ It follows that $W$ is tangent to $\Sigma,$
since $dp(W)=\gamma \alpha(W) = \gamma p=0$ at any point at which $p$
vanishes. An equivalent definition of $q \in \Sigma$ being a radial point
is that the vector field $W$ vanishes as $q$, as follows from
\eqref{HMV.179} and the nondegeneracy of $\alpha$.
  
\begin{Def} A radial point $\Rp\in\Sigma$ for a real-valued function
$p\in\CI(\scT^*_{\pa X}X)$ satisfying \eqref{HMV.177} is said to be
\emph{non-degenerate} if the vector field $W$, restricted to $\Sigma =\{p=
0\},$ has a non-degenerate zero at $\Rp.$ Note that this implies that a
non-degenerate radial point is necessarily isolated in the set of radial points.
\end{Def}

Since the vector field $W$ vanishes at a radial point $\Rp,$  its
linearization is well defined as a linear map, $A'$ on $T_{\Rp}\scT^*_{\pa
X}X,$ (later we will use the transpose, $A,$ as a map on differentials)
\begin{equation} A'v=[V,W](\Rp),
\label{HMV.181}\end{equation}
for any smooth vector field $V$ with $V(\Rp)=v;$ it is independent of the
choice of extension and can also be written in terms of the Lie derivative
\begin{equation}
A'v=-\mathcal{L}_WV(\Rp).
\label{HMV.182}\end{equation}
Since $Wp=\gamma p,$ $A'$ preserves the subspace $T_{\Rp}\Sigma.$ Since
$\alpha$ is normal to $T_{\Rp}\Sigma$, the restriction of $d\alpha$ to $T_{\Rp}\Sigma$
is a symplectic 2-form, $\omega_{\Rp}.$

\begin{lemma}\label{HMV2r.192} At a non-degenerate radial point for $p,$
  where $dp=\lambda \alpha,$ the linearization $A'$ acting on $T_{\Rp}\Sigma$ is such that 
\begin{equation*}
S \equiv A'-\frac12\lambda \Id \in \symp(2(n-1))
\label{HMV2r.193}\end{equation*}
is in the Lie algebra of the symplectic group with respect to $\omega _{\Rp}:$ 
\begin{equation*}
\omega _{\Rp}(Sv_1,v_2)+\omega _{\Rp}(v_1,Sv_2)=0,\ \forall\ v_1,v_2\in
T_{\Rp}\Sigma.
\label{HMV.186}\end{equation*}
\end{lemma}

\begin{proof} Observe that \eqref{HMV.179} implies that 
\begin{equation}
L_W\alpha =(d\alpha )(W,\cdot)+d(\alpha (W))=\gamma \alpha.
\label{HMV.183}\end{equation}
For two vector smooth vector fields $V_i,$ defined near $\Rp,$
\begin{equation}
\begin{gathered}
W(d\alpha (V_1,V_2))=L_W(d\alpha (V_1,V_2))\\
=(L_Wd\alpha )(V_1,V_2)+d\alpha
(L_WV_1,V_2)+d\alpha(V_1,L_WV_2).
\end{gathered}
\label{HMV.184}\end{equation}
The left side vanishes at $\Rp$ so using \eqref{HMV.182} 
\begin{equation}
\omega _{\Rp}(A'v_1,v_2)+\omega _{\Rp}(v_1,A'v_2)=\lambda \omega
_{\Rp}(v_1,v_2)\ \forall\ v_1,v_2\in T_{\Rp}\Sigma .
\label{HMV.185}\end{equation}
\end{proof}

It follows from Lemma~\ref{HMV2r.192}, see for example \cite{MR55:3504},
that $A'$ is decomposable into invariant subspaces of dimension $2$ and
$4$, with eigenvalues on the two-dimensional subspaces of the form $\lambda
r,$ $\lambda( 1 - r),$ $r \leq 1/2 $ real or $\lambda (1/2 + i a),$
$\lambda (1/2 - i a),$ with $a >0.$ 

Note that, by \eqref{HMV.179}, $d_\nu p(q) = -\gamma(q) = - \lambda$, so
from \eqref{eq:scHp-local}, the Hamilton vector field $\scH_p$ is equal to
$\lambda x \partial_x$ modulo vector fields of the form $f \cdot W'$ where
$W$ is tangent to $\{ x = 0 \}$ and $f(q) = 0.$ 
Therefore if $\lambda>0,$ then $x$ is increasing along
bicharacteristics of $p$ in the interior of $\scT^* X,$ \ie the
bicharacteristics leave the boundary, \ie `come in from infinity' if
$\pa X$ is removed, while if $\lambda<0$, the
bicharacteristics approach the boundary, \ie `go out to infinity'.
Correspondingly we make the following definition.

\begin{Def}
We say that a non-degenerate radial point $\Rp$ for $p$ with
$dp(\Rp)=\lambda\alpha(\Rp)$ is outgoing
if $\lambda<0$, and we say that it is incoming if $\lambda>0$.
\end{Def}

For $p = |\zeta|^2+V_0-\ev$, we have $\lambda = -\partial_\nu p =
-2\nu$. Hence, radial points are outgoing for $\nu > 0$ and incoming for
$\nu < 0$ in this case.  We next discuss the form the linearization takes
for $p=|\zeta|^2+V_0-\ev$.

\begin{lemma}\label{lemma:nondeg} For the function $p=|\zeta|^2+V_0-\ev$ with
  $V_0$ Morse, the radial points are all nondegenerate and the linear
  operator $S$ associated with each has only two-dimensional
  invariant symplectic subspaces.
\end{lemma}

\begin{rem} In view of the non-occurrence of non-decomposable invariant
  subspaces of dimension $4$ in this case we will exclude them from further
  discussion below.
\end{rem}

\begin{proof} Choose Riemannian normal coordinates $y_j$ on $\pa X,$ so the
metric function $h$ satisfies $h-|\mu|^2=\cO(|y|^2).$ Since the Hessian of
$V|_{\pa X}$ at a critical point is a symmetric matrix, it can be
diagonalized by a linear change of coordinates on $\pa X,$ given by a
matrix in $\mathrm{SO}(n-1)$, which thus preserves the form of the
metric. It follows that for each $j,$ $(dy_j,d\mu_j)$ is an invariant
subspace of $A.$
\end{proof}

Let $\cI$ denote the ideal of $\Cinf$ functions on $\scT^*_{\pa X}X$
vanishing at a given radial point, $\Rp.$ The linearization of $W$ then
acts on $T^*_{\Rp}\left(\scT^*_{\pa X} X\right)=\cI/\cI^2;$ $dp(\Rp),$ or 
equivalently $\alpha_{\Rp},$ is necessarily an eigenvector of $A$ with
eigenvalue $0.$ Similarly, $\scH_p$ defines a linear map $\tilde A$ on
$T^*_{\Rp}\left(\scT^* X\right).$ By \eqref{eq:scHp-local},  $\tilde A$
preserves the conormal line, $\Span{dx}$ and the eigenvalue of $\tilde A$
corresponding to the eigenvector $dx$ is $\lambda.$ Thus $\tilde A$ acts on
the quotient
\begin{equation*}
T^*_{\Rp}\left(\scT^*_{\pa X} X\right)\equiv T^*_{\Rp}\left(\scT^*
X\right)/\Span{dx},
\label{HMV.180a}\end{equation*}
and this action clearly reduces to $A.$ 

By Darboux's theorem we may make a local contact diffeomorphism of
$\scT^*_{\pa X}X$ and arrange that $\Rp=(0,0,0).$ Thus, as a module over
$\Cinf(\scT^*_{\pa X}X)$ in terms of multiplication of functions, $\cI$ is
generated by $\nu,$ $y_j$ and the $\mu_j,$ for $j=1,\ldots,n-1.$ Thus in
general we have the following possibilities for the two-dimensional
invariant subspaces of $A.$
\begin{enumerate}
\item There are two independent real eigenvectors with eigenvalues in
  $\lambda(\bbR\setminus[0,1]).$ 
\item There are two independent real eigenvectors with eigenvalues in
  $\lambda(0,1).$
\item There are no real eigenvectors and two complex eigenvectors with
  eigenvalues in $\lambda(\frac12+i(\bbR\setminus\{0\})).$
\item \label{Hessthres}There is only one non-zero real eigenvector with
  eigenvalue $\frac12\lambda.$\label{Hessian-threshold}
\end{enumerate}
Case \eqref{Hessthres} was called the `Hessian threshold' case in Part I. 
In all cases the sum of the two (generalized) eigenvalues is $\lambda.$ 

\begin{lemma}\label{HMV2r.191} 
By making a change of contact coordinates, i.e.\ a change of coordinates
on $\scT^*_{\pa X}X$
preserving the contact structure, near a radial point $\Rp$ for
$p\in\CI(\scT^*_{\pa X}X)$ for which the linearization has neither a
Hessian threshold subspace, \eqref{Hessian-threshold}, nor any
non-decomposable 4-dimensional invariant subspace, coordinates $y$ and
$\mu,$ decomposed as $y = (y', y'',y''')$ and $\mu =
(\mu', \mu'', \mu'''),$ may be introduced so that
\begin{enumerate}
\item 
\begin{equation}\label{y-decomp}(y',\mu')=(y_1,\dots,y_{s-1},\mu
  _1,\dots,\mu_{s-1})\end{equation}
where $e'_j=dy'_j$, $f'_j=d\mu'_j$ are eigenvectors of $A$ with eigenvalues
$\lambda r'_j,$ $\lambda(1-r'_j),$ $j=1,\ldots,s-1$ with $r'_j<0$ real and
negative.
\item $(y'',\mu'')=(y_s,\dots,y_{m-1},\mu _{s},\dots,\mu_{m-1})$ 
where $e''_j=dy''_j$, $f''_j=d\mu''_j$ are eigenvectors with eigenvalues $\lambda
r''_j,$ $\lambda(1-r''_j),$ $j=s,\ldots,m-1$ where $0 < r''_j \leq 1/2$ is
real and positive.
\item $(y''',\mu''')=(y_{m},\dots,y_{n-1},\mu _{m},\dots,\mu_{n-1}),$
where some complex combination $e'''_j,$ $f'''_j,$ of $dy'''_j$ and
$d\mu'''_j,$ $m \leq j \leq n-1,$ are eigenvectors with eigenvalues
$\lambda r'''_j$ and $\lambda (1 - r'''_j)$ with $r'''_j = 1/2 +
i\beta'''_j$, $\beta'''_j > 0.$ 
\end{enumerate} 
\end{lemma}
\noindent
Thus if we set $e = (e', e'', e''')$, $f = (f', f'',f''')$ the
eigenvectors of $A$ are $d\nu, e_j$ and $f_j,$ with respective eigenvalues
$0, \lambda r_j$ and $\lambda(1 - r_j);$ we will take the coordinates so
that the $r_j$ are ordered by their real parts.

\begin{rem}
We emphasize that the change of coordinates here is on the contact space,
$\scT^*_{\pa X}X$, and it is, in general, not induced by a change of
coordinates on $X$. Analytically it is implemented by a scattering FIO (see Section~\ref{lft}). 
\end{rem}

In coordinates in which the eigenspaces take this form it can be seen
directly that 
\begin{equation}
p=\lambda\big(-\nu+\sum_{j=1}^{m-1} r_j y_j\mu_j+\sum_{j=m}^{n-1}Q_j(y_j,\mu_j)+\nu g_1+g_2\big)
\label{HMV2r.189}\end{equation}
with the $Q_j$ elliptic homogeneous polynomials of degree 2, $g_1$ vanishing
at least linearly and $g_2$ to third order. 

\begin{rem}\label{rem:H-lin-ev}
For the function $p=|\zeta|^2+V_0-\ev$ with $V_0$ Morse, the eigenvalues of
$A$ at a radial point $\Rp$ are easily calculated in the coordinates
used in the proof of Lemma~\ref{lemma:nondeg}. Indeed, since the 2-dimensional
invariant subspaces decouple,
the results of \cite[Proof of Proposition~1.2]{HMV1}
can be used. The eigenvalues corresponding
to the 2-dimensional subspace in which the eigenvalue of the Hessian is
$2a_j$ are thus
\begin{equation*}
\lambda\left(\frac{1}{2}\pm\sqrt{\frac{1}{4}-\frac{a_j}{\ev-V_0(0)}}\right),
\Mwhere \lambda=-2\nu(\Rp).
\end{equation*}
\end{rem}

In fact, below we do not need the full power of Lemma~\ref{HMV2r.191}.
Essentially it suffices if we arrange that the eigenvectors corresponding
to the (in absolute value) larger eigenvalues, namely
$\lambda(1-r'_j)$, if $r'_j<0$, or $\lambda(1-r''_j)$,
if $r''_j\in(0,\half)$, are
in a model form on the two dimensional eigenspaces. The advantage of
the weaker conclusion is that one has more freedom in choosing the
contact change of coordinates.

\begin{lemma}\label{HMV2r.191p} (Weaker version of Lemma~\ref{HMV2r.191}.)
Suppose that $\half\lambda$ is not an eigenvalue of $A$.
By making a change of contact coordinates, i.e.\ a change of coordinates
on $\scT^*_{\pa X}X$
preserving the contact structure, near a radial point $\Rp$ for
$p\in\CI(\scT^*_{\pa X}X)$ for which the linearization has neither a
Hessian threshold subspace, \eqref{Hessian-threshold}, nor any
non-decomposable 4-dimensional invariant subspace, coordinates $y$ and
$\mu,$ decomposed as $y = (y', y'',y''')$ and $\mu =
(\mu', \mu'', \mu'''),$ may be introduced so that
\begin{enumerate}
\item 
\begin{equation}\label{y-decomp-p}(y',\mu')=(y_1,\dots,y_{s-1},\mu
  _1,\dots,\mu_{s-1}),\end{equation}
where some real linear combinations $e'_j$ of $d\mu'_j$ and $dy'_j$,
resp.\ $f'_j=d\mu'_j$ are eigenvectors of $A$ with eigenvalues
$\lambda r'_j,$ resp.\ $\lambda(1-r'_j),$
$j=1,\ldots,s-1$ with $r'_j<0$ real and
negative.
\item $(y'',\mu'')=(y_s,\dots,y_{m-1},\mu _{s},\dots,\mu_{m-1})$ 
where some real linear combinations $e''_j$ of $d\mu''_j$ and $dy''_j$,
resp.\ $f''_j=d\mu''_j$ are eigenvectors with eigenvalues $\lambda
r''_j,$ $\lambda(1-r''_j),$ $j=s,\ldots,m-1$ where $0 < r''_j < 1/2$ is
real and positive.
\item $(y''',\mu''')=(y_{m},\dots,y_{n-1},\mu _{m},\dots,\mu_{n-1}),$
where some complex combination $e'''_j,$ $f'''_j,$ of $dy'''_j$ and
$d\mu'''_j,$ $m \leq j \leq n-1,$ are eigenvectors with eigenvalues
$\lambda r'''_j$ and $\lambda (1 - r'''_j)$ with $r'''_j = 1/2 +
i\beta'''_j$, $\beta'''_j > 0.$ 
\end{enumerate} 
\end{lemma}
\noindent
Again, if we set $e = (e', e'', e''')$, $f = (f', f'',f''')$ the
eigenvectors of $A$ are $d\nu, e_j$ and $f_j,$ with respective eigenvalues
$0, \lambda r_j$ and $\lambda(1 - r_j);$ we will take the coordinates so
that the $r_j$ are ordered by their real parts. In these coordinates a
version of \eqref{HMV2r.189} still holds, namely if $a_j$ and $b_j$ are
any functions on $\scT^*_{\pa X}X$ vanishing at $(0,0,0)$ with
differential $e_j$, resp. $f_j$, $j=1,\ldots,m-1$ (so we may take
$b_j=\mu_j$, and we may take $a_j$ a $\Real$-linear combination of $y_j$ and
$\mu_j$) then
\begin{equation}\begin{split}
p&
=\lambda\big(-\nu+\sum_{j=1}^{m-1} r_j a_j b_j+\sum_{j=m}^{n-1}Q_j(y_j,\mu_j)+\nu g_1+g_2\big)\\
&=\lambda\big(-\nu+\sum_{j=1}^{m-1} r_j y_j\mu_j+\sum_{j=1}^{m-1}c_j\mu_j^2
+\sum_{j=m}^{n-1}Q_j(y_j,\mu_j)+\nu g_1+g_2\big),
\end{split}\label{HMV2r.189p}\end{equation}
where the $c_j$ are real, the $Q_j$ are elliptic homogeneous polynomials of
degree 2, $g_1$ vanishes at least linearly and $g_2$ to third order. 

As mentioned, Lemma~\ref{HMV2r.191p} is weaker than, hence is an immediate
consequence of, Lemma~\ref{HMV2r.191}. Although it is by no means
essential, this weaker result leaves more freedom in choosing the contact
map which is useful in making the choice rather explicit, if this is
desired. In fact, if $p=|\zeta|^2+V_0-\sigma,$ as in
Lemma~\ref{lemma:nondeg}, we immediately deduce the following.

\begin{lemma}\label{lemma:special-contact}
For the function $p=|\zeta|^2+V_0-\ev$ with
  $V_0$ Morse, the contact map in Lemma~\ref{HMV2r.191p} can be taken
as the composition of the contact map on $\scT^*_{\pa X}X$ induced
by a change of coordinates on $X$, with the canonical relation
of multiplication
by a function of the form $e^{i\phi/x}$, $\phi\in\Cinf(X)$.
\end{lemma}

\begin{rem}
The canonical relation of multiplication by $e^{i\phi/x}$ is given, in local
coordinates $(y,\nu,\mu)$, by the map
\begin{equation*}
\Phi_\phi:(y,\nu,\mu)\mapsto (y,\nu+\phi(y),\mu+\pa_y \phi(y)),
\end{equation*}
i.e.\ if we write $\Phi_\phi(y,\nu,\mu)=(\bar y,\bar \nu,\bar \mu)$, then
$\bar \mu_k=\mu_k+\pa_{y_k}\phi(y)$. Note that while $\phi$ is a function on
$X$, the canonical relation only depends on $\phi|_{\pa X}$, which is
why we simply regard $\phi$ as a function on $\pa X$ and write $\phi(y)$
here.
\end{rem}

\begin{proof}
As in the proof of Lemma~\ref{lemma:nondeg} we may assume, by a change of
coordinates on $X$, that the critical point of $V_0$ over which
the radial point $\Rp$ lies is $y=0$, that $h-|\mu|^2=\cO(|y|^2)$
and that the Hessian of $V_0$ at $0$
is diagonal, so for each $j$, $(dy_j,d\mu_j)$
is an invariant subspace of $A$. Note that in the coordinates $(y,\nu,\mu)$,
$\Rp=(0,\nu_0,0)$. With the notation of Remark~\ref{rem:H-lin-ev} above,
if $dy_j$ is an eigenvector
of the Hessian with eigenvalue $2a_j$ then the eigenvectors
of $A$ of eigenvalue $\lambda r_j$, resp.\ $\lambda(1-r_j)$, are
$\tilde e_j=\frac{\lambda}{2} (1-r_j)dy_j+d\mu_j$,
resp.\ $\tilde f_j=\frac{\lambda}{2} r_jdy_j+d\mu_j$, see
Remark~1.3 of \cite{HMV1}. In particular, if $r_j$ is real, so is
$\tilde f_j$.

Now, the contact map $\Phi_\phi$ induced by multiplication by $e^{i\phi/x}$
as above acts on $T^*\scT^*_{\pa X}X$ by pull-backs, namely
\begin{equation*}\begin{split}
\Phi_\phi^*&(\sum_k \bar y_k^*\, d\bar y_k+\bar\nu^*\,d\bar\nu
+\sum_k \bar\mu_k^*\,d\bar\mu_k)\\
&=\sum_k \bar y_k^*\, dy_k+
\bar\nu^* (d\nu+\sum_j(\pa y_j \phi)\,dy_j)+\sum_k \bar\mu_k^*
(d\mu_k+\sum_j\pa_{y_j}\pa_{y_k}\phi(y)\,dy_j).
\end{split}\end{equation*}
Thus, by the above remark, $\Phi$ will map $\Rp$ to $(0,0,0)$ provided
$\phi(0)=-\nu_0$, $\pa_{y_j}\phi(0)=0$ for all $j$. In this case,
moreover, the pull-back $\Phi_\phi^*$
will map $dy_k$ to $dy_k$,
$d\nu$ to $d\nu$ and $d\mu_k$ to $d\mu_k+\sum_j\pa_{y_j}\pa_{y_k}\phi(y)$.
Correspondingly, by letting $\phi(y)=-\nu_0+\sum_{j=1}^{m-1} b_j y_j^2$,
$b_j=\frac{\lambda}{4} r_j$,
$(\Phi_\phi^{-1})^*$ maps $\tilde f_j$ to $d\mu_j$, $j=1,\ldots,m-1$.
Since the Legendre vector field $W'$
of $(\Phi_\phi^{-1})^*p$ is the push-forward
of the Legendre vector field $W$ of $p$ under $\Phi_\phi$, it follows that
$d\mu_j$ is an eigenvector of the linearization of $W'$ with eigenvalue
$\lambda(1-r_j)$.
As $\Phi_\phi^*$ also maps the 2-dimensional subspaces
$(dy_j,d\mu_j)$ (at $(0,0,0)$) to the 2-dimensional subspaces $(dy_j,d\mu_j)$
(at $\Rp$), and the latter are invariant under $A$, so are the former under
the linearization of $W'$.
This proves the lemma.
\end{proof}


\section{Microlocal normal form}\label{sec:normal} Let $P \in \Psisc^{*,
-1}(X)$ be an operator with real principal symbol $p$ obeying
\eqref{HMV.177}, as in the previous section, and assume that $q$ is a
nondegenerate radial point for $p$. In this section we shall reduce $p$ to
a normal form, via conjugation with a scattering Fourier integral
operator. We first pause to define such operators.

\subsection{Scattering Fourier integral operators}\label{lft} Scattering
Fourier integral operators (FIOs) are defined in terms of conventional FIOs
via the local Fourier transform, as defined in \cite{MZ}. Let $X$ be a
manifold of dimension $n$ with boundary, and $(x,y)$ local coordinates
where $x$ is a boundary defining function. We can always identify a
neighbourhood $U \subset \partial X$ of $y_0 \in \partial X$ with an open
set $V \in S^{n-1}$, which we can think of as embedded in $\RR^n$ in the
standard way. Correspondingly we may identify the interior of a
neighbourhood $[0, \epsilon)_x \times U \subset X$ of $(0, y_0) \in X$ with
the an asymptotically conic open set $(\epsilon^{-1}, \infty) \times V
\subset \RR^n$ in $\RR^n$. If we choose a function $\phi \in C^\infty(X)$
supported in $[0, \epsilon)_x \times U$ which is identically $1$ in a
neighbourhood of $(0, y_0)$, then the operator $\Four$ with kernel
$$
e^{i z \cdot y/x} \phi(x,y) \frac{d\omega(y) dx}{x^{n+1}}
$$
is called a `local Fourier transform' on $X$. Here $z = (z_1, \dots, z_n)
\in \RR^n$, $z \cdot y$ denotes the inner product on $\RR^n$ and
$d\omega(y)$ denotes the standard measure on $S^{n-1}$ (pulled back to
$\partial X$ and then to $X$ via the identifications above). Of course, if
$X$ is the radial compactification of $\RR^n$ and the identification
between $U$ and $V$ is the identity, then $\Four$ really \emph{is} the
Fourier transform premultiplied by the cutoff function $\phi$.

It is shown in \cite{MZ} that $\Four$ induces a local bijection between
$\scT^*_{\pa X} X$ and the cosphere bundle of $\RR^n$. In fact, using our
identification between $U$ and $V \subset S^{n-1}$ we may represent points
in $\scT^*_{U} X$ as $(\hat z, \zeta)$ where $\hat z = z/|z| \in V$
represents a point in $U$ and $\zeta$ represents the point in the fibre
given by $(\nu, \mu)$ where $\nu$ is the parallel component of $\zeta$
relative to $\hat z$ and $\mu$ is the orthogonal component. The
identification is then given by the Legendre map
$$
L(\hat z, \zeta) = (\zeta, -\hat z) \in S^* \RR^n.
$$

In other words, $\Four$ sets up a bijection between scattering wavefront
set and conventional wavefront set. Moreover, it is shown in \cite{MZ} that
conjugation by $\Four$ maps the scattering pseudodifferential operators $A \in
\Psisc^{*,l}(X)$ microsupported near $(y_0, \nu_0, \mu_0)$ to the conventional
pseudodifferential operators microsupported near $L(y_0, \nu_0, \mu_0),$
with principal symbols related by
$$
\sigma^l(\Four A \Four^*)(L(q)) = a(q),
$$
where $a$ is the boundary symbol of $A$ (of order $l$). 

\begin{Def} A scattering FIO is an operator $E$ from $\CIdot(X)$ to
$\Dist(X)$ such that, for any local Fourier transforms $\Four_1$, $\Four_2$
on $X$, $\Four_2 E \Four_1^*$ is a conventional FIO on $\RR^n$.
\end{Def}

A simple example of a scattering FIO is multiplication by an oscillatory
factor $e^{i\psi(y)/x}$. Under conjugation by a local Fourier transform
this becomes a conventional FIO given by an oscillatory integral with phase
function $(z-z')\cdot \zeta + |\zeta|\psi(\zeta/|\zeta|)$. The scattering
resolvent kernel constructed in \cite{HV1} and \cite{HV2}, microlocalized
to the interior of the `propagating Legendrian', is another example.

It follows then that we can find a scattering FIO quantizing any given
contact transformation from a neighbourhood of a point $q \in \scT^*_{\pa
X} X$ to itself, since we may conjugate by a local Fourier transform and
reduce the problem to finding a conventional FIO quantizing a homogeneous
canonical transformation from a conic neighbourhood of $L(q) \in S^* \RR^n$
to itself. We can also use the local Fourier transform to import Egorov's
theorem into the scattering calculus. Namely, if $B \in \Psisc^{*, -1}(X)$
is a scattering pseudodifferential operator of order $-1$, with real
principal symbol, and $P \in \Psisc^{*, -1}(X)$ then also $e^{-iB} P e^{iB}
\in \Psisc^{*, -1}(X)$ is a scattering pseudodifferential operator of order
$-1$, whose symbol $p'$ is related to that of $P$ by the time $1$ flow of
the Hamilton vector field of $B$. This indeed is how we shall conjugate the
principal symbol $p$ of our operator to normal form.

\subsection{Normal form} In this section we put the principal symbol of $P$
into a normal form $p_{\norm}$. For later purposes we shall also need the
subprincipal symbol of $P$ in a normal form, but only only along the
`flow-out', \ie the unstable manifold, of $\Rp$, which can be done via
conjugation by a function; this is accomplished in
Lemma~\ref{effnonres-form}. (The model form of the subprincipal symbol only
plays a role in the polyhomogeneous, as opposed to just conormal, analysis,
which is the reason it is postponed to Section~\ref{sec:enr}.)

For this purpose, we only need to construct the principal
symbol $\sigma(B)$ of $B$ as in the first subsection.  This in turn can be
be written as $x^{-1}\tilde b$, $\tilde b\in\Cinf(\scT^* X)$, so we only
need to construct a function $b$ on $\scT^*_{\pa X}X$ such that the
pull-back $\Phi^* p$ of $p$ by the time $1$ flow $\Phi$ of $H_{x^{-1}\tilde
b}$ is the desired model form $p_{\norm}$, where $\tilde b$ is some
extension of $b$ to $\scT^* X;$ this property is independent of the chosen
extension. Thus {\em any} $B$ with $\sigma(B)=\tilde b$ will conjugate $P$
to an operator with principal symbol $p_{\norm}$. This construction is
accomplished in two steps, following Guillemin and Schaeffer
\cite{MR55:3504} in the non-resonant setting.  First we construct the
Taylor series of $b$ at $\Rp=(0,0,0)$, which puts $p$ into a model form
modulo terms vanishing to infinite order at $\Rp$. Next, we remove this
error {\em along the unstable manifold} of $\Rp$ by modifying an argument
due to Nelson \cite{Nelson1}.

Rather than using powers of $\cI$ to filter the Taylor series of $b,$ we
proceed as in \cite{MR55:3504} and assign degree $1$ to $y$ and $\mu$ but
degree two to $\nu$ in local coordinates as discussed above. Thus, let
$\frakh^j$ denote the space of functions
\begin{equation*}
\frakh^j=\sum_{2a+|\alpha|+|\beta|-2=j}
\nu^ay^\alpha\mu^\beta\Cinf(\scT^*_{\pa X}X)
\end{equation*}
Note that this is well-defined, independently of our choice of local
coordinates, since $-d\nu$ is the contact form $\alpha$ at $\Rp,$ so $\nu$
is well-defined up to quadratic terms. The Poisson bracket preserves this
filtration of $\cI$ in the following sense. If $\tilde a,$ $\tilde b$ are
some smooth extensions to $\scT^* X$ of elements $a\in\frakh^i$,
$b\in\frakh^j$ then
\begin{equation*}
x^{-1} \tilde c=\{x^{-1}\tilde a, x^{-1}\tilde b\}\Longrightarrow 
c=\tilde c|_{\scT^*_{\pa X}X}\in\frakh^{i+j}.
\end{equation*}
When this holds we write $c = \{\{a,b\}\}$; explicitly,
\begin{equation}
\{\{ a, b \}\} = W_a (b) + \frac{\partial a}{\partial \nu} b -
\frac{\partial b}{\partial \nu} a,
\label{bracket-defn}\end{equation}
with $W$ given by \eqref{W-coords}. Thus
\begin{equation}\label{eq:Poisson}
\{\{.,.\}\}:\frakh^i\times\frakh^j\mapsto\frakh^{i+j}.
\end{equation}
We then consider the quotient
\begin{equation*}
\frakg^j=\frakh^j/\frakh^{j+1},
\end{equation*}
so the bracket $\{\{.,.\}\}$ descends to 
\begin{equation*}
\frakg^i\times\frakg^j\to\frakg^{i+j}.
\end{equation*}

\begin{rem}
These statements remain true with $\frakh^j$ replaced by $\cI^j.$ However,
note that $p=-\nu$ in $\cI/\cI^2$, since $dp=-d\nu$ at $\Rp,$ but it is
{\em not the case} that $p=-\nu$ in $\frakg^0.$ In fact, $p$ is given by
\eqref{eq:p_0-def} below in $\frakg^0.$
\end{rem}

Using contact coordinates as discussed above, $\frakg^j$ may be freely
identified with the space of homogeneous functions of $\nu,y,\mu$ of degree
$j+2$ where the degree of $\nu$ is 2. Now let $p_0$ be the part of $p$ of
homogeneity degree two. In order to use Lemmas~\ref{HMV2r.191} and \ref{HMV2r.191p}, \emph{we assume throughout the paper from here on that case (iv) above Lemma~\ref{HMV2r.191} does not apply.} Hence  from \eqref{HMV2r.189}
\begin{equation}\label{eq:p_0-def}
p_0=\lambda\big(-\nu+\sum_{j=1}^{m-1} r_j y_j\mu_j+\sum_{j=m}^{n-1}Q_j(y_j,\mu_j)\big),
\ p-p_0\in\frakh^1.
\end{equation}
If we take $b\in\frakh^l,$ $l\geq 1$ and let
$\Phi$ be the time $1$ flow of $H_{x^{-1}b}$ then
\begin{equation}\label{eq:P'}
x\Phi^* (x^{-1}p)
=p+\{\{p,b\}\}=p+\{\{p_0,b\}\},\ \text{modulo}\ \frakh^{l+1}.
\end{equation}
This allows us to remove higher order term in the Taylor series of the symbol
successively provided we can solve the `homological equation'
\begin{equation*}
\{\{p_0,b\}\}=e\in\frakh^l,\ \text{modulo}\ \frakh^{l+1}.
\end{equation*}
Thus we need to consider the range of this linear map; its eigenfunctions are
easily found from the eigenfunctions of the linearization of $W.$

\begin{lemma} The (equivalence classes of the) monomials $p_0^a e^\alpha
f^\beta$ with $2a+|\alpha|+|\beta|=l+2$ satisfy
\begin{equation}
\begin{gathered}
\{\{p_0,p_0^a e^\alpha f^\beta\}\}=
R_{a,\alpha ,\beta}p_0^a e^\alpha f^\beta\text{ with eigenvalue}\\
R_{a,\alpha ,\beta}=\lambda\left(a-1+\sum_{j=1}^{n-1}\alpha_j
r_j+\sum\limits_{j=1}^{n-1}\beta_j(1-r_j)\right) 
\end{gathered} 
\label{HMV2r.194}\end{equation}
and give a basis of eigenvectors for $\{\{p_0,.\}\}$ acting on $\frakg^l.$

Here we identify the differentials $e_j$ and $f_j$ with linear functions
with these differentials.
\end{lemma}

\begin{rem}
In fact, the contact coordinates given by Lemma~\ref{HMV2r.191p} suffice
for the proof of this lemma; the additional information
in Lemma~\ref{HMV2r.191} is not needed. In this case, by
\eqref{HMV2r.189p},
\begin{equation}\label{eq:p_0-def-p}
p_0=\lambda\big(-\nu+\sum_{j=1}^{m-1}r_j e_jf_j+\sum_{j=m}^{n-1}Q_j(y_j,\mu_j)\big).
\end{equation}

We also remark that we could equally well use the eigenvector basis
for $\{\{p_0,.\}\}$ acting on $\frakg^l$ given by
$\nu^a e^\alpha
f^\beta$ with $2a+|\alpha|+|\beta|=l+2$. This follows from the lemma
using that $\nu=\sum_{j=1}^{m-1} r_j y_j\mu_j+\sum_{j=m}^{n-1}Q_j(y_j,\mu_j)
-\lambda^{-1}p_0$ in $\frakg^0$, and $y_j\mu_j$ as well as $Q_j(y_j,\mu_j)$
are eigenvectors with eigenvalue $\lambda(r_j+(1-r_j))=\lambda,$ and so
is $p_0$.
\end{rem}

\begin{proof}
Taking into account the eigenvalues and eigenvectors of $A,$ all
eigenvalues and eigenvectors of $\{\{p_0,.\}\}$ can be calculated
iteratively using the derivation property of the original Poisson bracket.
This implies
\begin{equation}\begin{split}\label{eq:product-bracket}
\{\{p_0,ab\}\}&=x\{x^{-1}p_0,x(x^{-1}a)(x^{-1}b)\}\\
&=x^{-1}\{x^{-1}p_0,x\}ab
+x\{x^{-1}p_0,x^{-1}a\}b+xa\{x^{-1}p_0,x^{-1}b\}\\
&=\lambda ab+\{\{p_0,a\}\}b
+a\{\{p_0,b\}\},
\end{split}\end{equation}
where each term within $\{.,.\}$ really uses a $\Cinf$ extensions of the
$a,$ $b,$ $p_0$ to $\scT^*X,$ followed by evaluation of the bracket and
then restriction to $\scT^*_{\pa X}X.$
Since
\begin{equation*}
\{\{p_0,a\}\}=x\{x^{-1}p_0,x^{-1}a\}=x\{x^{-1}p_0,x^{-1}\}a+\{x^{-1}p_0,a\}
=-\lambda a+\{x^{-1}p_0,a\},
\end{equation*}
on $\frakg^{-1}$ the eigenvectors of $\{\{p_0,.\}\}$ are the eigenvectors
$e_j$ and $f_j$ of $A$ with eigenvalues $-\lambda+\lambda r_j$ and
$-\lambda+\lambda(1-r_j).$ Moreover, in $\frakg^0$, $p_0$ is an eigenvector
of $\{\{p_0,.\}\}$ with eigenvalue $0.$ Thus, $e_j$, $f_j$ and $p_0$ satisfy
the claim of the lemma.
Since the other generators of $\frakg^0,$ as well
as generators of $\frakg^j$, $j\geq 1$, can be written as a products
of the $e_j$, $f_j$ and $p_0,$ the conclusion of the lemma
follows by induction.
\end{proof}

\begin{Def}\label{resdefn}
We call the multiindices in the set
\begin{equation}
I=\left\{(a,\alpha,\beta);R_{a,\alpha,\beta}=0\Mand
2a+|\alpha|+|\beta|\geq3\right\},
\label{HMV2r.197}\end{equation}
with $R_{a,\alpha ,\beta}$ given by \eqref{HMV2r.194}, {\em resonant}.
\end{Def}

Conjugation therefore allows us to remove, by iteration, all terms except
those with indices in $I.$ Expanding $p_0^a$ using \eqref{eq:p_0-def} we
deduce the following.

\begin{prop}\label{prop:linear}
If $P$ is as above and the leading term of $p=\sigma_{\pa,-1}(P)$ is
given by \eqref{eq:p_0-def} near a given radial point 
$\Rp$ then there exists a local contact diffeomorphism $\Phi$ near $\Rp$
such that 
\begin{equation}\begin{split}
\Phi^*p=&\lambda\big(-\nu+\sum_{j=1}^{m} r_j y_j\mu_j+\sum_{j=m+1}^{n-1}Q_j(y_j,\mu_j)
+\sum_{(a,\alpha,\beta)\in I} c_{a,\alpha,\beta}\nu^a e^\alpha f^\beta\big)\\
&\qquad\qquad \text{modulo}\ \cI^\infty=\frakh^\infty\text{ at }\Rp
\label{res-terms}\end{split}\end{equation}
with $I$ given by \eqref{HMV2r.197}. 
\end{prop}

\begin{proof} The Taylor series of $\Phi$
at $\Rp$ can be
  constructed inductively over the filtration $\frakh^j$ as indicated
  above. At the $j$th stage, the terms of weighted homogeneity $j$ can be
  removed from $p$ except for those in the null space of
  $\{\{p_0,\cdot\}\},$ \ie the resonant terms with $R_{a,\alpha ,\beta
  }=0.$
This
  leads to \eqref{res-terms} in the sense of formal power series. However,
  by use of Borel's Lemma a local contact diffeomorphism
can be found giving \eqref{res-terms}.
\end{proof}

Now a small extension of Nelson's proof of Sternberg's linearization
theorem can be used to remove the infinite order vanishing error along the
unstable manifold, \ie at $\nu=0,$ $\mu=0,$ $y''=0,$ $y'''=0.$

\begin{prop}\label{prop:Nelson}
Suppose that $X$ and $X_0$ are $\Cinf$ vector fields on $\Real^N$ with
$X_0(0)=0$ and $X_1=X-X_0$ vanishing to infinite order at $0.$ Suppose also that they are both
linear outside a compact set and equal there to their common linearization,
$DX(0),$ at $0$ which is assumed to have no pure imaginary eigenvalue. Let
$U(t),$ $U_0(t)$ be the flows generated by $X$ and $X_0.$ If $E$ is a linear
submanifold invariant under $X_0$ such that
\begin{equation}\label{eq:U_0-contracts}
\lim_{t\to\infty}U_0(t)x=0\ \forall\ x\in E
\end{equation}
then for all $j=0,1,2,\ldots$ and $x\in E$
\begin{equation}\label{eq:W-derivs}
\lim_{t\to\infty}D^j(U(-t)U_0(t))x
\end{equation}
exists, and is continuous in $x\in E,$ and
\begin{equation*}
W_-x=\lim_{t\to\infty}U(-t)U_0(t)x,\ x\in E
\end{equation*}
has a $\Cinf$ extension, $G,$ to $\Real^N$ which is the identity
to infinite order at $0$ and such that $(G^{-1})_*X=X_0$
to infinite order along $E$ in a neighbourhood of $0.$
\end{prop}

\begin{rem}
Note that the derivatives $D^j$ in \eqref{eq:W-derivs} refer to the
ambient space $\Real^N,$ and {\em not} merely to $E.$ This is useful
in producing the Taylor series of $G$ for the last part of the
conclusion.

Also, the limit $t\to\infty$ means $t\to +\infty,$ as in Nelson's book.
\end{rem}

\begin{proof} We follow the proof of Theorem 8 in \cite{Nelson1}. Indeed,
if $X_0$ was assumed to be linear then Nelson's theorem would apply
directly. Dropping this assumption has little effect on the proof; the main
difference is that a little more work is required to show the exponential
contraction property, \eqref{eq:U_0-t-est} below.

Since the real part of every eigenvalue of $DX(0)$ is non-zero,
$\Real^N=E_+\oplus E_-$ where $E_+,$ resp.\ $E_-,$ is the direct sum of the
generalized eigenspaces of $DX(0)$ with eigenvalues with positive, resp.\
negative, real parts. Since $E$ is invariant under $X_0,$ and hence under
$DX(0),$ necessarily $E\subset E_-.$ We actually apply the theorem with
$E=E_-,$ but, as in Nelson's discussion, the more general case is useful
for the inductive argument for the derivatives.

Let $e_j$ denote a basis of $E_-$ consisting of generalized eigenvectors of
$DX(0)$ with corresponding eigenvalue $\sigma_j;$ we shall consider the
$e_j$ as differentials of linear functions $f_j$ on $\Real^N.$ For $x\in
\Real^N,$ let $x(t)=U_0(t)x,$ $F_j(t)=f_j(x(t)).$ Then
$\frac{dF_j}{dt}|_{t=t_0}=(X_0 f_j)(x(t_0))$ where
\begin{equation*}
X_0 f_j(y)=DX(0)f_j(y)+\cO(\|y\|^2).
\end{equation*}
Moreover, for $y\in E_-$, $\|y\|^2\leq C_1\sum_j f_j^2$ for some $C_1>0.$
So, setting $\rho=\sum f_j^2,$ we deduce that
\begin{equation*}
X_0 \rho(y)=\sum_j 2\sigma_j f_j^2(y) +\cO(\rho(y)^{3/2}),
\end{equation*}
hence with $R(t)=\rho(x(t)),$ $c_0\in(\sup\sigma_j,0),$ there exists
$\delta>0$ such that for $\|R(t)\|\leq\delta,$
\begin{equation*}
\frac{dR}{dt}-2c_0 R \leq 0,
\end{equation*}
and hence $R(t)\leq e^{-2c_0 t}\|x\|$ for $t\geq 0,$ $\|r(x)\|\leq\delta$,
$x\in E_-.$ A corresponding estimate also holds outside a compact set,
as $X_0$ is given by $DX(0)$ there, so a patching argument and
\eqref{eq:U_0-contracts} yield the estimate $R(t)\leq C_0 e^{-2c_0 t}\|x\|$
for all $x\in E_-.$ Since $R(t)^{1/2}$ is equivalent to $\|.\|,$ we deduce that
there are constants $C,$ $c>0$ such that
\begin{equation}\label{eq:U_0-t-est}
\|U_0(t)x\|\leq Ce^{-ct}\|x\|\ \forall\ x\in E\text{ and }t\geq 0.
\end{equation}

For the remainder of the argument we can follow Nelson's proof even more
closely. Thus, let $\kappa$ be a Lipschitz constant for $X$ and $X_0,$ and
choose $m$ such that $cm>\kappa.$ Note that there exists $c_0>0$ such that
for all $x\in\Real^N,$ 
\begin{equation}\label{eq:X_1-est-8}
\|X_1(x)\|\leq c_0\|x\|^m,\ X_1=X-X_0.
\end{equation}
For $t_1\geq t_2\geq 0,$ $t_1=t_2+t,$ $x\in E,$
\begin{equation*}\begin{split}
I=\|U(-t_1)U_0(t_1)x-U(-t_2)U_0(t_2)x\|&=\|U(-t_2)\left(U(-t)U_0(t)-\Id\right)
U_0(t_2)x\|\\
&\leq e^{\kappa t_2}\|(U(-t)U_0(t)-\Id)U_0(t_2)x\|
\end{split}\end{equation*}
by the Lipschitz condition (see \cite[Theorem~5]{Nelson1}). But with
$X=X_0+X_1$, by \cite[Proof of Theorem~6, (5)]{Nelson1}
\begin{equation*}
\|U(-t)U_0(t)y-y\|\leq\int_0^t e^{\kappa s}\|X_1(U_0(s)y)\|\,ds.
\end{equation*}
Applying this with $y=U_0(t_2)x$, we deduce that
\begin{equation}\label{eq:I-est-8}
I\leq
e^{\kappa t_2}\int_0^t e^{\kappa s}\|X_1(U_0(s+t_2)x)\|\,ds.
\end{equation}
Thus, by \eqref{eq:X_1-est-8} and \eqref{eq:U_0-t-est},
\begin{equation*}\begin{split}
I&\leq e^{\kappa t_2}\int_0^t e^{\kappa s}c_0 C^m e^{-cm(s+t_2)}\|x\|^m\,ds\\
&\leq e^{\kappa t_2}\int_0^\infty e^{\kappa s}c_0 C^m e^{-cm(s+t_2)}\|x\|^m\,ds
=\frac{c_0 C^m e^{-(cm-\kappa)t_2}\|x\|^m}{cm-\kappa}.
\end{split}\end{equation*}
Letting $t_2\to\infty$ shows that $W_- x=\lim_{t\to\infty} U(-t)U_0(t)x$
exists, with convergence uniform on compact sets, hence $W_-$ is continuous
in $x\in E.$ Moreover, applying the estimate with $t_2=0$ shows that
$W_-(x)-x=\cO(\|x\|^m).$ Since $m$ is arbitrary, as long as it is sufficiently
large, this shows that $W_-$ is the identity to infinite order at $0,$
provided it is smooth, as we proceed to show.

Smoothness can be seen by a similar argument, although we need to put
a slight twist into Nelson's argument. Namely, first consider the first
derivatives, or rather the 1-jet. Thus, we work on $\Real^N\oplus
\cL(\Real^N).$ Let $(x,\xi)$ denote the components 
with respect to this decomposition. These evolve under the 
flow $U'(t),$ resp.\ $U'_0(t),$ given by
\begin{equation*}
X'(x,\xi)=(X(x),DX(x)\cdot\xi),\ X'_0(x,\xi)=(X_0(x),DX_0(x)\cdot\xi),
\end{equation*}
where $DX(x)$ and $\xi$ are considered as elements of $\cL(\Real^N),$
and $\cdot$ is composition of operators. Note that
the second, $\cL(\Real^N)$, component of these vector fields is a homogeneous
degree zero vector field, i.e.\ it is invariant under push-forward
by the natural $\Real^+$-action (by dilations).

The twist, as compared to Nelson's work, is that we identify $\cL(\Real^N)$
with $\Real^{N^2},$ which we radially compactify to a (closed) ball
$B^{N^2},$ which we further embed as the closed unit ball in $\Real^{N^2}$
in such a fashion that the smooth structure of the ball agrees with the
restriction of the smooth structure from $\Real^{N^2}$. Let
$\iota:\Real^{N^2}\to\Real^{N^2}$ be this map with range the interior of
$B^{N^2}$. Then the push-forward under $\iota$ of a homogeneous degree zero
vector field, such as $DX(x)\cdot\xi$ is for each $x\in\Real^N,$ extends to
a $\Cinf$ vector field on the closed ball $B^{N^2},$ which by homogeneity
is tangent to the boundary. Furthermore, if
$\iota_1=\mathrm{id}_{\Real^N}\times \iota,$ then $(\iota_1)_*X'$ and
$(\iota_1)_*X'_0$ extend to $\Cinf$ vector field on $\Real^N\times B^{N^2}$
tangent to the boundary and their difference, $(\iota_1)_*X'_1$, in
addition vanishes to infinite order at $\{0\}\times B^{N^2}$. Thus
$(\iota_1)_*X'$ and $(\iota_1)_*X'_0$ are Lipschitz with some Lipschitz
constant $\kappa'$: this is automatic over a compact subset of
$\Real^N\times B^{N^2},$ which in fact suffices here, but in fact holds on
all of $\Real^N\times B^{N^2}$ since outside the inverse image of a compact
subset of $\Real^N\times B^{N^2}$, $X'$ and $X'_0$ are linear, so in
particular their $B^{N^2}$ component is independent of $x.$

To minimize confusion about the `change of coordinates', we write
the coordinates on $\Real^N\times B^{N^2}$ as $(x,\eta)$ below.
With $c$ as in \eqref{eq:U_0-t-est}, choose
$m$ such that $cm>\kappa'.$ Then the infinite order vanishing of
$(\iota_1)_*X'_1$ at $x=0$ yields
\begin{equation}\label{eq:X_1p-est-8}
\|((\iota_1)_*X_1')(x,\eta)\|\leq c'_0\|x\|^m
\end{equation}
for all $(x,\eta)$.
Let $U'(t),$ $U_0'(t)$ denote the evolution groups generated by 
$(\iota_1)_*X'$ and $(\iota_1)_*X'_0,$ respectively.
Thus, for all real $t,$
\begin{equation}\label{eq:Up_0-t-est}
\|U'(t)(x,\eta)\|\leq e^{\kappa't}\|(x,\eta)\|,
\end{equation}
see \cite[Theorem~5]{Nelson1}. So \eqref{eq:I-est-8} still applies, with
$X_1$ replaced by $(\iota_1)_*X_1',$
$\kappa$ replaced by $\kappa',$ etc. Thus,
by \eqref{eq:U_0-t-est} and \eqref{eq:Up_0-t-est},
\begin{equation*}\begin{split}
I'&=\|U'(-t_1)U_0'(t_1)(x,\eta)-U'(-t_2)U_0'(t_2)(x,\eta)\|\\
&\leq e^{\kappa' t_2}\int_0^t e^{\kappa' s}\|((\iota_1)_*X'_1)(U'_0(s+t_2)
(x,\eta))\|\,ds\\
&\leq e^{\kappa' t_2}\int_0^t e^{\kappa' s}c'_0 C^m e^{-cm(s+t_2)}\|x\|^m\,ds\\
&\leq e^{\kappa' t_2}\int_0^\infty e^{\kappa' s}c'_0 C^m e^{-cm(s+t_2)}\|x\|^m
\,ds=\frac{c'_0 C^m e^{-(cm-\kappa')t_2}\|x\|^m}{cm-\kappa'}.
\end{split}\end{equation*}
Thus, $\lim_{t\to\infty} U'(-t)U'_0(t)x$ exists, with convergence uniform on
compact sets, so the limit depends continuously on $(x,\xi)$ for $x\in E.$

The higher derivatives can be handled similarly. The resulting Taylor
series about $E$ can be summed asymptotically, giving $G$: this part of the
argument of Nelson is unchanged.
\end{proof}

\subsection{Effective resonance and nonresonance}
Next we apply this general result to the symbol $p.$ Following
Lemma~\ref{HMV2r.191}, when resonances occur we cannot remove all error
terms even in the sense of formal power series. Consequently we do not attempt to
get a full normal form in a neighbourhood of the critical point, but only
along the submanifold
\begin{equation}
S=\{\nu=0,\ y''=0,\ y'''=0,\ \mu=0\},
\label{HMV2r.198}\end{equation}
which is the unstable manifold for $W_0.$ After reduction to normal form,
errors which are polynomial in the normal directions to $S$ will
remain. For later purposes, we divide these into two parts.

\begin{Def}\label{Def:effres}
With $I$ as in Definition~\ref{resdefn}, let
\begin{equation}\begin{split}\label{eq:effres-I}
&I_{\effr}=I_{\effr}'\cup I_{\effr}'',\\
&I'_{\effr}=\{(a,\alpha,\beta)\in I:\ \alpha=(\alpha',\alpha'',\alpha'''),
\ \beta=(\beta',\beta'',\beta'''),\\
&\qquad\qquad a=0,\ \alpha'''=0,\ \beta'''=0,
\ \alpha''=0,\ \beta''=0,\ |\beta'|=1\},\\
&I''_{\effr}=\{(a,\alpha,\beta)\in I:\ \alpha=(\alpha',\alpha'',\alpha'''),
\ \beta=(\beta',\beta'',\beta'''),\\
&\qquad\qquad\ a=0,\ \alpha'''=0,\ \beta'''=0,
\ \alpha'=0,\ \beta'=0\}.
\end{split}\end{equation}
An {\em effectively
resonant} function is a polynomial of the form
\begin{equation*}
r_{\effr} = \sum_{(a,\alpha,\beta)\in I_{\effr}}c_{a,\alpha,\beta}
\,p_0^a e^\alpha f^\beta,
\end{equation*}
or equivalently
\begin{equation*}
r_{\effr} = \sum_{(a,\alpha,\beta)\in I_{\effr}}c_{a,\alpha,\beta}
\,\nu^a e^\alpha f^\beta.
\end{equation*}
\end{Def}

Thus, elements of $I_{\effr}$
satisfy $(0,\alpha,\beta)\in I$ (i.e. are resonant, see Definition~\ref{resdefn}), with $\alpha=(\alpha',\alpha'',0)$,
$\beta=(\beta',\beta'',0)$, and either $\alpha''=0$, $\beta''=0$, $|\beta'|=1$,
or $\alpha'=0$, $\beta'=0$.

Moreover, an effectively resonant function has the form
\begin{equation}
\sum_{\alpha', |\beta'| = 1} c_{\alpha' \beta'} (e')^{\alpha'}
(f')^{\beta'} + \sum_{\alpha'', \beta''} c_{\alpha'' \beta''}
(e'')^{\alpha''} (f'')^{\beta''}.
\label{effres-span}\end{equation}

For a fixed critical point of a fixed operator $P$ (e.g. $P = x^{-1}(\Delta + V - \sigma)$ for a fixed $\sigma$), the set $I_{\effr}$ is finite. Thus,
only a finite number of terms can occur in
\eqref{effres-span}, and hence restricting to polynomials in the definition
of effectively resonant functions (rather than infinite formal sums)
is in fact not a restriction.
To see this, note that in the expression for $R_{a, \alpha, \beta}$ in \eqref{HMV2r.197}, we have $a = 0, \, \alpha''' = \beta''' = 0$ and either (i) $\alpha'' = \beta'' = 0$ and $|\beta'| = 1$ or (ii) $\alpha' = \beta' = 0$. In case (i), if $\beta'_j = 1$ then to have $R_{a, \alpha, \beta} = 0$ we need $\sum \alpha_k' r_k' = r_j'$, which is only possible for $|\alpha'| \leq |r_j'|/\min_k |r_k'|$. In case (ii), we need 
$\sum \alpha''_j r''_j + \sum \beta''_j (1 - r''_j) = 1$, which is only possible for $|\alpha''| \leq 1/\min r_k''$ and $|\beta''| \leq 2$. (Actually in case (ii) we must have $|\beta''| \leq 1$ in order to satisfy the condition $2a + |\alpha| + |\beta| \geq 3$ in \eqref{HMV2r.197}.)

\begin{Def}\label{Def:effnonres}
Let $\cJS$ denote the ideal of $\Cinf$ functions on $\scT^*_{\pa X}X$ which
vanish on $S$ and set 
\begin{equation}
J''=\left\{(\alpha'',\beta'');
\sum_{j=s}^{m-1}r''_j\alpha''_j+(1-r''_j)\beta''_j\in(1,2)\right\}.
\label{eq:I''-def}\end{equation}
An {\em effectively nonresonant} function is an element of $\cJS$ of the form
\begin{multline}
r_{\effnr}=\sum_{j=1}^{s-1} h_jf'_j+\sum_{(\alpha'',\beta'')\in I''}
h''_{\alpha'',\beta''}e^{\alpha''}f^{\beta''}+\sum_{j,k} h'''_{jk}e'''_j f'''_k\\
h_j\in\cJS,\ j=0,1,\ldots,s,
\ h''_{\alpha'',\beta''}\in\Cinf(\scT^*_{\pa X}X),
\ (\alpha'',\beta'')\in I'',\\
h'''_{jk}\in \cJS,\ j,k=m,\ldots,n-1.
\label{effnonres-span}\end{multline}
\end{Def}

Note that $J''$ is finite, hence all sums in the definition are finite.

\begin{thm}\label{thm:model} Using the notation of Lemma~\ref{HMV2r.191}
for coordinates near a radial point of $\Rp$ of $p$ there is a local contact
diffeomorphism $\Phi$ from a neighbourhood of $(0,0,\ldots,0)$ to a
neighbourhood of $\Rp$ such that $\Phi^*p=p_{\norm}$ such that
\begin{equation} 
\lambda^{-1}p_{\norm}=-\nu+\sum_j r_j y_j\mu_j+\sum_{j=m}^{n-1}Q_j(y_j,\mu_j) 
+ r_{\effnr} + r_{\effr},
\label{res-terms-2}\end{equation}
with $r_{\effnr}$ of the form \eqref{effnonres-span} and $r_{\effr}$ of the
form \eqref{effres-span}; in addition at a non-resonant critical point, \ie
if $I=\emptyset,$ then we may take $r_{\effnr}=r_{\effr}=0$ near $\Rp.$
\end{thm}

\begin{rem}\label{rem:model-op}
If $F$ is an elliptic Fourier integral operator with canonical relation $\Phi$
then $\tilde P=F^{-1}PF$
satisfies $\sigma_{\pa,-1}(\tilde P)=p_{\norm}.$
\end{rem}

\begin{rem} As will be seen below, of the two error terms, only $r_{\effr}$
has any effect on the leading asymptotics of microlocal solutions. The
construction below shows that modulo $\cI^\infty,$ $r_{\effnr}$ may be
chosen to consist of resonant terms only, \ie to be an asymptotic sum of resonant
terms. However, this plays no role in the paper; all the relevant
information is contained in the statement of the theorem.
\end{rem}

\begin{rem}\label{rem:special-contact} We do not need the full power of
Lemma~\ref{HMV2r.191} to find $\Phi$ as in this theorem;
Lemma~\ref{HMV2r.191p} suffices. Indeed, the terms $\sum_{j=1}^{m-1}
c_j\mu_j^2$ in \eqref{HMV2r.189p} can be absorbed in $r_{\effnr}$.

Similarly any term $\nu^a\mu^\beta y^\alpha$ with $a+|\beta|\geq 2$ and
$a\neq 0,$ or with $|\beta|\geq 3$ can be included in $r_{\effr}$ or
$r_{\effnr}.$ The same is true for any term with $|\beta|\geq 2$ such that
$\beta_j\neq 0$ for some $j$ with $\re r_j\neq \frac{1}{2}.$ In particular,
if $\re r_j\neq \frac{1}{2}$ for any $j,$ the only terms which need to be
removed have $a+|\beta|\leq 1.$ The conjugating Fourier integral operator
can therefore also be arranged to have such terms only and thus to be of
the form $e^{iB},$ with $B=Z+(f/x)$ where $Z$ is a vector field on $X$
tangent to its boundary and $f$ is a real valued
smooth function on $X.$ Correspondingly, the normal form may be
achieved by conjugation of $P$ by an oscillatory function, $e^{if/x},$
followed by pull-back by a local diffeomorphism of $X,$ \ie a 
change of coordinates.
However, if $\re r_j= \frac{1}{2}$ for some $j,$ some
quadratic terms in $\mu$ would also need to be removed for the model form,
but since they play a role analogous to $r_{\effr},$ the arguments of
Section~\ref{sec:modules}, giving conormality, are unaffected, and only the
polyhomogeneous statements of Section~\ref{sec:enr} would need
alterations. However, the contact diffeomorphism (\ie FIO conjugation)
approach we present here is both more unified and more concise.

If $p=|\zeta|^2+V_0-\ev$,
the model form of Lemma~\ref{HMV2r.191p} also only required
a change of coordinates and multiplication by an oscillatory function
(see Lemma~\ref{lemma:special-contact}),
the model form of this theorem can be obtained
by these two operations, starting from the original operator $P$ with symbol
$p$.
\end{rem}

\begin{proof} First we apply Proposition~\ref{prop:linear}. Next we need to
show that $r_{\effr}$ as in \eqref{effres-span} and $r_{\effnr}$ as in
\eqref{effres-span} can be chosen to have Taylor series at $0$ given
exactly by the error term in \eqref{res-terms}.

So, consider a monomial $\nu^a e^\alpha f^\beta$ with $(a,\alpha,\beta)\in
I.$ If $\alpha'''\neq 0$ then $\beta'''\neq 0$ since $\im r'''_j>0,$ and
only the eigenvalues of $f'''_j$ have negative imaginary parts, and
conversely. In addition, $2a+|\alpha|+|\beta|\geq 3$ implies that a
monomial with $\alpha'''\neq 0$ or $\beta'''\neq 0$ has the form $\nu^a
e^{\tilde\alpha} f^{\tilde\beta} e'''_jf'''_k$ for some $j,k$ with
$2a+|\tilde\alpha|+|\tilde\beta|\geq 1$ and
\begin{equation*}
\re(a+\sum r_l\tilde\alpha_l+\sum (1-r_l)\tilde\beta_l)=0.
\end{equation*}
Since $\re (1-r_l)>0$ for all $l$ and $\re r_l>0$ for $l\geq s,$ while
$r_l<0$ for $l\leq s-1,$ we must have $\tilde\alpha'\neq 0$ (\ie
$\tilde\alpha_l\neq 0$ for some $l\leq s-1$) and correspondingly
$a+|\tilde\alpha''|+|\tilde\alpha'''|+ |\tilde\beta|>0$. Due to the latter,
$\nu^a e^{\tilde\alpha} f^{\tilde\beta}$ vanishes on $S,$ so the terms with
$\alpha'''\neq 0$ or $\beta'''\neq 0$ appear in $r_{\effnr}$.

So we may assume that $\alpha'''=\beta'''=0.$ If $a\neq 0,$ the monomial is
of the form $\nu^{\tilde a} e^{\tilde\alpha} f^{\tilde\beta} \nu,$ $\tilde
a=a-1,$ $2\tilde a+|\tilde\alpha|+|\tilde\beta|\geq 1$ with
\begin{equation*}
\tilde a+\sum r_j\tilde\alpha_j+\sum (1-r_j)\tilde\beta_j= 0.
\end{equation*}
Arguing as in the previous paragraph we deduce that the terms with $a\neq
0$ also appear in $r_{\effnr}.$

So we may now assume that $a=0,$ $\alpha'''=\beta'''=0$. If $\beta'\neq 0,$
the monomial is of the form $\nu^{a} e^{\tilde\alpha} f^{\tilde\beta} f_j$
for some $j,$ and $2 a+|\tilde\alpha|+|\tilde\beta|\geq 2,$
\begin{equation*}
a+\sum r_l\tilde\alpha_l+\sum (1-r_l)\tilde\beta_l=r_j< 0.
\end{equation*}

We can still conclude that $\tilde\alpha'\neq 0,$ but it is not automatic
that $a+|\tilde\alpha''|+ |\tilde\beta|>0.$ However, if
$a+|\tilde\alpha''|+|\tilde\beta|>0$ then $\nu^{a} e^{\tilde\alpha}
f^{\tilde\beta} f_j$ is again included in $r_{\effnr},$ while if
$a+|\tilde\alpha''|+|\tilde\beta|=0,$ then the monomial is included in
$r_{\effr}.$

Finally then, we may assume that $a=0,$ $\beta'=0,$ $\alpha'''=\beta'''=0.$
Since $r'_j<0$ for all $j=1,\ldots,s-1$
\begin{equation*}
\sum (r''_j\alpha''_j+(1-r''_j)\beta''_j)\geq
\sum r'_j\alpha'_j+\sum (r''_j\alpha''_j+(1-r''_j)\beta''_j)=1.
\end{equation*}
Moreover, the equality holds if and only if $\alpha'=0,$ in which case this
term is included in $r_{\effr}.$ The terms with $\alpha'\neq 0$ can be
included in
$h''_{\tilde\alpha'',\tilde\beta''}e^{\tilde\alpha''}f^{\tilde\beta''}$ for
some $\tilde\alpha''\leq\alpha''$, $\tilde\beta''\leq\beta'',$ chosen by
reducing $\alpha''$ and/or $\beta''$ to make
\begin{equation*}
\sum (r''_j\tilde\alpha''_j+(1-r''_j)\tilde\beta''_j)\in(1,2).
\end{equation*}
This can be done since $r''_j,$ $1-r''_j\in (0,1).$

It follows that $\lambda^{-1}p$ can be conjugated to the form
\begin{equation} 
-\nu+\sum_j r_j y_j\mu_j+\sum_{j=m}^{n-1}Q_j(y_j,\mu_j) 
+ r_{\effnr} + r_{\effr} + r_\infty,
\label{res-terms-1b}\end{equation}
where $ r_{\effnr}, r_{\effr}$ are as in \eqref{effnonres-span},
\eqref{effres-span}, with both vanishing if $\Rp$ is non-resonant, and
$r_{\infty}$ vanishes to infinite order at $(0,0,0).$ Thus, it remains to
show that we can remove the $r_\infty$ term in a neighbourhood of the origin.

To do this we apply Proposition~\ref{prop:Nelson}. Let $X'$ be the
Legendre vector field of \eqref{res-terms-1b}, and let $X'_1$ be the
Legendre vector field of $r_\infty,$ while $X'_0=X'-X'_1.$  Let $\tilde X$
be the linear vector field with differential equal to $DX(0),$ let $\chi$
be compactly supported, identically $1$ near $0$, and let $X=-(\chi X'+(1-\chi)\tilde X),$  $X_0=-(\chi X'_0+(1-\chi)\tilde X).$
The overall minus sign is due to $S$ being the unstable manifold
of $X'_0$ near the origin, hence the stable manifold of $-X'_0.$
Let $E$ be the subspace $S$ of $\Real^{2n-1},$ defined by
\eqref{HMV2r.198}. Then Proposition~\ref{prop:Nelson} is applicable, and
$G$ given by it may be chosen as a contact diffeomorphism since $U(t),$
$U_0(t)$ are such, see \cite[Section~3, Theorem~4]{MR55:3504}.
\end{proof}

\subsection{Parameter-dependent normal form}
We also need a parameter-dependent version of this theorem. Namely,
suppose that $p$ depends smoothly on a parameter $\ev,$ can we make the
normal form depend smoothly on $\ev$ as well? This problem can be
approached in at least two different ways. One can consider $\ev$ simply
as a parameter, so $p\in\Cinf((\pa\scT^*X)\times I)=\Cinf((\scT^*_{\pa
X}X)\times I)$ and then try to carry out the reduction to normal form
uniformly. Alternatively, one identify $p$ with the function $p'$
on the larger space $\pa\scT^*(X\times I)$ arising by the pull-back under
the natural projection
\begin{equation*}
p'=\pi^* p,\ \pi:\scT^*_{\pa X\times I}(X\times I)\to(\scT^*_{\pa X}X)\times I
\end{equation*}
and then carry out the reduction to a model on the larger space. Whilst the
second approach may be more natural from a geometric stance, we will adopt
the first, since it is closer to the point of view of spectral theory of
\cite{HMV1}. Clearly the difficulty in obtaining a uniform normal form is
particularly acute near a value of $\ev$ at which the effectively resonant
terms do not vanish. Fortunately in the case of central interest here, and
in other cases too, the set of points at which such problems arise is
discrete.

\begin{lemma}\label{lemma:er-energies}
If $P=P(\ev)=x^{-1}(\Delta+V-\ev),$ $\Rp=\Rp(\ev)$ is a radial point of $P$
lying over the critical point $z = \pi(\Rp)$ of $V_0$ and $I(\ev),$ resp.\ $I_{\effr}(\ev),$ are the sets
\eqref{HMV2r.197}, resp.\ \eqref{eq:effres-I},
for $p(\ev)$ then the set $\cR_{z} = \cR_{\hesst, z} \cup \cR_{\effr,z}$, defined by 
\begin{equation}\begin{split}
\cR_{\hesst, z} &= \big\{\ev\in(V_0(z),+\infty) \mid \exists j \text{ such that } r_j = \frac1{2} \big\}, \\
\cR_{\effr,z} &= \big\{\ev\in(V_0(z),+\infty) \mid I_{\effr}(\ev)
\neq\emptyset\big\},
\label{HMV2r.199}
\end{split}\end{equation}
i.e., the set of energies $\sigma$ which are either a Hessian threshold (see Lemma~\ref{HMV2r.191}) or such that $q(\sigma)$ has a nontrivial effectively resonant error term
(see Definition~\ref{Def:effres}), is discrete in $(V_0(\pi(\Rp)),+\infty).$
\end{lemma}

\begin{rem}
It follows that if $K\subset (V_0(z),+\infty)$ is compact then
$K\cap\cR_{z}$ is finite. Thus, to prove properties such as
asymptotic completeness, one can ignore all $\ev\in K$ which are Hessian thresholds or effectively resonant. 

Note also that by the definition of $I_{\effr}(\ev),$
\begin{multline*}
\cR_{\effr,z}= \big\{\ev\in(V_0(z),+\infty) \mid
\text{either }\exists\ (0,(\alpha',0,0),(\beta',0,0))\in I(\ev)
\Mwith |\beta'|=1\\
\text{ or }\exists\ (0,(0,\alpha'',0),(0,\beta'',0))\in I(\ev)\big\}.
\end{multline*}
\end{rem}

\begin{proof} Using Remark ~\ref{rem:H-lin-ev},  the set $\cR_{\hesst, z}$ of Hessian thresholds is given by $\{ V_0(z) + 4 a_j \}$ where $a_j$ is an eigenvalue of the Hessian of $V_0$ at $z$ and hence has cardinality at most $n-1$, so this set is trivially discrete.

Let $K$ be a compact subset of $(V_0(z),+\infty).$ 
The set
$K\cap\cR_{\effr,z}$ of effectively resonant energies in $K$ is the the
union of zeros of a finite number of analytic functions (none of which are
identically zero). Indeed,  $\cR_{\effr,z}$
is given by the union of the set of zeros of the countable collection of
functions
\begin{equation*}
-1+\sum_{j=s}^{m-1} \alpha''_j r''_j(\ev)+\beta''_j(1-r''_j(\ev)),
\quad -1+(1-r_k)+\sum_{j=1}^{s-1}\alpha'_j r'_j(\ev)
\end{equation*}
as $k=1,\ldots,s-1,$ while $\alpha',$ $\alpha'',$ $\beta''$ are
multiindices.  But if $c>0$ is large enough then $c^{-1}>|r_j(\ev)|>c$ for
all $j$ and for all $\ev\in K$ as $K$ is compact and the $r_j$ do not
vanish there.  Correspondingly, for $|\alpha'|>\frac{2}{c^2},$
\begin{equation*}
-1+(1-r_k)+\sum_{j=1}^{s-1}\alpha'_j r'_j(\ev)<-r_k-|\alpha'|c<-c^{-1},
\end{equation*}
and analogously for $|\alpha''|+|\beta''|>\frac{2}{c},$
\begin{equation*}
-1+\sum_{j=s}^{m-1} \alpha''_j r''_j(\ev)+\beta''_j(1-r''_j(\ev))
>-1+(|\alpha''|+|\beta''|)c>1.
\end{equation*}
Thus, there are only a finite number of these analytic functions that may
vanish in $K,$ as claimed.
\end{proof}

If $q(\ev)$ are the radial points corresponding to $z \in \Cv(V)$, and $\ev \notin \cR_{\effr,z}$,  then we will say that $q(\ev)$ is effectively nonresonant, or that $\ev$ is an effectively nonresonant energy for $z$. We now prove that, away from effectively resonant energies and Hessian thresholds, we have a normal form for $p(\sigma)$ of the form \eqref{res-terms-2} with $r_{\effr} = 0$ and depending smoothly on $\sigma$. Thus, for a given critical point $z$ of $V_0$, consider an open interval $O\subset
(V_0(z),+\infty)\setminus\cR_{z}.$ Apart from the
coefficients $h_j,$ $h''_{\alpha'',\beta''},$ etc., in
\eqref{effnonres-span} the only part of the model form depending on $\ev$
is
\begin{equation*}
J''(\ev)=\{(\alpha'',\beta'');
\ \sum_{j=s}^{m-1} r''_j(\ev)\alpha''_j+(1-r''_j(\ev))\beta''_j\in(1,2)\}.
\end{equation*}
We note that on compact subsets $K$ of $O$, there is a $c>0$
such that $r''_j(\ev)>c$ for $\ev\in K$, and then
for $|\alpha''|+|\beta''|>2c^{-1}$,
\begin{equation*}
\fraks_{\alpha''\beta''}(\ev)
=\sum_{j=s}^{m-1} r''_j(\ev)\alpha''_j+(1-r''_j(\ev))\beta''_j>2,
\end{equation*}
so if we let
\begin{equation*}
J_K=\cup_{\ev\in K} J''(\ev),
\end{equation*}
then $J_K$ is a finite set of multiindices. For each multiindex
$(\alpha'',\beta'')$ we let 
\begin{equation}
O_{\alpha'',\beta''}=\fraks_{\alpha''\beta''}^{-1}((1,2)),
\label{HMV2r.201}\end{equation}
which is thus an open subset of $O$.

For the parameter dependent version of the Theorem~\ref{thm:model} we
introduce  
\begin{equation}
\cS=\{(y,\nu,\mu,\ev);\nu=0,\ y''=0,\ y'''=0,\ \mu=0,\ \ev\in O\},
\label{HMV2r.200}\end{equation}
in place of $S$ \eqref{HMV2r.198}. 

\begin{thm}\label{thm:model-smooth}
Suppose that $p\in\Cinf(\scT^*_{\pa X}X\times O),$
$O\subset(V_0(z),+\infty)\setminus\cR_{z}$ is open, that the
symplectic map $S$ induced by the linearization $A'$ of $p$ at $\Rp(\ev)$ (see
Lemma~\ref{HMV2r.192}) can be smoothly decomposed (as a function of $\ev\in
O$) into two-dimensional invariant symplectic subspaces and that there exists
$c>0$ such that $r''_j(\ev)>c$ for $\ev\in O$. Then $\Phi(\ev)$ and $F(\ev)$
can be chosen smoothly in $\ev$ so that $p_{\norm}(\ev)=\sigma_1(\Pt(\ev)),$
$\Pt(\ev)=F(\ev)^{-1}P(\ev)F(\ev),$ is of the form in Theorem~\ref{thm:model}, with $r_{\effr} \equiv 0$,
with the sum over $J''$ replaced by a locally finite sum (the sum is over
$J_K$ over compact subsets $K\subset O$), the $h_j$, etc., in
\eqref{effnonres-span} depending smoothly on $\ev$, \ie they are in
$\Cinf(\scT^*_{\pa X}X\times O)$, vanishing at $\cS$ as in
Theorem~\ref{thm:model}, and with the $h''_{\alpha''\beta''}$ supported in
$\scT^*_{\pa X}X\times O_{\alpha''\beta''}$ in terms of \eqref{HMV2r.201}.
\end{thm}

\begin{rem}
For $P=x^{-1}(\Delta+V-\ev)$ the conditions of the theorem are satisfied
for any bounded $O=I$ disjoint from the discrete set of effectively
resonant $\ev$, since in local coordinates $(y,\mu)$ on $\Sigma(\ev)$, the
eigenspaces of $S$ are independent of $\ev$ as shown in the proof of
Lemma~\ref{lemma:nondeg}, and the $r''_j$ are bounded below by
Remark~\ref{rem:H-lin-ev}.
\end{rem}

\begin{proof}
Since the invariant subspaces depend smoothly on $\ev$ by assumption, so do the
eigenvalues of the linearization, and there is smooth
family of local contact diffeomorphisms, \ie coordinate changes, under
which $p(\ev)$ takes the form
\eqref{HMV2r.189}, \ie
\begin{equation}
p(\ev)
=\lambda(\ev)\big(-\nu+\sum_{j=1}^{m-1} r_j(\ev) y_j\mu_j+\sum_{j=m}^{n-1}Q_j(\ev,y_j,\mu_j)+\nu g_1+g_2\big)
\label{HMV2r.189-2}\end{equation}
the $Q_j(\ev,.),$ are homogeneous polynomials of degree 2,
$g_1$ vanishes at least
linearly and $g_2$ to third order, all depending smoothly on $\ev$.

For the rest of the argument it is convenient to reduce the size of the
parameter set $O$ as follows.
For $\ev\in O$, let
\begin{equation}\begin{split}
\hat O(\ev)=&\left(\bigcap_{(\alpha'',\beta''):\ \fraks_{\alpha'',\beta''}(\ev)\in(1,2)}
\fraks_{\alpha'',\beta''}^{-1}((1,2))\right)\\
&\cap 
\left(\bigcap_{(\alpha'',\beta''):\ \fraks_{\alpha'',\beta''}(\ev)\in(-\infty,1)}
\fraks_{\alpha'',\beta''}^{-1}((-\infty,1))\right),
\label{eq:hat-O-def}\end{split}\end{equation}
an open set (as it is a finite intersection of open sets) that includes $\ev$.
Thus, $\{\hat O(\ev):\ \ev\in O\}$ is an open cover of $O$.
Take a locally finite subcover and a partition of unity subordinate to it.
It suffices now to show the theorem for
each element $\hat O(\ev_0)$ of the subcover in place of $O$,
for we can then paste together
the models $p_{\norm}$ we thus obtain using the partition of unity.
Thus, we may assume that $O=\hat O(\ev_0)$ for some $\ev_0\in O$,
and prove the theorem with the sum over $J''$ replaced by a sum over
$J''(\ev_0)$.
Hence, on $O$, for any $(\alpha'',\beta'')$ either 
\begin{enumerate}
\item[a)]
$\fraks_{\alpha''\beta''}(\ev_0)>1$, and then for some $(\tilde\alpha'',
\tilde\beta'')\in J''(\ev_0)$, $(\alpha'',\beta'')\geq(\tilde\alpha'',
\tilde\beta'')$ (reduce $|\alpha''|+|\beta''|$ until
$\fraks_{\tilde\alpha'',\tilde\beta''}\in(1,2)$ -- this will happen
as $r_j\in(0,1/2)$) hence $\fraks_{\alpha''\beta''}(\ev)\geq
\fraks_{\tilde\alpha''\tilde\beta''}(\ev)>1$ for all $\ev\in O$ by
the definition of $\hat O(\ev_0)$, or
\item[b)]
$\fraks_{\alpha''\beta''}(\ev_0)<1$, and then $\fraks_{\alpha''\beta''}(\ev)<1$
for all $\ev\in O$ by the definition of $\hat O(\ev_0)$.
\end{enumerate}

In order to make $\Phi(\ev)$ smooth in $\ev$,
we slightly modify the construction of the local contact diffeomorphism
$\Phi_1(\ev)$ in Proposition~\ref{prop:linear}
so that for any given $\ev$ we do not necessarily remove every term we can
(\ie which are non-resonant for that particular $\ev$). Namely,
we choose the set $I'$ of multiindices $(a,\alpha,\beta)$
{\em which we do not
remove by $\Phi_1(\ev)$} so that $I'$ is
independent of $\ev$, and such that $I'$ contains every multiindex which is resonant for \emph{some} $\ev \in O$, \ie $I'\supset\cup_{\ev\in O}I(\ev)$,
with $I(\ev)$ denoting the set of multiindices corresponding
to resonant terms for $p(\ev)$, as in Proposition~\ref{prop:linear}.
With any such choice of $I'$, the local contact diffeomorphism
of Proposition~\ref{prop:linear},
$\Phi_1(\ev)$, can be chosen smoothly in $\ev$
such that $\lambda^{-1}\Phi_1^*p$ is of the form
\begin{equation*}
-\nu+\sum_{j=1}^{m} r_j(\ev) y_j\mu_j+\sum_{j=m+1}^{n-1}Q_j(\ev,y_j,\mu_j)
+\sum_{I'} c_{a\alpha\beta}(\ev)\nu^a e^\alpha f^\beta
\ \text{modulo}\ \cI^\infty=\frakh^\infty\text{ at }\Rp,
\end{equation*}
with $c_{a\alpha\beta}$ depending smoothly on $\ev$.

The requirement $I'\supset\cup_{\ev\in O}I(\ev)$ means
that for $(a,\alpha,\beta)\nin I'$, $R_{a,\alpha,\beta}(\ev)$ must not vanish
for $\ev\in O$. Here we recall that $R_{a,\alpha,\beta}(\ev)$ is the
eigenvalue of $\{\{p_0,.\}\}$ defined by \eqref{HMV2r.194}, namely
\begin{equation}
R_{a,\alpha,\beta}(\ev)=\lambda\left(a-1+\sum_{j=1}^{n-1}\alpha_j
r_j(\ev)+\sum\limits_{j=1}^{n-1}\beta_j(1-r_j(\ev))\right) 
\label{eq:HMV2r.194-2}\end{equation}

Keeping this in mind,
we choose $I'$ by defining its complement $(I')^c$ to consist of
multiindices $(a,\alpha,\beta)$ with $2a+|\alpha|+|\beta|\geq 3$
such that either
\begin{enumerate}
\item\label{case1}
$a+|\beta'|=1$ and $\alpha''=0$, $\alpha'''=0$, $\beta''=0$, $\beta'''=0$, or
\item\label{case2}
$|\alpha'''|\geq 1$, $\beta'''=0$, or
\item\label{case3}
$|\beta'''|\geq 1$, $\alpha'''=0$, or
\item\label{case4}
$a=0$, $\beta'=0$, $|\alpha'''|+|\beta'''|=2$, $\alpha''=0$, $\beta''=0$,
or
\item\label{case5}
$a=0$, $\beta'=0$, $\alpha'''=\beta'''=0$,
$\fraks_{\alpha''\beta''}(\ev)<1$ (for one, hence all,
$\ev\in O$, as remarked above).
\end{enumerate}

We next show that
multiindices in $(I')^c$ are indeed non-resonant.
In cases \eqref{case2}--\eqref{case3}, $\im R_{a,\alpha,\beta}(\ev)\neq 0$
since the imaginary part of all terms in \eqref{eq:HMV2r.194-2} (with
nonzero imaginary part) has the same sign, and there is at least one term
with non-zero imaginary part, so $(a,\alpha,\beta)$ is non-resonant.

In case \eqref{case5}, the non-resonance follows from
\begin{equation*}
\lambda^{-1} R_{a,\alpha,\beta}(\ev)\leq-1+\fraks_{\alpha''\beta''}(\ev)<0,
\end{equation*}
since $\lambda^{-1} R_{a,\alpha,\beta}(\ev)
=-1+\fraks_{\alpha''\beta''}(\ev)+\sum_{j=1}^{s-1}
\alpha_j r_j$, and each term in the last summation is non-positive.

In case \eqref{case1}, if $a=1$, $\beta'=0$ then
$\lambda^{-1}R_{a,\alpha,\beta}(\ev)=\sum_{j=1}^{s-1} r_j\alpha_j<0$ since
$|\alpha'|\geq 1$ due to $2a+|\alpha|+|\beta|\geq 3$. Also in
case \eqref{case1}, if $a=0$, $|\beta'|=1$, with say $\beta_l=1$, then
\begin{equation*}
\lambda^{-1}R_{a,\alpha,\beta}(\ev)=-r_l+\sum_{j=1}^{s-1}\alpha_j r_j
\end{equation*}
which does not vanish since otherwise $(a,\alpha,\beta)$ 
would be effectively resonant -- it would correspond to
one of the terms in the first summation
in \eqref{effres-span}.

Finally, in case \eqref{case4},
\begin{equation*}
\lambda^{-1}\re R_{a,\alpha,\beta}(\ev)=\sum_{j=1}^{s-1}\alpha_j r_j<0
\end{equation*}
since $\alpha'\neq 0$ due to $2a+|\alpha|+|\beta|\geq 3$.

Thus, all terms corresponding to multiindices in $(I')^c$ can be
removed from $p(\ev)$ by a local contact diffeomorphism $\Phi_1(\ev)$ that
is $\Cinf$ in $\ev.$ So we only need to remark that any term
corresponding to a multiindex in $I'$ can be absorbed into $r_{\effnr}(\ev).$
In fact, such a multiindex has either
\begin{enumerate}
\item[1)]
$a+|\beta'|\geq 2$, or
\item[2)]
$a+|\beta'|=1$ and $|\alpha''|+|\alpha'''|+|\beta''|+|\beta'''|\geq 1$, or
\item[3)]
$|\alpha'''|+|\beta'''|\geq 3$ (with neither $\alpha'''$ nor $\beta'''$ zero), or
\item[4)]
$a=0$, $\beta'=0$, $|\alpha'''|=1$, $|\beta'''|=1$,
$|\alpha''|+|\beta''|\geq 1$, or
\item[5)]
$a=0$, $\beta'=0$, $\alpha'''=0$, $\beta'''=0$, $\fraks_{\alpha''\beta''}>1$.
\end{enumerate}
The first two cases can be incorporated into the $h_0$ or $h_j$ terms in
\eqref{effnonres-span}. The third and fourth ones can be incorporated into
the $h'''_{jk}$
term. Finally, in the fifth case, any infinite linear combination of these
monomials can be written as
\begin{equation*}
\sum_{(\tilde\alpha'',\tilde\beta'')\in J''(\ev_0)}
h''_{\tilde\alpha'',\tilde\beta''} (e'')^{\tilde\alpha''}(f'')^{\tilde\beta''},
\end{equation*}
as remarked in (a) after \eqref{eq:hat-O-def}.

We thus obtain
\begin{equation*} 
\lambda(\ev)\big(-\nu+\sum_j r_j(\ev) y_j\mu_j+\sum_{j=m}^{n-1}Q_j(y_j,\mu_j) 
+ r_{\effnr}(\ev) + r_\infty\big),
\end{equation*}
with $ r_{\effnr}$ as in \eqref{effnonres-span},
and $r_{\infty}$ vanishes to infinite order at $(0,0,0)$.
Finally, we can remove the $r_\infty$ term in a neighbourhood of the
origin, smoothly in $\ev,$ using Proposition~\ref{prop:Nelson} as in the
proof of Theorem~\ref{thm:model}, thus completing the proof of this theorem.
\end{proof}


\section{Microlocal solutions}\label{sec:micro-sol}

In \cite[Equation~(0.15)]{HMV1}
microlocally outgoing solutions were defined using the
global function $\nu$ on $\scT^*_{\pa X}X.$ This is increasing along $W$
and plays the role of a time function; microlocally incoming and outgoing
solution are then determined by requiring the wave front set to lie on one
side of a level surface of $\nu.$ In the present study of microlocal
operators, no such global function is available. However there are always
microlocal analogues, denoted here by ${\rho},$ defined in appropriate
neighbourhoods of a critical point.

\begin{lemma}\label{lemma:cO-exists} There is a neighbourhood $\cO_1$ of
$\Rp$ in $\scT^*_{\pa X}X$ and a function ${\rho}\in\Cinf(\cO_1)$ such that
$\cO_1$ contains no radial point of $P$ except $\Rp,$ $\rho(\Rp)=0,$ and
$W{\rho}\geq 0$ on $\Sigma\cap\cO$ with $W{\rho}>0$ on
$\Sigma\cap\cO_1\setminus\{\Rp\}.$
\end{lemma}

\begin{proof}
This follows by considering the linearization of $W.$ Namely, if $P$ is
conjugated to the form \eqref{HMV2r.189}, then for outgoing radial points
$\Rp$ take $\rho=|y'|^2-(|y''|^2+|y'''|^2+|\mu|^2),$ defined in a
coordinate neighbourhood $\cO_0,$ for incoming radial points take its
negative. On $\Sigma,$ $W\rho\geq c(|y|^2+|\mu|^2)+h$ for some $c>0$
and $h\in \cI^3.$  As $(y,\mu)$ form a coordinate system on $\Sigma$ near
$\Rp,$ it follows that $W\rho\geq \frac{c}{2}(|y|^2+|\mu|^2)$ on a
neighbourhood $\cO'$ of $\Rp$ in $\Sigma.$ Now let $\cO_1\subset\cO_0$ be
such that $\cO\cap\Sigma=\cO'.$ Then $W\rho(p)=0,$ $p\in\cO_1,$
implies $p=\Rp,$ so there are indeed no other radial points in $\cO_1,$
finishing the proof.
\end{proof}

\begin{rem}\label{rem:cO-def}
Below it is convenient to replace $\cO_1$ by a smaller neighbourhood
$\cO$ of $\Rp$ with $\overline{\cO}\subset\cO_1,$ so $\rho$ is
defined and increasing on a neighbourhood of $\overline{\cO}.$
\end{rem}

Consider the structure of the dynamics of $W$ in $\cO.$ First, $\rho$ is
increasing (\ie `non-decreasing') along integral curves $\gamma$ of $W,$
and it is strictly increasing unless $\gamma$ reduces to $\Rp.$  Moreover,
$W$ has no non-trivial periodic orbits and

\begin{lemma}\label{lemma:gamma-cO} Let $\cO$ be as in
Remark~\ref{rem:cO-def}. If $\gamma:[0,T)\to \cO$ or $\gamma:[0,+\infty)\to
\cO$ is a maximally forward-extended bicharacteristic, then either $\gamma$
is defined on $[0,+\infty)$ and $\lim_{t\to+\infty}\gamma(t)=\Rp,$ or
$\gamma$ is defined on $[0,T)$ and leaves every compact subset $K$ of
$\cO,$ \ie there is $T_0<T$ such that for $t>T_0$, $\gamma(t)\nin K.$

An analogous conclusion holds for maximally backward-extended
bicharacteristics.
\end{lemma}

\begin{proof} If $\gamma:[0,+\infty)\to\cO$ then $\lim_{t\to+\infty}
\rho(\gamma(t))=\rho_+$ exists by the monotonicity of $\rho,$ and any
sequence $\gamma_k:[0,1]\to\Sigma,$ $\gamma_k(t)=\gamma(t_k+t),$ $t_k\to
+\infty,$ has a uniformly convergent subsequence, which is then an integral
curve $\tilde\gamma$ of $W$ in $\Sigma$ with image in $\overline{\cO},$
hence in $\cO_1$ along which $\rho$ is constant. The only such
bicharacteristic segment is the one with image $\{\Rp\},$ so
$\lim_{t\to+\infty}\gamma(t)=\Rp.$ The claim for $\gamma$ defined on
$[0,T)$ is standard.
\end{proof}

As in \cite{HMV1} we make use of open neighbourhoods of the critical points
which are well-behaved in terms of $W.$

\begin{Def}\label{HMV.92} By a $W$-balanced neighbourhood of a
non-degenerate radial point $\Rp$ we shall mean a neighbourhood, $O,$ of
$\Rp$ in $\scT^*_{\pa X}X$ with $\overline{O}\subset\cO$ (in which $\rho$
is defined) such that $O$ contains no other radial point, meets
$\Sigma(\ev)\cap O$ in a $W$-convex set (that is, each integral curve of $W$
meets $\Sigma(\ev)$ in a single interval, possibly empty) and is such that
the closure of each integral curve of $W$ in $O$ meets
${\rho}={\rho}(\Rp).$
\end{Def}
\noindent The existence of $W$-balanced neighbourhoods follows as in
\cite[Lemma~1.8]{HMV1}.

If $\Rp$ is a radial point for $P$ and $O$ a $W$-balanced neighbourhood of
$\Rp$ we set
\begin{multline}
\tilde E_{\mic,\out}(O,P)=
\big\{u\in\CmI(X);O\cap\WFsc(Pu)=\emptyset,\\
\Mand \WFsc(u)\cap O\subset\{{\rho} \geq {\rho}(\Rp)\}\big\},
\label{HMV.42}\end{multline}
with $\tilde E_{\mic,\inc}(O,P)$ defined by reversing the inequality.

\begin{lemma}\label{lemma:E-mic-structure-8}
If $O\ni\Rp$ is a $W$-balanced neighbourhood then every $u\in \tilde
E_{\mic,\pm}(O,P)$ satisfies $\WFsc(u)\cap O\subset\Phi_\pm(\{\Rp\});$
furthermore, for $u\in \tilde E_{\mic,\pm}(O,P)$
\begin{equation*}
\WFsc(u)\cap O=\emptyset\Longleftrightarrow \Rp\nin\WFsc(u).
\label{HMV.196}\end{equation*}
\end{lemma}

\noindent Thus, we could have defined $\tilde E_{\mic,\pm}(O,P)$ by
strengthening the restriction on the wavefront set to $\WFsc(u)\cap
O\subset\Phi_\pm(\{\Rp\}).$ With such a definition there is no need for $O$
to be $W$-balanced; the only relevant bicharacteristics would be those
contained in $\Phi_\pm(\{\Rp\}).$ Moreover, with this definition ${\rho}$
does not play any role in the definition, so it is clearly independent
of the choice of ${\rho}$.

\begin{proof} For the sake of definiteness consider $u\in \tilde
E_{\mic,\out}(O,P)$; the other case follows similarly. Suppose $\zeta\in
O\setminus\{\Rp\}.$ If ${\rho}(\zeta)<{\rho}(\Rp),$ then
$\zeta\nin\WFsc(u)$ by the definition of $\tilde E_{\mic,\out}(O,P),$ so we
may suppose that ${\rho}(\zeta)\geq{\rho}(\Rp).$ Since $\Rp\in
\Phi_+(\{\Rp\})$ we may also suppose that $\zeta\neq\Rp.$

Let $\gamma:\Real\to\Sigma$ be the bicharacteristic through $\zeta$ with
$\gamma(0)=\zeta.$ As $O$ is $W$-convex, and $\WFsc(Pu)\cap O=\emptyset,$
the analogue here of H\"ormander's theorem on the propagation of
singularities shows that
\begin{equation*}
\zeta\in\WFsc(u)\Rightarrow \gamma(\Real)\cap O\subset\WFsc(u).
\end{equation*}
As $O$ is $W$-balanced, there exists $\zeta'\in
\overline{\gamma(\Real)}\cap O$ such that ${\rho}(\zeta')={\rho}(\Rp).$ If
$\rho(\zeta)=\rho(\Rp)=0,$ we may assume that $\zeta'=\zeta.$ Fr{o}m this
assumption, and the fact that $\rho$ is increasing along the segment of
$\gamma$ in $\cO,$ and $O$ is $W$-convex, we conclude that $\zeta'\in
\overline{\gamma((-\infty,0])}\cap O.$

If $\zeta'=\gamma(t_0)$ for some $t_0\in\Real,$
then for $t<t_0,$ ${\rho}(\gamma(t))<{\rho}(\gamma(t_0))={\rho}(\Rp),$ and for
sufficiently small $|t-t_0|,$ $\gamma(t)\in O$ as $O$ is open. Thus,
$\gamma(t)\nin\WFsc(u)$ by the definition of $\tilde E_{\mic,\out}(O,P),$
and hence we deduce that $\zeta\nin\WFsc(u).$

On the other hand, if $\zeta'\nin\gamma(\Real),$ then as $O$ is open
$\gamma(t_k)\in O$ for a sequence $t_k\to-\infty,$ and as $O$ is
$W$-convex, $\gamma|_{(-\infty,0]}\subset O.$ Then, again from the
propagation of singularities and Lemma~\ref{lemma:gamma-cO}, $\zeta'=\Rp.$
\end{proof}

We may consider $\tilde E_{\mic,\pm}(O,P)$ as a space of microfunctions, $E
_{\mic,\out}(\Rp,P),$ by identifying elements which differ by functions
with wavefront set not meeting $O$:
\begin{equation*}
E_{\mic,\pm}(\Rp,P)=\tilde E_{\mic,\pm}(O,P)/
\{u\in\CmI(X);\WFsc(u)\cap O=\emptyset\}.
\end{equation*}
The result is then independent of the choice of $O,$
as we show presently.

If $O_1$ and $O_2$ are two $W$-balanced neighbourhoods of $\Rp$ then
\begin{equation}
O_1\subset O_2\Longrightarrow
\tilde E_{\mic,\pm}(O_2,P)\subset\tilde E_{\mic,\pm}(O_1,P).
\label{HMV.84}\end{equation}
Since $\{u\in\CmI(X);\WFsc(u)\cap O=\emptyset\}\subset \tilde
E_{\mic,\pm}(O,P)$ for all $O$ and this linear space decreases with $O,$
the inclusions \eqref{HMV.84} induce similar maps on the quotients  
\begin{equation}
\begin{gathered}
E_{\mic,\pm}(O,P)=\tilde E_{\mic,\pm}(O,P)/
\{u\in\CmI(X);\WFsc(u)\cap O=\emptyset\},\\
O_1\subset O_2\Longrightarrow
E_{\mic,\pm}(O_2,P)\longrightarrow E_{\mic,\pm}(O_1,P).
\end{gathered}
\label{HMV.85}\end{equation}

\begin{lemma}\label{HMV.86} Provided $O_i,$ for $i=1,$ $2,$ are $W$-balanced
neighbourhoods of $\Rp,$ the map in \eqref{HMV.85} is an isomorphism.
\end{lemma}

\begin{proof} We work with $E_{\mic,+}$ for the sake of definiteness.

The map in \eqref{HMV.85} is injective since any element $u$ of its kernel
has a representative $\tilde u\in\tilde E_{\mic,+}(O_2,\ev)$ which
satisfies $\Rp\nin\WFsc(\tilde u),$ hence $\WFsc(\tilde u)\cap
O_2=\emptyset$ by Lemma~\ref{lemma:E-mic-structure-8}, so $u=0$ in
$E_{\mic,+}(O_2,\ev).$

The surjectivity follows from H\"ormander's existence theorem in the real
principal type region \cite{MR48:9458}. First, note that
\begin{equation*}
R=\inf\{\rho(p):\ p\in\Phi_+(\{\Rp\})\cap(\cO\setminus O_1)\}>0=\rho(\Rp)
\end{equation*}
since in $\cO,$ $\rho$ is increasing along integral curves of $W,$ and
strictly increasing away from $\Rp.$ Let $U$ be a neighbourhood of
$\Phi_+(\{\Rp\})\cap\overline{O_1}$ such that $\overline{U}\subset\cO,$ and
$\rho>R_0= R/2$ on $U\setminus O_1.$ Let $A\in\Psisc^{-\infty,0}(\cO)$ be
such that $\WFsc'(\Id-A)\cap \overline{O_1}\cap\Phi_+(\{\Rp\})=\emptyset$
and $\WFsc'(A)\subset U.$ Thus, $\WFsc(Au)\subset U$ and $\WFsc(PAu)
\subset U\setminus O_1,$ so in particular $\rho>R_0$ on $\WFsc(PAu).$ We
have thus found an element, namely $\tilde u=Au,$ of the equivalence class
of $u$ with wave front set in $\cO$ and such that $\rho>R_0>0=\rho(\Rp)$ on
the wave front set of the `error', $P\tilde u.$

The forward bicharacteristic segments from $U\setminus
O_1$ inside $\cO$ leave $\overline{O_2}$ by the remark after
Lemma~\ref{lemma:cO-exists}; since $\overline{O_2}\setminus O_1$ is
compact, there is an upper bound $T>0$ for when this happens. Thus,
H\"ormander's existence theorem allows us to solve $Pv=P\tilde u$ on $O_2$
with $\WFsc(v)$ a subset of the forward bicharacteristic segments emanating
from $U\setminus\overline{O_1}.$ Then $u'=\tilde u-v$ satisfies
$\WFsc(u')\subset\cO\cap\{\rho\geq 0=\rho(\Rp)\},$ $\WFsc(Pu')\cap
O_2=\emptyset,$ so $u'\in E_{\mic,+}(O_2,P),$ and $\Rp\nin\WFsc(u'-u).$ Thus
$\WFsc(u'-u)\cap O_1=\emptyset,$ hence $u$ and $u'$ are equivalent in
$\tilde E_{\mic,+}(O_1,P).$  This shows surjectivity.
\end{proof}

It follows from this Lemma that the quotient space $E_{\mic,\pm}(\Rp,P)$ in
\eqref{HMV.85} is well-defined, as the notation already indicates, and each
element is determined by the behaviour microlocally `at' $\Rp.$ When $P$ is
the operator $x^{-1}(\Delta+V-\ev),$ then we will denote this space
\begin{equation}
E_{\mic,\pm}(\Rp,\ev). 
\label{Emic-defn}\end{equation}

\begin{Def}\label{HMV2r.204} By a \emph{microlocally outgoing
solution} to $Pu=0$ at a radial point $\Rp$ we shall mean either
an element of $\tilde E_{\mic,\out}(O,P),$ where $O$ is a $W$-balanced neighborhood of $\Rp,$ or of $E_{\mic,\out}(\Rp,P).$
\end{Def}


\section{Test Modules}\label{sec:modules}

Following \cite{HMV1}, we use test modules of pseudodifferential operators
to analyze the regularity of microlocally incoming solutions near radial points.
This involves microlocalizing near the critical point with errors which are
well placed relative to the flow. For readers comparing this discussion
to \cite{HMV1}, we mention that the microlocalizer $Q$ in the following
definition corresponds to the microlocalizer $Q$ in
\cite[Equation~(6.27)]{HMV1}; the orders in the commutator are different
as now we are working with $P\in \Psisc^{*,-1}(X)$.

\begin{Def}\label{HMV2r.202} An element $Q\in\Psisc^{*,0}(X)$ is a
\emph{forward microlocalizer} in a neighbourhood $O\ni\Rp$ of a radial
point $\Rp\in\scT^*_{\pa X}X$ for $P\in\Psisc^{*,-1}(X)$ if it is elliptic
at $\Rp$ and there exist $B,$ $F\in\Psisc^{0,0}(O)$ and $G\in
\Psisc^{0,1}(X)$ such that
\begin{equation}
i[Q^*Q,P]= (B^*B + G)+ F\Mand \WFscp(F)\cap\Phi_+(\{\Rp\})=\emptyset.
\label{eq:[Q,P]}\end{equation}
\end{Def}

Using the normal form established earlier we can show that such forward
microlocalizers exist under our standing assumption that 
\begin{equation}
\begin{gathered}
\text{the linearization
has neither a Hessian threshold subspace, \eqref{Hessian-threshold},}\\
\text{nor any non-decomposable 4-dimensional invariant subspace.}
\end{gathered}
\label{Standingassumption}\end{equation}

\begin{prop}\label{Forward-Microlocalizer}
A forward microlocalizer exists in any neighbourhood of any non-degenerate
outgoing radial point $\Rp\in\scT^*_{\pa X}X$ for $P\in\Psisc^{*,-1}(X)$ at which
the linearization satisfies \eqref{Standingassumption}.
\end{prop}

\begin{proof} Since the conditions \eqref{eq:[Q,P]} are microlocal and
  invariant under conjugation with an elliptic Fourier integral operator,
  it suffices to consider the model form in Theorem~\ref{thm:model} which
  holds under the same conditions \eqref{Standingassumption}.

Let $R=|\mu'|^2+|y''|^2+|y'''|^2+|\mu''|^2+|\mu'''|^2,$ and
\begin{equation*}
S=\{p_{\norm}=0,\ R=0\},
\end{equation*}
so $S$ is the flow-out of $\Rp.$ We shall choose
$Q\in\Psisc^{-\infty,0}(X)$ such that
\begin{equation}
\sigma_{\pa}(Q)=q=\chi_1(|y'|^2)\chi_2(R)\psi(p_{\norm}),
\end{equation}
where $\chi_1,$ $\chi_2,$ $\psi\in\Cinf_c(\Real),$ $\chi_1,$ $\chi_2\geq 0$
are supported near $0,$ $\psi$ is supported near $0,$ $\chi_1,$
$\chi_2\equiv 1$ near $0$ and $\chi_1'\leq 0$ in $[0,\infty).$ Choosing all
  supports sufficiently small ensures that $Q\in\Psisc^{-\infty,0}(O).$
Note that $\supp d(\chi_2\circ R)\cap S=\emptyset.$
On the other hand,
\begin{equation}
\scH_p \chi_1(\sum_j (y'_j)^2)=2 \sum_j y'_j(\scH_p y'_j)\chi_1'(|y'|^2)
=2\lambda y'_j(r'_j y'_j+h_j)\chi_1'(|y'|^2),
\end{equation}
with $h_j$ vanishing quadratically at $\Rp.$ Moreover, on
$\supp\chi_1'\circ (|.|^2),$ $y'$ is bounded away from $0.$ Since $r'_j<0,$
$-\sum_j r'_j (y'_j)^2>0$ on $\supp\chi_1'\circ (|.|^2).$ The error terms
$h_j$ can be estimated in terms of $|y'|^2,$ $R$ and $p_{\norm}^2,$ so,
given any $C>0,$ there exists $\delta>0$ such that the $-\sum_jy'_j(r'_j
y'_j+h_j)>0$ if $\supp\chi_1\subset(-\delta,\delta),$ $R/|y'|^2<C$ and
$|p_{\norm}|/|y'|<C.$ In particular, taking $C=2,$ $-\sum_jy'_j(r'_j
y'_j+h_j)>0$ on $S\cap\supp\chi_1'\circ (|.|^2),$ for $R=p_{\norm}=0$ on $S.$
Thus \eqref{eq:[Q,P]} is satisfied (with $B$ appropriately specified,
microsupported near $S),$ provided that $\chi_1$ is chosen so that
$(-\chi_1\chi_1')^{1/2}$ is smooth.

More explicitly, letting $\chi\in\Cinf_c(\Real)$ be supported in
$(-1,1)$ be identically equal to $1$ in $(-\frac{1}{2},\frac{1}{2})$ with
$\chi'\leq 0$ on $[0,\infty),$ $\chi\geq 0,$
$\chi_1=\chi_2=\psi=\chi(./\delta).$ Indeed, for any choice of $\delta\in(0,1),$
$|y'|^2\geq\delta/2$ on $\supp \chi_1'\circ|.|^2,$ hence $R/|y'|^2<2$,
$|p_{\norm}|/|y'|<2$ on $\supp q\cap\supp \chi_1'\circ|.|^2.$ With $C=2,$ choosing
$\delta\in(0,1)$ as above, we can write
\begin{equation}\begin{gathered}
\sigma_{\pa}(i[Q^*Q,P])=-\scH_p q^2=-4\lambda \tilde b^2+\tilde f,\\
\tilde b=(\sum_j y'_j(r'_j y'_j+h_j)\chi_1'(|y'|^2)\chi_1(|y'|^2))^{1/2}
\chi_2(R)\psi(p_{\norm}),\ \supp\tilde f\cap S=\emptyset,
\end{gathered}\label{eq:Q-saddle}\end{equation}
which finishes the proof since $\lambda<0$ for an outgoing radial point.
\end{proof}

A test module in an open set $O \subset \scT^*_{\pa X}X$ is, by definition,
a linear subspace $\cM \subset \Psisc^{*, -1}(X)$ consisting of operators
microsupported in $O$ which contains and is a module over
$\Psisc^{*,0}(X),$ is closed under commutators, and is algebraically
finitely generated. To deduce regularity results we need extra conditions
relating the module to the operator $P.$

\begin{Def}\label{5.5.2004.1} If $P\in\Psisc^{*,-1}(X)$ has real
  principal symbol near a non-degenerate outgoing radial point $\Rp$ then a
  test module $\cM$ is said to be $P$-positive at $\Rp$ if it is supported
  in a $W$-balanced neighbourhood of $\Rp$ and
\begin{enumerate}
\item $\cM$ is generated by $A_0=\Id,$ $A_1,\dots,$ $A_N=P$ over
  $\Psisc^{*,0}(X),$ 
\item for $1\le i\le N-1,$ $0\le j\le N$ there exists
  $C_{ij}\in\Psi^{*,0}_{\scl}(X),$ such that
\begin{equation}
i[ A_i,xP]=\sum\limits_{j=0}^N xC_{ij} A_j
\label{eq:A_i-comm-p}\end{equation}
where $\sigma_{\pa}(C_{ij})(\tilde\Rp)=0,$ for all $0\not=j<i,$
and $\re\sigma_{\pa}(C_{jj})(\tilde\Rp)\geq 0.$
\end{enumerate}
\end{Def}

As shown in \cite{HMV1}, microlocal regularity of solutions of a
pseudodifferential equation can be deduced by combining such a $P$-positive
test module with a microlocalizing operator as discussed above. We recall
and slightly modify this result.

\begin{prop}\label{prop:pos-comm-condition}(Essentially Proposition 6.7 of
\cite{HMV1}, see Proposition~\ref{prop:HMV1-mod} below for a slightly
modified statement and a corrected
proof).
Suppose that $P\in \Psisc^{*,-1}(X)$ has real
principal symbol, $\Rp$ is a non-degenerate outgoing radial point for $P,$
\begin{equation}\label{eq:xP-formally-sa}
\sigma_{\pa,1}(xP-(xP)^*)(\Rp)=0,
\end{equation}
$\cM$ is a $P$-positive test module at $\Rp,$ $Q,Q'\in\Psisc^{*,0}(X)$ are
forward microlocalizers for $P$ at $\Rp$ with $\WFscp(Q')$ being a subset
of the
elliptic set of $Q$. Finally suppose that
for some $s<-\frac1{2},$ $u\in\Hsc^{\infty,s}(X)$ satisfies
\begin{equation}
\WFsc(u)\cap O\subset\Phi_+(\{\Rp\})\text{ and }Pu\in\CIdot(X).
\label{improving}\end{equation}
Then $u\in\Isc^{(s)}(O',\cM)$ where $O'$ is the elliptic set of $Q'.$
\end{prop}

\begin{proof} As already noted this is essentially Proposition~6.7 of
\cite{HMV1}, with a small change to the statement and the
proof given in Proposition~\ref{prop:HMV1-mod} below. However,
there are some small differences to be noted. In \cite{HMV1} (and here
in the Appendix),
the condition in \eqref{eq:A_i-comm-p} was $j>i$; here we changed to $j
< i$ for a more convenient ordering. Since the labelling is arbitrary, this
does not affect the proof of the Proposition.

Also, in \cite{HMV1} the proposition was stated for the $0$th order
operators such as $\Delta+V-\ev$, which are formally self-adjoint with
respect to a scattering metric. This explains the appearance of $xP$ both
in \eqref{eq:xP-formally-sa} and in \eqref{eq:A_i-comm-p} here, even though
in the applications below, $[A_i,x]$ could be absorbed in the $C_{i0}$
term. In particular, $s<-1/2$ in \eqref{improving} arises from a pairing
argument that uses the formal self-adjointness of $xP$, modulo terms that
can be estimated by $[x^s A^\alpha,xP],$ $s>0,$ $A^\alpha$ a product of the
$A_j.$

The proposition in \cite{HMV1} is proved with \eqref{eq:xP-formally-sa}
replaced by $(xP)=(xP)^*$, but \eqref{eq:xP-formally-sa} is sufficient for
all arguments to go through, since $B=(xP)-(xP)^*$ would contribute error
terms of the form $x^s A^\alpha B$ with $\sigma_{\pa,1}(B)(\Rp)=0$, which
can thus be handled exactly the same way as the $C_{jj}$ term in
\eqref{eq:A_i-comm-p}.

In fact \eqref{eq:xP-formally-sa} can always be arranged for any
$P_0\in\Psisc^{*,-1}(X)$ with a non-degenerate radial point and real
principal symbol. Indeed, we only need to conjugate by $x^k$ giving
\begin{equation*}
P=x^kP_0 x^{-k}, \
k=\frac{-\sigma_{\pa,1}(B)(\Rp)}{2i\lambda}\in\bbR
\end{equation*}
satisfies \eqref{eq:xP-formally-sa}; here $dp|_{\Rp}=\lambda\alpha|_{\Rp}$,
with $\alpha$ the contact form. Microlocal solutions $P_0u_0=0,$
correspond to microlocal solutions $Pu=0$ via $u=x^k u_0,$ so
$u\in\Hsc^{\infty,s}(X)$ is replaced by $u_0\in\Hsc^{\infty,s-k}(X).$
\end{proof}

Thus, iterative regularity with respect to the module essentially reduces
to showing that the positive commutator estimates \eqref{eq:A_i-comm-p}
hold. For each critical point $\Rp$ satisfying \eqref{Standingassumption} a
suitable (essentially maximal) module is constructed below, so microlocally
outgoing solutions to $Pu=0$ have iterative regularity under the module;
that is, that
\begin{equation}
u\in\Isc^{(s)}(O,\cM)=
\{u;\cM^mu\subset\Hsc^{\infty,s}(X) \text{ for all } m \}.
\label{Isc-defn}\end{equation}
The test modules are elliptic off the forward flow out $\Phi_+(\Rp)$ which
is an isotropic submanifold of $\Sigma.$ Thus, it is natural to expect that
$u$ is some sort of an isotropic distribution. In fact the flow out (in the
model setting just the submanifold $S)$ has non-standard homogeneity
structure, so these distributions are more reasonably called
`anisotropic'.

First we construct a test module for the model operator when there are no
resonant terms. Thus, we can assume that the principal symbol is
\begin{equation*}
p_0=\lambda\big(-\nu+\sum_{j=1}^{m-1} r_j y_j\mu_j+\sum_{j=m}^{n-1}Q_j(y_j,\mu_j)\big).
\end{equation*}
Then let $\cM$ be the test module generated by $\Id$ and operators with
principal symbols
\begin{equation}
x^{-1} f'_j, \quad x^{-r''_j} e''_j, \quad x^{-(1-r''_j)} f''_j, \quad
x^{-1/2} e'''_j, \quad x^{-1/2} f'''_j \quad \text{and } \ x^{-1} p_0 
\label{module-defn-nr}\end{equation}
over $\Psisc^{*,0}(X).$ 

Note that the order of the generators is given by the negative of the
normalized eigenvalue (\ie the eigenvalue in Lemma~\ref{HMV2r.191} divided
by $\lambda$) subject to the conditions that if the order would be $<-1,$
it is adjusted to $-1,$ and if it would be $>0,$ it is omitted. The latter
restrictions conform to our definition of a test module, in which all terms
of order $0$ are included and there are no terms of order less than $-1.$ These
orders can be seen to be optimal (\ie most negative) by a principal symbol
calculation) of the commutator with $A$ in which the corresponding
eigenvalue arises.

\begin{lemma}\label{nonres-commutators}
Suppose $P$ is nonresonant at $\Rp.$ Then the module $\cM$ generated
by \eqref{module-defn-nr} is closed under commutators and satisfies
condition \eqref{eq:A_i-comm-p}.
\end{lemma}

\begin{proof} It suffices to check the commutators of generators to show
that $\cM$ is closed.  In view of \eqref{eq:scHp-local} (applied with $a$
in place of $p$), $\{a,b\}=\scH_a b$, this can be easily done.  Property
\eqref{eq:A_i-comm-p} follows readily from \eqref{bracket-defn}. Indeed, we
have the stronger property
\begin{equation*}
i [A_i, P(\ev)] = c_i A_i + G_i,\ G_i \in \Psi^{*,0}(X),\ \re c_i\geq 0
\end{equation*}
where $A_i$ is any of the generators of $\cM$ listed in \eqref{module-defn-nr}. 
\end{proof}

\begin{rem} We may take generators of $\cM$ to be the operators
\begin{equation}\begin{gathered}
D_{y_j'}, \ x^{-r''_j} y''_j, \ x^{r''_j} D_{y_j''}, \ x^{-1/2}
y'''_j, \ x^{1/2} D_{y'''_j}\text{ and} \\  
xD_x+\sum_{j=1}^{m-1} r_j y_j D_{y_j} +\sum_{j=m}^{n-1}Q_j(x^{-1/2}y_j,
x^{1/2}D_{y_j}).
\end{gathered}\label{module-defn-nr-ops}\end{equation}
\end{rem}

Combining this with Proposition~\ref{prop:pos-comm-condition} proves that,
in the nonresonant case, if $u$ is a microlocal solution at $\Rp,$ and if
$\WFsc^s(u)$ is a subset of the $W$-flowout of $\Rp,$ then
$u\in\Isc^{(s)}(O,\cM)$ for all $s<-1/2.$

The discrepancy between the `resonance order' of polynomials in $\nu^a
e^\beta f^\gamma$, given by $a + \sum_j \beta_j r_j + \sum_k \gamma_k (1 -
r_k)$ and the `module order' given by the sum of the orders of the
corresponding module elements is closely related to arguments which allow
us to regard
most resonant terms as `effectively nonresonant'. To give an explicit
example, take a resonant term of the form $y'_i \mu_j' (y'')^{\beta''},$
corresponding to a term like $x^{-1} y'_i (y'')^{\beta''} (x D_{y'_j})$ in
$P.$ Resonance requires that $r'_i + (1 - r'_j) + \sum_k \beta''_k r''_k =
1$ and $|\beta''| > 0.$ In the module, this corresponds to a product of
module elements with an additional factor of $x^\epsilon$ \emph{with}
$\epsilon > 0,$ since we can write it
$$
x^\epsilon y'_i \prod_k (x^{-r''_k}y''_k)^{\beta''_k} D_{y'_i}, \quad
\epsilon = \sum_k \beta''_k r''_k > 0.
$$
Since, by Proposition~\ref{prop:pos-comm-condition}, the eigenfunction $u$
remains in $x^s L^2(X)$, for all $s < -1/2$, under application of products
of elements of $\mathcal{M},$ this term applied to $u$ yields a factor
$x^\epsilon$, and therefore it can be treated as an error term in
determining the asymptotic expansion of $u;$ see the proof of
Theorem~\ref{thm:enr-structure}. Only the terms with the module order
equal to the resonance order affect the expansion of $u$ to leading order,
and it is these we have labelled `effectively resonant'.

Next we consider the general resonant case. To do so, we need to enlarge
the module $\cM$ so that certain products of the generators of $\cM,$ such
as those in the resonant terms of Theorem~\ref{thm:model}, are also
included in the larger module $\tcM.$ For a simple example, see section 8
of Part I.  It is convenient to replace $P_0$ by $xD_x$ as the last
generator of $\cM$ listed in \eqref{module-defn-nr-ops}, though this is not
necessary; all arguments below can be easily modified if this is not done.
Let us denote the generators of $\cM$ by
\begin{equation}\begin{split}\label{eq:effr-gen}
&A_0 = \Id, A_1=x^{-s_1}B_1, \ldots, A_{N-1}=x^{-s_{N-1}}B_{N-1},
\ A_N = xD_x=x^{-1}B_N,\\
&\qquad\ s_i = - \operatorname{order}(A_i),\ B_i\in\Psisc^{-\infty,0}(O).
\end{split}\end{equation}
Note that for each $i=1,\ldots,N,$ $d\sigma_{\pa,0}(B_i)$ is an eigenvector
of the linearization of $W$; we denote the eigenvalue by $\sigma_i.$ Thus,
\begin{equation*}
s_i=\min(1,\sigma_i)>0\ \text{for}\ i=1,\ldots,N.
\end{equation*}
For any multiindex $\alpha\in\Nat^N$ (with $\Nat=\{1,2,\ldots\}$) let
\begin{equation*}
s(\alpha)=\min(\sum s_i\alpha_i,1),
\ \tilde s(\alpha)=\sum_i s_i\alpha_i-s(\alpha)=\max(0,\sum_i s_i\alpha_i-1),
\end{equation*}
and let
\begin{equation*}
A^\alpha = A_1^{\alpha_1} A_2^{\alpha_2} \dots A_N^{\alpha_N}.
\end{equation*}
Let $e_i$ be the multiindex $e_i=(0,\ldots,0,1,0,\ldots,0),$ where the $1$
is in the $i$th slot, if $i=1,\ldots,N,$ and let $e_0=(0,\ldots,0).$

To deal with resonant terms, we define a module $\cM_k$ generated (over
$\Psisc^{-\infty,0}(O)$) by the operators
\begin{equation}
x^{\tilde s(\alpha)} A^\alpha\in\Psisc^{-\infty,-s(\alpha)}(O),
\ |\alpha|\leq k.
\label{Aalpha}
\end{equation}
Note that $\alpha=0$ gives $\Id$ as one of the generators. Thus, the order
of the generators in \eqref{Aalpha} is `truncated' so that it is always
between $0$ and $-1$; in particular
$\cM_k\subset\Psisc^{-\infty,-1}(O).$ Since in computations below we will
think of $\Psisc^{-\infty,0}(O)$ as the submodule of $\cM_k$ consisting 
of trivial elements, it is convenient to work modulo such terms, so we
use what is essentially the principal symbol equivalence relation on $\cM_k$
where $P\sim Q$ if $P-Q\in\Psisc^{-\infty,0}(O).$

While it appears that the ordering in the factors in the product $A^\alpha$
matters, this is not the case. Indeed, if $\sigma$ is a permutation
of $\{1,\ldots,|\alpha|\}$, and $j:\{1,\ldots,|\alpha|\}\mapsto \{1,\ldots,N\}$
which takes $\alpha_m$-times the value $m$, $m=1,\ldots,N,$ then
\begin{equation*}
x^{\tilde s(\alpha)} A_{j(1)}\ldots A_{j(|\alpha|)}
\sim x^{\tilde s(\alpha)} A_{j(\sigma(1))}\ldots A_{j(\sigma(|\alpha|))},
\end{equation*}
for this is clear if $\sigma$ interchanges $n$ and $n+1,$ as
\begin{equation*}\begin{split}
&x^{\tilde s(\alpha)}A_{j(1)}\ldots A_{j(n-1)}[A_{j(n)},A_{j(n+1)}]
A_{j(n+2)}\ldots A_{j(|\alpha|)}\\
&\qquad\qquad\in\Psisc^{-\infty,
\tilde s(\alpha)+1-\sum s_i\alpha_i}(O)\subset \Psisc^{-\infty,0}(O)
\end{split}\end{equation*}
since $\tilde s(\alpha)+1-\sum_i s_i\alpha_i=1-s(\alpha)\geq 0.$

In addition, for $Q\in\Psisc^{-\infty,0}(O)$,
\begin{equation*}
x^{\tilde s(\alpha)} QA_{j(1)}\ldots A_{j(|\alpha|)}
\sim x^{\tilde s(\alpha)} A_{j(1)}\ldots A_{j(m)}QA_{j(m+1)}\ldots
A_{j(|\alpha|)}.
\end{equation*}
Similarly, one can shift powers of $x$ from in front of the
product to in between factors, so in fact the generators can be written
equivalently, modulo $\Psisc^{-\infty,0}(O)$, as
\begin{equation}\label{Balpha}
x^{s(\alpha)} B^\alpha\in\Psisc^{-\infty,-s(\alpha)}(O),
\ |\alpha|\leq k,
\end{equation}
where $B^\alpha=B_1^{\alpha_1}\ldots B_N^{\alpha_N}$.

Moreover,
there is an integer $J$ such that $\cM_k = \cM_J$ if $k \geq J$; indeed this
is true for any $J \geq 2(r''_s)^{-1}$, where $r''_s$ is the smallest
positive eigenvalue of the operator in Lemma~\ref{lemma:nondeg}
(or $J\geq 4$ if no eigenvalue lies in $(0,\frac{1}{2}]$), since then
adding new elements to the product simply has the effect of multiplying by an
element of $\Psisc^{*,0}(X)$.

In particular, note that the generators in \eqref{Aalpha} or \eqref{Balpha}
are usually not linearly independent: some $B_{\alpha_j}$ may be absorbable
into a $\Psisc^{*,0}(O)$ factor without affecting $s(\alpha)$. We could
easily give a linearly independent (over $\Psisc^{*,0}(O)$) subset of the
generators, but this is of no importance here.

Suppose that $\tilde P$, the normal operator for $P(\ev)$ at $\Rp$,
contains resonant terms. Then Lemma~\ref{nonres-commutators} is replaced by 

\begin{lemma}\label{res-commutators}
Let $>$ be a total order on multiindices $\alpha$ satisfying
\begin{enumerate}
\item
$|\alpha'|>|\alpha|$ implies $\alpha'>\alpha$;
\item
$|\alpha'|=|\alpha|$ and $\sum_k s_k\alpha'_k>\sum_k s_k\alpha_k$
imply $\alpha'>\alpha$;
\item
$|\alpha'|=|\alpha|=1$, $\alpha'=e_i$, $\alpha=e_j$, $s_i=s_j=1$,
$\sigma_i>\sigma_j$ imply that $\alpha'>\alpha$.
\end{enumerate}
With the corresponding ordering of the generators $x^{-\tilde s(\alpha)}
A^\alpha$,
the module $\cM_J$ is a test module for $\tilde P$ at $\Rp$ satisfying
\eqref{eq:A_i-comm-p}.
\end{lemma}

\begin{rem}
(ii) and (iii) could be replaced by (ii)': $|\alpha'|=|\alpha|$ and
$\sum_k \sigma_k\alpha'_k>\sum_k \sigma_k\alpha_k$
imply $\alpha'>\alpha$, which would simplify the statement of the lemma.
However, the proof is slightly simpler with the present statement. Note
that (ii)+(iii) is not equivalent to (ii)', \ie the ordering of
the generators may be different, but either ordering gives
\eqref{eq:A_i-comm-p}.
\end{rem}

\begin{proof}
We first observe that $\cM_J$ is closed under commutators. Indeed,
not only is $\cM$ closed under commutators, but the
commutators $[A_i,A_j]$ can be written as $\sum_{l=0}^N C_l A_l$ with $C_l\in
\Psisc^{-\infty,0}(X)$ and $C_l=0$ unless $s_l\leq s_i+s_j-1$.
Expanding
\begin{equation*}
[x^{\tilde s(\alpha)}Q_\alpha A^\alpha,x^{\tilde s(\beta)}Q_\beta A^\beta],
\ Q_\alpha, Q_\beta\in\Psisc^{-\infty,0}(O),
\end{equation*}
and ignoring momentarily the commutators with powers of $x$ and with $Q_\alpha$
and $Q_\beta$, gives a sum of
terms of the form
\begin{equation*}
x^{\tilde s(\alpha)+\tilde s(\beta)}Q_\alpha Q_\beta
A^{\alpha'}A^{\beta'}[A_i,A_j]A^{\alpha''}
A^{\beta''}
\end{equation*}
with $\alpha=\alpha'+\alpha''+e_i$, and similarly for $\beta$. Substituting in
$[A_i,A_j]=\sum_{l=0}^N C_l A_l$ shows that this is an element of the
module and is indeed equivalent, modulo $\Psisc^{-\infty,0}(O)$, to
\begin{equation}\begin{split}\label{eq:comm-cl-8}
\sum_{l:s_l\leq s_i+s_j-1}
&\left(C_l x^{\tilde s(\alpha)+\tilde s(\beta)-\tilde s(\gamma^{(l)})}\right)
x^{\tilde s(\gamma^{(l)})} A^{\gamma^{(l)}},\\
&\qquad \gamma^{(l)}=\alpha'+\alpha''+\beta'+\beta''+e_l=\alpha+\beta-e_i-e_j+e_l,
\end{split}\end{equation}
provided that
\begin{equation}\label{eq:ts-ineq}
\tilde s(\gamma^{(l)})\leq \tilde s(\alpha)+\tilde s(\beta).
\end{equation}
But $\tilde s(\alpha)+\tilde s(\beta)\geq (\sum s_i\alpha_i-1)
+(\sum s_i\beta_i-1)=\sum s_i\gamma^{(l)}_i+s_i+s_j-s_l-2\geq
\sum s_i\gamma^{(l)}_i-1$ as $s_i+s_j-s_l\geq 1$.
Moreover, $\tilde s(\alpha)+\tilde s(\beta)\geq 0$, so
\begin{equation*}
\tilde s(\alpha)+\tilde s(\beta)\geq\max(\sum s_k\gamma^{(l)}_k-1,0)=\tilde
s(\gamma^{(l)}),
\end{equation*}
proving \eqref{eq:ts-ineq}.

The commutators
\begin{equation}\label{eq:comm-cl-16}
x^{\tilde s(\beta)}Q_\beta[x^{\tilde s(\alpha)}Q_\alpha,A^\beta]A^\alpha,
\ x^{\tilde s(\alpha)}Q_\alpha[A^\alpha,x^{\tilde s(\beta)}Q_\beta]A^\beta
\end{equation}
also lie in $\cM_J$. Indeed, $[A_i,x^\rho Q]=x^{\rho-s_i+1}Q'$ for some
$Q'\in\Psisc^{-\infty,0}(O)$, so they are sums of terms of the form
$x^{\tilde s(\alpha)+\tilde s(\beta)-s_i+1}Q'A^\gamma$ with
$\gamma=\alpha+\beta-e_i$. Now,
\begin{equation*}
\tilde s(\gamma)\leq \tilde s(\alpha)+\tilde s(\beta)-s_i+1
\end{equation*}
since $\tilde s(\alpha)+\tilde s(\beta)-s_i+1\geq 0$ as $1\geq s_i$
as well as $\tilde s(\alpha)+\tilde s(\beta)-s_i+1\geq(\sum_k s_k\alpha_k-1)
+(\sum_k s_k\beta_k-1)-s_i+1=\sum_k s_k\gamma_k-1$, so
$\tilde s(\alpha)+\tilde s(\beta)-s_i+1\geq\max(\sum_k s_k\gamma_k-1,0)
=\tilde s(\gamma)$ indeed, proving that \eqref{eq:comm-cl-16} is in $\cM_J$.
The commutators
\begin{equation}\label{eq:comm-cl-24}
[x^{\tilde s(\alpha)}Q_\alpha,x^{\tilde s(\beta)}Q_\beta]A^\alpha A^\beta
\end{equation}
can be shown to lie in $\cM_J$ by a similar argument,
this time using $\gamma=\alpha+\beta$, and $\tilde s(\gamma)\leq
\tilde s(\alpha)+\tilde s(\beta)+1$.
Thus, we conclude that
$[x^{\tilde s(\alpha)}Q_\alpha A^\alpha,x^{\tilde s(\beta)}Q_\beta
A^\beta]\in\cM_J$, and hence
$\cM_J=\cM_{J+1}=\ldots$
is closed under commutators.

Modulo $\Psisc^{-\infty,0}(O)$,
$x^{\tilde s(\gamma^{(l)}))} A^{\gamma^{(l)}}$ may be replaced by
$x^{-s(\gamma^{(l)})}B^{\gamma^{(l)}}$.
If $|\gamma^{(l)}|>J$ in \eqref{eq:comm-cl-8},
then this is written in terms of one of the
generators listed in \eqref{Balpha} (or equivalently,
modulo $\Psisc^{-\infty,0}(O)$, in \eqref{Aalpha}),
only after some of the factors in $B^{\gamma^{(l)}}$, which we
may always take from $B_l B^{\beta'}B^{\beta''}$, are moved to the front
and are incorporated in $C_l$, \ie they are simply regarded
as $0$th order operators and $C_l$ is replaced
by $\tilde C_l=C_l B_l B^{\beta'}B^{\beta''}$. Notice the principal
symbol of $\tilde C_l$ always vanishes at $\Rp$ in this case.
Analogous conclusions hold for the terms in \eqref{eq:comm-cl-16}
and \eqref{eq:comm-cl-24}.

On the other hand, if $|\gamma^{(l)}|\leq J$, then
$x^{\tilde s(\gamma^{(l)})}A^{\gamma^{(l)}}$ is one of the generators in
\eqref{Aalpha}, and $|\gamma^{(l)}|=|\alpha|+|\beta|-1$ if $l\geq 1$, and
$|\gamma^{(l)}|=|\alpha|+|\beta|-2$ if $l=0$.
Moreover, if $\sum_k s_k\beta_k>1$ then
\begin{equation}\label{eq:comp-gamma-alpha}
\sum_k s_k\gamma^{(l)}_k=\sum_k s_k\alpha_k+\sum_k s_k\beta_k-s_i-s_j+s_l
\geq \sum_k s_k\alpha_k+\sum_k s_k\beta_k-1> \sum_k s_k\alpha_k.
\end{equation}
For the terms in \eqref{eq:comm-cl-16}
and \eqref{eq:comm-cl-24}, if $|\gamma|\leq J$,
we always get $|\gamma|\geq |\alpha|+|\beta|-1$ since $\gamma=\alpha+\beta$
or $\gamma=\alpha+\beta-e_i$ for some $i$.

Now we turn to \eqref{eq:A_i-comm-p}. First, with $\tilde P$ replaced
by $P_0$, \eqref{eq:A_i-comm-p} is certainly satisfied, exactly
as in the non-resonant case, since the $\sigma_{\pa,0}(B^\alpha)$ are
eigenvectors of
the linearization of $W$ with eigenvalue given in Section~\ref{sec:normal}.
Thus,
\begin{equation}\label{eq:comm-cl-32}
i[A^\alpha,x^{-1}P_0]\sum_\gamma C'_\gamma A^\gamma,
\ C'_\gamma\in\Psisc^{-\infty,0}(O),
\end{equation}
with $\sigma_{\pa,0}(C'_\gamma(\Rp))=0$ if $\alpha\neq\gamma$ and
$\re\sigma_{\pa,0}(C'_\alpha(\Rp))\geq 0$.
So it remains to show that it also holds
for the resonant terms. If $x^{-s(\beta)}Q_\beta B^\beta$ is a resonant term,
then
$s(\beta)=1$. Moreover,
\begin{enumerate}
\item
if $|\beta|=1$, then $x^{-1}Q_\beta B^\beta=\sum_{\mu'} (y')^{\mu'} D_{y'_k}$
for some $\mu'$ and some $k$; in particular it is a summand of $r_{\effr}$;
\item
if $|\beta|=2$, then either
$x^{-1}Q_\beta B^\beta= B_j D_{y'_k}$ for some $j>0,k$, or
$x^{-1}Q_\beta B^\beta$ is associated to the sum over $J''$ in
\eqref{effnonres-span}; in either case $\sum s_k\beta_k>1$.
\end{enumerate}

We claim that for a resonant term $x^{-s(\beta)}Q_\beta B^\beta$,
\begin{equation}\label{eq:res-comm-8}
[x^{-s(\alpha)}B^\alpha,x^{-s(\beta)}Q_\beta B^\beta]\sim\sum_\gamma
\tilde C_\gamma x^{-s(\gamma)}B^\gamma,
\ \tilde C_\gamma\in\Psisc^{-\infty,0}(X),
\end{equation}
and each term on the right hand side has the following property:

\begin{enumerate}
\item
Either $\sigma_{\pa,0}(\tilde C_\gamma)(\Rp)=0$,
\item
or $|\gamma|>|\alpha|$,
\item
or $|\gamma|=|\alpha|$, $\sum_k s_k\gamma_k>\sum_k s_k\alpha_k$,
\item
or $|\gamma|=|\alpha|=1$, $\gamma=e_k$, $\alpha=e_j$, $s_j=s_k=1$ and
$\sigma_k>\sigma_j$.
\end{enumerate}

Indeed, if $|\beta|\geq 3$, then either (i) or (ii) holds, depending
on whether any factors $A_k$ had to be cancelled to write the commutator
in terms of the generators in \eqref{Aalpha}. If $|\beta|=2$, then
$\sum s_k\beta_k>1$. Thus,
again, either (i) or (ii) holds, or $|\gamma|=|\alpha|$ and
$\sum_k s_k\gamma_k>\sum_k s_k\alpha_k$ by \eqref{eq:comp-gamma-alpha},
so (iii) holds. Finally, if $|\beta|=1$, then
$x^{-1}Q_\beta B^\beta=\sum_{\mu'} (y')^{\mu'} D_{y'_k}$
for some $\mu'$ and some $k$. Since
$r_1\leq r_2\leq\ldots\leq r_{s-1}<0$, and the resonance condition is
$\sum_{l=1}^{s-1}\mu'_l r_l+(1-r_k)=1$ with $|\mu'|+1\geq 3$,
we immediately deduce that $\mu'_l=0$ for $l\leq k$.
Thus, not only do powers of $x$
commute with $x^{-1}Q_\beta B^\beta$, but all $A_i$ commute with $D_{y'_k}$
and $[A_i,(y')^{\mu'}]=0$ unless $A_i=D_{y'_j}$ and $\mu'_j\neq 0$ for some
$j$, which in turn implies that $j>k$, so $1-r_k>1-r_j$, hence (iv) holds.
This completes the proof
of \eqref{eq:res-comm-8}.

By the assumption on the ordering of the multiindices $\alpha$, we deduce
that for all resonant terms $x^{-s(\beta)}B^\beta$,
\begin{equation*}
i[A^\alpha,x^{-s(\beta)}B^\beta]=\sum_\gamma C_\gamma A^\gamma,
\ C_\gamma\in\Psisc^{-\infty,0}(O),
\end{equation*}
and either $\sigma_{\pa,0}(C_\gamma)(\Rp)=0$, or $\gamma>\alpha$.
Combining this with \eqref{eq:comm-cl-32}, we deduce that
$\cM_J$ satisfies \eqref{eq:A_i-comm-p}.
This establishes the lemma.
\end{proof}

\begin{cor}\label{cor:Isc-reg}
Let $\cM=\cM_J$ be as in the previous lemma.
Suppose that 
\begin{equation}
s < -\frac1{2}, \quad u\in\Hsc^{\infty,s}(X),
\ \tilde P u \in \CIdot(X)\Mand \WFsc(u)\cap O\subset\Phi_+(\{\Rp\}).
\end{equation}
Then $u\in\Isc^{(s)}(O,\cM)$.
\end{cor}

Regularity with respect to $\cM$ can be understood more geometrically
as follows. Suppose $\delta>0$ is sufficiently small so that $(x,y',y'',y'')$
define local coordinates on the region $U$ given by
$0\leq x<\delta$, $|y_j|<\delta$ for all $j$.
Let
\begin{equation}\label{eq:def-Phi}
\Phi:U^\circ\to\Real^n_+,\ 
\Phi(x,y',y'',y''')=(x,y',Y'',Y'''),\ 
Y''_j=\frac{y''_j}{x^{r_j}},\ Y'''_j=\frac{y'''_j}{x^{1/2}}.
\end{equation}
Thus, $\Phi$ is a diffeomorphism onto its range $\Phi(U^\circ)$ with
\begin{equation*}
\Phi^{-1}(x,y',Y'',Y''')=(x,y'_j,x^{r_j}Y''_j,x^{1/2}Y''').
\end{equation*}
Note that $\overline{\Phi(U^\circ)}$ is not compact; $Y''$ and $Y'''$ are
`global' variables.
Thus $\Phi^{-1}$ is actually continuous on $\overline{\Phi(U^\circ)}$
since $r''_j>0$. Thus, $\Phi$ is a blow-up and $\Phi^{-1}$ is a somewhat
singular blow-down map.
In the coordinates $(x,y',Y'',Y''')$ the Riemannian density
takes the form
\begin{equation*}
a x^{-n-1}\,dx\,dy=a x^{-n+\sum r''_j+(n-m)/2-1}\,dx\,dy'\,dY''\,dY''',
\end{equation*}
$a>0$, $a\in\Cinf(X)$. We thus conclude that (for $O$ small)
$u\in \Isc^{(s)}(O,\cM)$ if and only if for any
$Q\in\Psisc^{-\infty,0}(O)$ with Schwartz kernel supported in
$U\times U$, its microlocalization $Qu$ satisfies
\begin{equation}\begin{split}\label{eq:conormal-reg}
(Y'')^{\gamma''}(Y''')^{\gamma'''}
&(xD_x)^a D_{y'}^{\beta'}D_{Y''}^{\beta''}D_{Y'''}^{\beta'''}Qu\\
&\in
x^{s+n/2-\sum r''_j/2-(n-m)/4}L^2(x^{-1}\,dx\,dy'\,dY''\,dY'''),
\end{split}\end{equation}
for every $a$, $\beta$, $\gamma''$ and $\gamma'''$,
\ie if and only if microlocally $u$ is conormal in $(x,y')$ with values in
Schwartz functions in $(Y'',Y''')$, with the weight given by
$s+n/2-\sum r''_j/2-(n-m)/4$.

We also recall that for conormal functions, the $L^2$ and the $L^\infty$
spaces are very close, namely they are included in each other with a loss
of $x^\ep$. Thus, $u\in \Isc^{(s)}(O,\cM)$
implies that
\begin{equation*}\begin{split}
(Y'')^{\gamma''}(Y''')^{\gamma'''}
&(xD_x)^a D_{y'}^{\beta'}D_{Y''}^{\beta''}D_{Y'''}^{\beta'''}Qu\\
&\in
x^{s+n/2-\sum r''_j/2-(n-m)/4-\ep}L^\infty(x^{-1}\,dx\,dy'\,dY''\,dY'''),
\end{split}\end{equation*}
for every $\ep>0$.


\section{Effectively nonresonant operators}\label{sec:enr}

We now assume that the normal form $p_{\norm}$ for $\sigma_1(P(\ev))$ at
$\Rp$ is such that the term $r_{\effr}$ in Theorem~\ref{thm:model}
vanishes. If this is true, we shall call $p_{\norm}$ \emph{effectively
nonresonant}, and $\ev$ an \emph{effectively nonresonant energy} for $\Rp$.
The significance of the notion of effective resonance in general is that
the form of the asymptotics of microlocally outgoing solutions of $Pu=f$,
$f\in\dCinf(X)$, is independent of $r_{\effnr}$; only $r_{\effr}$ changes
this form slightly.  Moreover, effective non-resonance is a more typical
condition than non-resonance.  We deal with the effectively nonresonant
case in this section and treat the effectively resonant case in the
following section.  In both cases, it is convenient to reduce $P,$ and not
only its principal symbol, to model form. This is accomplished in the
following lemma. We recall here our ongoing assumption \eqref{Standingassumption}.  

\begin{lemma}\label{effnonres-form}
Let $p_{\norm}$ be as in Theorem~\ref{thm:model} and $\tilde P$ as in
Remark~\ref{rem:model-op}, i.e.\ $\sigma_{\pa,-1}(\tilde P)=p_{\norm}$. Then
$\tilde P$ can be conjugated by a smooth function to the form
\begin{equation}\begin{split}\label{eq:Pt-form}
P_{\norm}=&\lambda\left(xD_x+\sum_{j=1}^{m-1} r_j y_j D_{y_j} +\sum_{j=m}^{n-1}Q_j(x^{-1/2}y_j, x^{1/2}D_{y_j}) + R_{\effr}+ b + R\right)\\
&R_{\effr} = \sum_{j=1}^{s-1} \cP_{j}(y')
D_{y_j} + \sum_{j=s}^{m-1} \cP_j(y'') D_{y_j}+\cP_0(y''),
\end{split}\end{equation}
where $b$ is a constant, $Q_j$ is a real elliptic homogeneous quadratic
polynomial (\ie a harmonic oscillator), $\cP_j$ and $\cP_0$
are homogeneous polynomials of
degree $r_j$ resp.\ $1$, when $y_k$ is assigned degree $r_k$,
and $R \in x^{\epsilon} (\cM)^j$ for some $j \in \NN$ and $\epsilon > 0$.
In addition, for $s\leq j\leq m-1$, $\cP_j$ is
actually a polynomial in $y_s,\ldots,y_{j-1}$ (\ie is independent
of $y_j,\ldots,y_{m-1}$) without constant or linear terms,
while for $j\leq s-1$,
$\cP_j$ is a polynomial in $y_{j+1},\ldots,y_{s-1}$.

We call $P_{\norm}$ a {\em normal form} for $P.$ If $p_{\norm}$ is effectively
non-resonant then $R_{\effr}=0$.
\end{lemma}

\begin{rem}\label{rem:Q_j-quantization}
Note that $Q_j(x^{-1/2}y_j,x^{1/2}D_{y_j})$ is not completely well-defined
since $Q_j$ is a homogeneous quadratic polynomial, and
$y_j$ and $D_{y_j}$ do not commute. However, any two choices
for the quantization $Q_j$ differ by a constant multiple of the
commutator $[x^{-1/2}y_j,x^{1/2}D_{y_j}]=[y_j,D_{y_j}]$, hence by
a constant.

In particular, with the notation of the previous section,
$Q_j(Y_j,D_{Y_j})$ may be arranged to be self-adjoint with respect to
$dY_j$, by symmetrizing if necessary, which changes $Q_j$ at most by
a constant.
\end{rem}

\begin{proof}
With the notation of Lemma~\ref{res-commutators},
any effectively resonant monomial (defined in Definition~\ref{Def:effres})
gives rise to a term of the form
$x^{-1}Q_\beta B^\beta$ with $\sum_k s_k\beta_k=1$, while the
effectively non-resonant terms (defined in Definition~\ref{Def:effnonres})
are of the form
$x^{-1}Q_\beta B^\beta$ with $\sum_k s_k\beta_k>1$.
This is indeed the key point in categorizing resonant terms as effectively
resonant or nonresonant; see
the proof of Theorem~\ref{thm:enr-structure}. But if $\ep=\sum s_k\beta_k-1>0$,
we can rewrite $x^{-1}Q_\beta B^\beta\sim x^{\ep}Q_\beta A^\beta$ (\ie the
difference of the two sides is in $\Psisc^{-\infty,0}(X)$), and
$Q_\beta A^\beta\in\cM^{|\beta|}$. Since there are only finitely
many effectively non-resonant terms in \eqref{effnonres-span}, we deduce that
any $\tilde P$ with
$\sigma_1(\tilde P)=p_{\norm}$ may be written 
\begin{equation}\label{eq:Pt-form-0}
\lambda^{-1}\tilde P = xD_x +\sum_{j=1}^{m-1} r_j y_j D_{y_j} +\sum_{j=m}^{n-1}Q_j(x^{-1/2} y_j,x^{1/2}D_{y_j}) + R_{\effr}+ B + \tilde R,
\end{equation}
where $R_{\effr}$ is as in \eqref{eq:Pt-form},
$\tilde R \in x^{\epsilon} \cM_J$ for some $\epsilon > 0$, and $B \in \Psisc^{*,0}(X)$. Note that $\cP_j$ and $\cP_0$ are polynomials, and the homogeneity
claim is the meaning of the resonance condition
Proposition~\ref{prop:linear}.
For $s\leq j\leq m-1$, $\cP_j$ is
independent of $y_j,\ldots,y_{m-1}$ since
$0<r_s\leq r_{s+1}\leq\ldots\leq r_{m-1}$; $y_j$ itself cannot appear
in $\cP_j$ due to the restriction $2a+|\beta|+|\gamma|\geq 3$
in Proposition~\ref{prop:linear}.
Similarly, for $j\leq s-1$,
$\cP_j$ is independent of $y_1,\ldots,y_j$ as $r_1\leq
r_2\leq\ldots\leq r_{s-1}<0$.
This also shows that the polynomials $\cP_j$, $j\neq 0$,
have no constant or linear terms.

Let $B$ have symbol $b(\nu, y, \mu)$. Modulo terms in $x^{\epsilon} \cM^j$,
this can be reduced to the symbol $b'(0, (y', 0, 0), 0)$. Finally, by conjugating $P_{\norm}$ by a function $e^{if(y')}$, we can remove the $y'$-dependence of $b'$. Indeed, the Taylor series of $f$ can be constructed iteratively. Let 
$\cI'$ denote the ideal of functions of $y'$ that vanish at $0$. Conjugating
$\tilde P$ by $e^{if}$ produces the terms $\sum_{j=1}^{s-1} r'_j y'_j D_{y'_j} f$, as well as terms from $R_{\effr}$, which map $(\cI')^k\to(\cI')^{k+1}$.
For $k\geq 1$, $f\mapsto \sum_{j=1}^{s-1} r'_j y'_j \pa_{y'_j} f$
defines a linear map on $(\cI')^k$, $k\geq 1$,
with all eigenvalues negative since
$r'_j<0$ for $j=1,\dots,s-1$. Thus, this map is invertible, and this shows
that $b'-b'(0)$ can be conjugated away in Taylor series. Then it is
straightforward to check that the infinite order vanishing error can also be
removed.
\end{proof}

Later in this section we show that if $p_{\norm}$ is effectively non-resonant,
the leading asymptotics of microlocally outgoing solutions
for \eqref{eq:Pt-form} and for the completely explicit operator
\begin{equation}
P_0 = xD_x+\sum_{j=1}^{m-1} r_j y_jD_{y_j}+\sum_{j=m}^{n-1}Q_j(x^{-\frac{1}{2}}y_j,x^{\frac1{2}}D_{y_j}) + b , \quad b \text{ constant}
\label{effnonres-modelop}\end{equation}
are the same, if $R \in x^{1 + \epsilon} \cM^j$ for some $\epsilon > 0$,
\ie $R$ is indeed an `error term'. An analogous conclusion holds
in the effectively resonant case, with $R_{\effr}$ included
in the right hand side of \eqref{effnonres-modelop}.

First, however, we study the asymptotics
of approximate solutions of $P_0 u=0$. The constant $b$ simply introduces
a power $x^{-ib}$ into the asymptotics, as can be seen by conjugation of
$P_0$ by $x^{-ib}$. Here it is convenient to have the asymptotics for the
ultimately relevant case, where the operator $xP$ is self-adjoint, stated
explicitly, so we assume that $xP_0$ is formally self-adjoint on
$L^2_{\scl}(X)$, which
amounts to
\begin{equation}\label{eq:im-b-form-sa}
\im b = \frac{n-1}{2} - \frac1{2} (\sum_{j=1}^{s-1}r'_j+\sum_{j=s}^{m-1} r''_j)
- \frac{n-m}{2},
\end{equation}
provided that we have already made $Q_j$ self-adjoint as stated in
Remark~\ref{rem:Q_j-quantization}. Note that $\frac{n-m}{2}=\sum_{j=m}^{n-1}
\re r'''_j$.

For convenience, we separate the case where $\Rp$ is a source/sink
of $W$, hence of the contact vector field of $P_0$.
Recall from the previous section that
\begin{equation}\label{eq:def-Y''}
Y''_j = x^{-r''_j} y''_j, \quad Y''' = x^{-1/2} y''', 
\end{equation}
and define the exponents
\begin{equation}\label{eq:def-exps}
\tilde b=b-i\frac{n-m}{4},
\ a_{\beta'} = -\sum_{j=1}^{s-1}r_j\beta'_j-i\tilde b.
\end{equation}
Notice that $\re a_{\beta'} \to \infty$ as $|\beta'| \to \infty$. 

\begin{prop}\label{prop:enr-min} Suppose that the radial point $\Rp$ is a source/sink of $W$, and \eqref{eq:im-b-form-sa} holds.
Suppose that $u \in I^{(s)}(O,\cM)$, and $P_0 u \in I^{(s')}(O,\cM)$ where  $s < -1/2 < s'$. Then $u$ takes the form
\begin{equation}\label{eq:homog-summand-1}
u = \sum_{k} x^{-i\tilde b - i\kappa_{k}} w_{k}(Y'')v_{k}(Y''')
+u'
\end{equation}
where the sum is over $k \in \NN$, $v_{k}(Y)$ is an $L^2$-normalized eigenfunction of the harmonic oscillator
\begin{equation}
\sum_{j=m}^{n-1}\tilde Q_j(Y_j,D_{Y_j}),\ \tilde Q_j(Y_j,
D_{Y_j})=Q_j(Y_j,D_{Y_j})-\frac{1}{4}(Y_j D_{Y_j}+D_{Y_j}Y_j),
\ Y_j=\frac{y_j}{x^{1/2}},
\label{model-efn}\end{equation}
with eigenvalue $\kappa_{k},$ $w_{k}$ are Schwartz functions with each
seminorm rapidly decreasing in $k$, and $u' \in I^{(s' - \epsilon)}(O,\cM)$
for every $\epsilon > 0.$ 

Conversely, given any rapidly decreasing Schwartz sequence, $w_{k},$ in
$Y'',$ meaning one for which all seminorm rapidly decreasing in $k,$
and given any $f\in\Isc^{(s')}(O,\cM)$,
there exists $u\in\cap_{s<-1/2}\Isc^{(s)}(O,\cM)$ of the form
\eqref{eq:homog-summand-1} with $\WFsc(P_0 u-f)\cap O=\emptyset$.
\end{prop}

\begin{rem}
The result is true if we only assume $s<s'$. However, if $s\geq -1/2$,
we can replace $s$ by $\tilde s>-1/2$, apply the proposition with $\tilde s$
in place of $s$, and then use $u\in\Isc^{(s)}(O,\cM)$ to show that each
$w_k$ vanishes. On the other hand, if $s'\geq -1/2$, the proof of the
proposition shows that $u\in\Isc^{(s)}(O,\cM)$ implies
$u\in\Isc^{(s'-\ep)}(O,\cM)$ for every $\ep>0$.
\end{rem}

\begin{prop}\label{prop:enr-saddle} Suppose that $\Rp$ is a saddle point of
$W$, and \eqref{eq:im-b-form-sa} holds.
Suppose $u \in I^{(s)}(O,\cM)$, and $P_0 u \in I^{(s')}(O,\cM)$ for some $s < s' < \infty$. Then $u$ takes the form 
\begin{equation}\label{eq:homog-summand-2}
u = \sum_{\beta', k} x^{a_{\beta'} - i\kappa_{k}}
(y')^{\beta'}w_{\beta', k}(Y'')v_{k}(Y''')
+u'
\end{equation}
where the sum is over $k \in \NN$ and a finite set of multiindices $\beta'$, $v_{k}(Y)$ and $\kappa_k$ are as above,
$w_{\beta', k}$ is a rapidly decreasing  Schwartz sequence and $u' \in
I^{(s'-\ep)}(O,\cM)$ for every $\ep>0.$

Conversely, given any rapidly decreasing sequence of Schwartz functions
$w_{\beta', k}$, finite in $\beta'$ and any
$f\in\Isc^{(s')}(O,\cM)$ there exists
$u\in\cap_{s<-1/2}\Isc^{(s)}(O,\cM)$ of the form
\eqref{eq:homog-summand-2} with $\WFsc(P_0 u-f)\cap O=\emptyset$.
\end{prop}

\begin{rem}\label{eq:rem-orders-at-saddle}
As shown later, $x^2D_x$ gives rise to the terms in $\tilde Q-Q$ after
the change of variables $(x,y_j)\mapsto (x,\frac{y_j}{x^{1/2}})$.
If $Q_j$ is self-adjoint on $L^2(\Real, dY_j)$ then
$\tilde Q_j$ has the same property.

Also, with
\begin{equation*}
B=\frac{n-1}{2} - \frac1{2} \sum_j r''_j - \frac{n-m}{4},
\end{equation*}
the $(\beta',k)$ summand in \eqref{eq:homog-summand-2}
is in $\Isc^{(\re a_{\beta'}-B-\frac{1}{2}-\ep)}(O,\cM)$ for
every $\ep>0$. We show below that $\im\tilde b=B+d$,
$d=-\frac{1}{2}\sum r'_j>0$, so the $(\beta',k)$ summand is in
$\Isc^{(d-\sum r_j\beta'_j-\frac{1}{2}-\ep)}(O,\cM)$ for
every $\ep>0$, and in view of
the rapid decay in $k$, the same is true after the $k$ summation.
Thus, for $u$ as in \eqref{eq:homog-summand-2}, $u\in
\Isc^{(d-\frac{1}{2}-\ep)}(O,\cM)$ provided $s'>d-\frac{1}{2}$, \ie decays
by a factor $x^d$ faster than the microlocal solutions at sources/sinks of $W.$
\end{rem}

\begin{proof}[Proof of Proposition~\ref{prop:enr-min}]
Suppose that $P_0u = f \in I^{(s')}(O,\cM)$ for some $s' > -1/2$.
Let $O'$ be a $W$-balanced neighbourhood of $\Rp$ with
$\overline{O'}\subset O$, and let $Q\in\Psisc^{-\infty,0}(X)$ satisfy
$\WFsc'(Q)\subset O$ (\ie $Q\in\Psisc^{-\infty,0}(O)$) and
$\WFsc'(\Id-Q)\cap \overline{O'}=\emptyset$, with Schwartz kernel
supported in $U\times U$,
\begin{equation*}
U=\{0\leq x<\delta,\ |y_j|<\delta\ \text{for all}\ j\}.
\end{equation*}
(See \eqref{eq:def-Phi} for the definition of the diffeomorphism $\Phi$,
the coordinates $Y_j$, etc.)
Then, as noted in \eqref{eq:conormal-reg}, by the definition of
$\Isc^{(s)}(O,\cM)$, $\tilde u=Qu$ satisfies
\begin{equation*}
(Y'')^{\gamma''}(Y''')^{\gamma'''}(xD_x)^a D_{Y''}^{\beta''}
D_{Y'''}^{\beta'''} \tilde u\in x^s L^2_{\scl}(X)
\end{equation*}
for all $a$, $\beta''$, $\beta'''$, $\gamma''$ and $\gamma'''.$
Here $\tilde u$ is a microlocalization of $u$ since
$\WFsc(u-Qu)\subset\WFsc'(\Id-Q)$, so $\WFsc(u-Qu)\cap O'=\emptyset.$
Moreover,
\begin{equation*}
P_0 (Qu)=Q P_0 u+[P_0,Q]u=Qf+f',\ f'\in\dCinf(X),
\end{equation*}
since $\WFsc(u)\cap O\subset\{\Rp\}$, while $\WFsc'([P_0,Q])\subset
\WFsc'(Q)\cap\WFsc'(\Id-Q)\subset O\setminus\overline{O'}$,
so $\WFsc(u)\cap\WFsc'([P_0,Q])=\emptyset$. Thus, with $\tilde f=Qf+f'$,
\begin{equation}\begin{split}\label{eq:def-tilde-f}
&P_0 \tilde u=\tilde f,\\
&(Y'')^{\gamma''}(Y''')^{\gamma'''}
(xD_x)^a D_{Y''}^{\beta''} D_{Y'''}^{\beta'''} \tilde f
\in x^{s'} L^2_{\scl}(X),
\end{split}\end{equation}
for all $a$, $\beta''$, $\beta'''$, $\gamma''$ and $\gamma'''$.

To prove first part of the proposition, it thus suffices to show that,
with the notation of \eqref{eq:homog-summand-1},
\begin{equation}\label{eq:homog-summand-1p}
\tilde u = \sum_{k} x^{-i\tilde b - i\kappa_{k}} w_{k}(Y'')v_{k}(Y''')
+u'.
\end{equation}

Writing the operator $P_0$ in the coordinates $x, Y'', Y'''$ we have
\begin{equation}\label{eq:P_0-enr-min}
P_0 = x D_x |_{Y} + \sum_j \tilde Q_j(Y'''_j, D_{Y'''_j}) + \tilde b
\end{equation}
with $\tilde b=b-i\frac{n-m}{4}$ as in \eqref{eq:def-exps}.
Formal self-adjointness of $x P_0$, \ie \eqref{eq:im-b-form-sa},
requires that 
\begin{equation}
\im \tilde b = \frac{n-1}{2} - \frac1{2} \sum_j r''_j - \frac{n-m}{4} \equiv B.
\label{Beqn}\end{equation}

As already remarked, \eqref{eq:def-tilde-f}, which states that $\tilde f$
is conormal in $x$, and Schwartz in $Y'', Y'''$, and belongs to $x^{s'}
L^2(dx dy/x^{n+1})$, or in terms of the $Y$ coordinates, to $x^{s' + n/2 -
\sum r''_j/2 - (n-m)/4}~L^2(dx dY/x)$, implies (by conormality) that
$$
\tilde f\in x^{s' + 1/2 + B - \epsilon} L^\infty
$$
for every $\epsilon> 0$, where $B$ is defined by \eqref{Beqn}. More
precisely, for all $a,$ $\beta$, $\gamma''$ and $\gamma''',$
\begin{equation*}
(Y'')^{\gamma''}(Y''')^{\gamma'''}
(xD_x)^a D_{Y''}^{\beta''} D_{Y'''}^{\beta'''} \tilde f
\in x^{s' + 1/2 + B - \epsilon} L^\infty
\end{equation*}
for every $\epsilon > 0$.
Conversely these conditions imply that
$\tilde f $ satisfies \eqref{eq:def-tilde-f} with $s'$ replaced
by $s'-\ep$ for every $\epsilon >0$.

Writing
$\tilde f$ in the form 
\begin{equation*}
\tilde f(x, Y'', Y''') = \sum_k f_k(x, Y'') v_k(Y'''),
\end{equation*}
where $f_k$ is conormal in $x$, rapidly decreasing as a Schwartz sequence
in $Y'',$ a particular solution to $P_0 \tilde u = \tilde f$, is given by
\begin{equation}\begin{gathered}
\tilde u = \sum_k u_k(x, Y'') v_k(Y'''), \\
u_k = - i x^{-i\tilde b - i \kappa_k} \int_0^x f_k(t, Y'') t^{i\tilde b + i \kappa_k} \, \frac{dt}{t}.
\end{gathered}\label{particular}\end{equation}
Since $s' + 1/2 > 0$, this integral is convergent and $\tilde u \in I^{(s'
  - \epsilon)}(O,\cM)$ for every $\epsilon > 0.$ 

On the other hand, the general solution to $P_0\tilde u = 0$ with
$\tilde u$ Schwartz in $Y''$ and $Y'''$ is given by 
\begin{equation}\label{eq:homog-summand-8}
\tilde u = \sum_{k} x^{-i\tilde b - i\kappa_{k}} w_{k}(Y'')v_{k}(Y'''),
\end{equation}
where $w_k$ is rapidly decreasing in $k$. Since any solution is the sum of
the particular solution \eqref{particular} and some homogeneous solution,
the first half of the proposition follows.

In fact, the second half also follows by defining
\begin{equation*}
\tilde u=\sum_k u_k(x, Y'') v_k(Y''')+
\sum_{k} x^{-i\tilde b- i\kappa_{k}} w_{k}(Y'')v_{k}(Y'''),
\end{equation*}
with $u_k$ as in \eqref{particular}. Multiplying by a cutoff function
$\phi\in\Cinf(X)$ which is identically $1$ near $(0,0,\ldots,0)$, it
follows that $u=\phi\tilde u$ satisfies all requirements.
\end{proof}

\begin{proof}[Proof of Proposition~\ref{prop:enr-saddle}]
We use a similar argument to prove this result. Let $O'$, $Q$, etc.,
be as in the previous proof.
With $\tilde u=Qu$, as noted in \eqref{eq:conormal-reg},
\begin{equation}\label{eq:tilde-u-conormal}
(Y'')^{\gamma''}(Y''')^{\gamma'''}(xD_x)^a D_{y'}^{\beta'}
D_{Y''}^{\beta''} D_{Y'''}^{\beta'''} \tilde u\in x^{s} L^2_{\scl}(X),
\end{equation}
for all $a$, $\beta$, $\gamma''$ and $\gamma'''$.
One of the main differences with the proof of Proposition~\ref{prop:enr-min}
is that microlocalization introduces a non-trivial error, i.e.\ $P_0\tilde u$
is not globally well-behaved (not as good as $f$ was microlocally).
However, the error is supported away from $y'=0$. Indeed,
now $\WFsc(u)\cap O\subset S$, and
\begin{equation*}
\tilde f=P_0\tilde u=Qf+f',\ f'=[P_0,Q]u.
\end{equation*}
Here $\WFsc'([P_0,Q])\cap S\subset\{|y'|> \delta_0\}$ for some
$\delta_0>0$, so $f'\in\Isc^{(s)}(O,\cM)$ in fact satisfies
\begin{equation*}
(Y'')^{\gamma''}(Y''')^{\gamma'''}(xD_x)^a D_{y'}^{\beta'}
D_{Y''}^{\beta''} D_{Y'''}^{\beta'''} f'\in x^s L^2_{\scl}(X)
\end{equation*}
for all $a$, $\beta'$, $\beta''$ and $\beta'''$, $\gamma''$ and $\gamma'''$,
with the improved conclusion
\begin{equation*}
\phi(y')(Y'')^{\gamma''}(Y''')^{\gamma'''}(xD_x)^a D_{y'}^{\beta'}
D_{Y''}^{\beta''} D_{Y'''}^{\beta'''} f'\in \dCI(X)
\end{equation*}
if $\phi$ is supported in $|y'|<\delta_0$. Correspondingly,
\begin{equation}\label{eq:def-tilde-f-2}
\phi(y')(Y'')^{\gamma''}(Y''')^{\gamma'''}(xD_x)^a D_{y'}^{\beta'}
D_{Y''}^{\beta''} D_{Y'''}^{\beta'''} \tilde f\in x^{s'} L^2_{\scl}(X).
\end{equation}

The operator $P_0$ in the coordinates $x, y', Y'', Y'''$ now takes the form 
\begin{equation}\label{eq:P_0-enr-saddle}
P_0 = x D_x |_{y',Y'',Y'''} + \sum_j r'_j y'_j D_{y'_j} + \sum_j \tilde Q_j(Y'''_j, D_{Y'''_j}) + \tilde b,
\end{equation}
with $\tilde b=b-i\frac{n-m}{4}$ as in \eqref{eq:def-exps}.
Again, \eqref{eq:def-tilde-f-2} implies that $\tilde f$
is conormal in $x$, smooth in $y'$, and Schwartz in $Y'', Y'''$,
and belongs to $x^{s' + 1/2 + B - \epsilon} L^\infty$ for every $\epsilon > 0$, where $B$ is defined by \eqref{Beqn}, in the precise sense that
for all
$a$, $\beta$, $\gamma''$ and $\gamma'''$,
\begin{equation*}
\phi(y')(Y'')^{\gamma''}(Y''')^{\gamma'''}
(xD_x)^a D_{Y''}^{\beta''} D_{Y'''}^{\beta'''} \tilde f
\in x^{s' + 1/2 + B - \epsilon} L^\infty
\end{equation*}
for every $\epsilon > 0$.
However, now formal self-adjointness of $xP_0$ requires that
\begin{equation}
\im \tilde b = B + d, \quad d = -\frac1{2} \sum_j r'_j > 0,
\label{Beqn2}\end{equation}
so there is a discrepancy of $d$ compared with the previous proposition.  Write $\tilde f$ in the form 
\begin{equation*}
\tilde f(x, Y'', Y''') = \sum_k f_k(x, y', Y'') v_k(Y'''),
\end{equation*}
where $f_k$ is rapidly decreasing sequence which is conormal in $x,$ smooth
in $y'$ and Schwartz in $Y''.$

We start by describing solutions of the homogeneous equation $P_0\tilde u=0$
in $U$ which in addition satisfy \eqref{eq:tilde-u-conormal}.
Decomposing $\tilde u$ in terms of the $v_k$, and factoring
out a power of $x$ for convenience, i.e.\ writing
$\tilde u=\sum_k x^{-i\tilde b-i\kappa_k}u_k(x,y',Y'')v_k(Y''')$,
we see that the coefficients
$u_k$ satisfy
\begin{equation*}
(x\pa_x|_{y',Y'',Y'''}+\sum_j r'_j y'_j\pa_{y'_j})u_k
=0.
\end{equation*}
Since $\tilde u$ is smooth in the interior of $U$,
$P_0\tilde u=0$ amounts to
demanding that $u_k$ be constant
along each integral curve segment of the vector field
$x\pa_x+\sum_j r'_j y'_j\pa_{y'_j}$, with the value of $\tilde u$ depending
smoothly on the choice of the integral curve.
(We remark that $U$ is convex for this vector field; $|y'|$ is increasing
as $x\to 0$.) Thus, $u_k(x,y',Y'')=\hat u_k(Y',Y'')$
with $\hat u_k$ smooth in $Y'$ and Schwartz in $Y''$.
Here $Y'_j=y'_j/x^{r'_j}$; note that $r'_j<0$. Expanding $\hat u_k$ in
Taylor series around $Y'=0$ to order $N$, we see that
\begin{equation*}
u_k(x,y',Y'')-\sum_{|\beta'|\leq N-1} x^{-\sum_j r'_j\beta'_j} (y')^{\beta'}
w_{\beta',k}(Y'')
\end{equation*}
is a finite sum of terms of the form $x^{-\sum_j r'_j\beta'_j}(y')^{\beta'}
\hat u_{k,\beta'}(Y',Y'')$ with $\hat u_{k,\beta'}$ smooth (Schwartz in
$Y'''$), where the sum runs over $\beta'$ with $|\beta'|=N$. Thus,
given any $s''$ (e.g.\ $s''=s'$), we can choose $N$ sufficiently large so
this difference lies in $\Isc^{(s'')}(O,\cM)$, which means it is
ignorable for our purposes. Thus, the general solution to $P_0\tilde u = 0$
in $U$ which satisfies \eqref{eq:tilde-u-conormal} is given by
\begin{equation}\label{eq:homog-summand-16}
\tilde u = \sum_{\beta', k} x^{a_{\beta'} - i\kappa_{k}}
(y')^{\beta'}w_{\beta', k}(Y'')v_{k}(Y'''),
\end{equation}
modulo any $\Isc^{(s'')}(O,\cM)$ (where the sum is understood as a finite
one, due to the remark above),
where the seminorms of $w_{\beta', k}$ are rapidly decreasing in $k$ for each $\beta'$.

In expressing a particular solution $\tilde u$ of $P_0 \tilde u = f$
in terms of $f$, we
need to integrate along integral curves of the vector field
$x \partial_x + \sum_j r'_j y'_j \partial_{y'_j}$, and
since $r'_j <0$, $|y'| \to \infty$ as $x \to 0$ along such curves
(unless $y' = 0$); in fact $|y'|$ is increasing as $x\to 0$ as mentioned above.
So we cannot integrate down to $x=0$. Instead we fix an $x_0 > 0$ and use
the formula
\begin{equation}\begin{gathered}\label{eq:particular-2}
u_k(x, y', Y'') = (\frac{x}{x_0})^{-i\tilde b - i \kappa_k} u_k(x_0, (\frac{x}{x_0})^{-r'_j} y'_j, Y'') \\ + i x^{-i\tilde b - i \kappa_k} \int_{x_0}^x  f_k(t, (\frac{x}{t})^{-r'_j} y'_j, Y'') t^{i\tilde b + i \kappa_k} \, \frac{dt}{t}.
\end{gathered}\end{equation}
Notice that $u_k(x_\sharp,y'_\sharp,Y''_\sharp)$
depends only on $f_k$ evaluated at points $(x,y',Y'')$ with
$|y'|\leq |y'_\sharp|$. Thus, \eqref{eq:def-tilde-f-2} can be used to
deduce properties of $u_k$, hence of $\tilde u$, in $|y'|<\delta_0$.

If $s' < -1/2 + d$, then \eqref{eq:particular-2} gives
$\phi(y')\tilde u \in I^{(s' - \epsilon)}(O,\cM)$
for every $\epsilon > 0$, with $\phi$ as in
\eqref{eq:def-tilde-f-2}. If $s' \geq -1/2 + d$, then 
$\phi(y')\tilde u \in I^{(-1/2 + d - \epsilon)}(O,\cM)$
for every $\epsilon > 0$.
However, this is actually a sum of terms solving the homogeneous
equation, plus a function in $I^{(s' - \epsilon)}(O,\cM)$ for
every $\epsilon > 0$. For simplicity we show this only in the case that
$-1/2 + d < s' < -1/2 + d + |r'_{s-1}|$. Then we observe that
$(x/x_0)^{-i\tilde b - i \kappa_k}\tilde u(x_0, 0, Y'')$ is a solution of the
homogeneous equation, while the difference
$$\begin{gathered}
(\frac{x}{x_0})^{-i\tilde b - i \kappa_k} \tilde u(x_0, (\frac{x}{x_0})^{-r'_j} y'_j, Y'')
- (\frac{x}{x_0})^{-i\tilde b - i \kappa_k} \tilde u(x_0, 0, Y'')  \\ = 
 \sum_j (\frac{x}{x_0})^{-r'_j} 
\int_0^1 y'_j \partial_{y'_j} \big( \tilde u(x_0, \tau(\frac{x}{x_0})^{-r'_j} y'_j, Y'') \big) \, d\tau
\end{gathered}$$
has decay at least $x^{-r'_{s-1}}$ better, hence yields a term in
$I^{(s' - \epsilon)}(O,\cM)$ for every $\epsilon > 0$. Similarly, if we
replace $f_k(t, (\frac{x}{t})^{-r'_j} y'_j, Y'')$ in the integral by
$f_k(t, 0, Y'')$ then we get a homogeneous term, while the difference gives
a term in $I^{(s' - \epsilon)}(O,\cM)$ for every $\epsilon > 0$. The
argument can be repeated, removing more and more terms in the Taylor series
for $\tilde u$ and $\tilde f$, for larger values of $s'$.
Since any solution is the sum of
the particular solution above and the general solution, the first half of the
proposition follows with $O$ replaced by a smaller neighbourhood $O''$
of $\Rp$. However, we recover the original statement by using the
real principal type parametrix construction of Duistermaat and H\"ormander
\cite{FIOII}.

The second half can be proved as in the previous proposition.
Fix some $x_0>0$, and let $u_k$ be given by the second term on the right hand side
of \eqref{eq:particular-2}, and let
$\hat u=\sum_k u_k(x, Y'') v_k(Y''')$. Then $P_0 \hat u=f$, and
as shown above,
$\hat u$ has the form
\begin{equation}\label{eq:homog-summand-2b}
\hat u = \sum_{\beta', k} x^{a_{\beta'} - i\kappa_{k}}
(y')^{\beta'}\hat w_{\beta', k}(Y'')v_{k}(Y''')
+\hat u',
\end{equation}
with $\hat u'\in\Isc^{(s'-\ep)}(O,\cM)$ for all $\ep>0$.
Then with
\begin{equation*}
\tilde u=\sum_k u_k(x, Y'') v_k(Y''')+
\sum_{\beta',k} x^{a_{\beta'}- i\kappa_{k}}
(w_{\beta',k}(Y'')-\hat w_{\beta', k}(Y''))v_{k}(Y'''),
\end{equation*}
$u=\phi\tilde u$, $\phi\in\Cinf(X)$ identically $1$ near $(0,\ldots,0)$,
$u$ satisfies all requirements.
\end{proof}

These results on the explicit normal form $P_0$ then allow us to
parameterize microlocally outgoing solutions for every effectively
nonresonant critical point.

\begin{thm}\label{thm:enr-structure} Suppose that $P(\ev)$ is effectively
  nonresonant at $\Rp$, with normal form $P_{\norm}$ near $\Rp$ as in
  Lemma~\ref{effnonres-form}, and \eqref{eq:im-b-form-sa} holds.

\begin{enumerate}
\item If in addition $\Rp$ is a source/sink of $W$, then
any microlocally outgoing solution $u$ of $P_{\norm}$ has the form
\eqref{eq:homog-summand-1}, and conversely given any Schwartz sequence of
Schwartz functions $w_k$ there is a microlocally outgoing solution $u$ of
$P_{\norm}$ which has the form \eqref{eq:homog-summand-1}. Thus, microlocal
solutions at a source/sink of $W$ are parameterized by Schwartz functions
of the variables $(Y'', Y''')$.

\item
If $\Rp$ is a saddle point of $W$, then all microlocally outgoing
solutions are in $x^{-1/2 + \epsilon} L^2$ for some $\epsilon > 0$.
For each monomial $(y')^\beta$ in the variables $y'$, each $k \in \NN$ and each
Schwartz function $w(Y'')$ there is a microlocally outgoing solution of the form 
\begin{equation}
u = \sum_k x^{a_{\beta'} - i\kappa_{k}}
(y')^{\beta'}w(Y'')v_{k}(Y''')
+u', 
\label{eq:saddle-exp-2}
\end{equation}
where $u'$ is in a strictly smaller weighted $L^2$ space than $u,$
and every microlocally outgoing solution is a sum of such solutions, with
the $w=w_{k,\beta'}$ rapidly decreasing as $k\to\infty$ in every seminorm.
\end{enumerate}
\end{thm}

\begin{proof}
First, $P_{\norm}=\lambda (P_0+R)$, $R\in x^\ep\cM^j$, $\ep>0$. Thus, if $O$ is
a neighbourhood of $\Rp$ as above, $\WFsc(P_{\norm} u)\cap O=\emptyset$,
then $u\in\Isc^{(s)}(O,\cM)$ for all $s<-1/2$, so
$Ru\in \Isc^{(s')}(O,\cM)$ for some $s'>1/2$. Hence
$P_0u=\lambda^{-1}P_{\norm}u-Ru\in \Isc^{(s')}(O,\cM)$.

If $\Rp$ is a source/sink of $W$, then
Proposition~\ref{prop:enr-min} is applicable, and we deduce
that $u$ is microlocally of the form \eqref{eq:homog-summand-1}.
Moreover, if $\Rp$ is a source/sink of $W$, then given any
Schwartz sequence of Schwartz functions $w_k$, let
$u_0\in\cap_{s<-1/2}\Isc^{(s)}(O,\cM)$ be of the form \eqref{eq:homog-summand-1}
with $P_0 u_0\in\dCinf(X)$. We construct $u_k\in
\cap_{r<-1/2-k\ep}\Isc^{(r)}(O,\cM)$, $k\geq 1$, inductively so that
$P_0 u_k+Ru_{k-1}\in\dCinf(X)$ for $k\geq 1$; this can be done
by the second half of Proposition~\ref{prop:enr-min}. Asymptotically
summing $\sum_k u_k$ to some $u\in\cap_{s<-1/2}\Isc^{(s)}(O,\cM)$ gives
a microlocally outgoing solution with the prescribed asymptotics,
completing the proof of the theorem in this case.

If $\Rp$ is a saddle point of $W$, we apply
Proposition~\ref{prop:enr-saddle} with $s'>-1/2$ as in the first
paragraph of the proof. If $\ep'>0$ is sufficiently small,
all of the terms in \eqref{eq:homog-summand-2}
are in $\Isc^{(-1/2+\ep')}(O,\cM)$ proving the first claim. To show
the next, let
$u_0=x^{a_{\beta'} - i\kappa_{k}}
(y')^{\beta'}w(Y'')v_{k}(Y''')$, so $P_0u_0=0$ and $u_0\in\Isc^{(s)}(O,\cM)$
for any $s<-1/2+d$.
We construct $u_k$ inductively as above, using
Proposition~\ref{prop:enr-saddle}, to obtain $u$.
\end{proof}

\begin{rem}\label{original-exp} Fr{o}m \eqref{eq:homog-summand-1} or
\eqref{eq:saddle-exp-2} it is not hard to derive the asymptotic expansion of
eigenfunctions of the original operator $\Delta + V - \ev$; we need only
apply the Fourier integral operator $F^{-1}$ arising by composing
any Fourier integral operators with canonical relation given by
the contact maps in
Lemma~\ref{HMV2r.191} and Theorem~\ref{thm:model} to these expansions.
In fact, as mentioned in Remark~\ref{rem:special-contact}, this Fourier
integral operator can be taken to be a composition of a change of
coordinates with multiplication by an oscillatory function if $\Rp$
is either a source/sink (so $q \in \Min_+(\ev)$) or the linearization
of $W$ has no non-real eigenvalues (so there are no $y'''$ variables).

In the
case of a radial point $q \in \Min_+(\ev)$, {\em in appropriate coordinates
$y$ on $\pa X$}, the expansion takes the form
\begin{equation}\label{eq:homog-summand-111}
u = e^{i\Phi(y)/x}\sum_{k} x^{-i\tilde b - i\kappa_{k}} w_{k}(Y'')v_{k}(Y''')
+u', \ u' \in I^{-\frac1{2} + \epsilon}(O, \cM) \, \text{ for some } \epsilon > 0
\end{equation}
where $\Phi$ is a smooth function (it parameterizes the Legendrian
submanifold which is the image of the zero section under the canonical
relation of $F^{-1}$). For a given $\ev$, only a finite number of terms in
the Taylor series for $\Phi$ are relevant. Similarly in the case of radial
points $q \in \RP_+(\ev) \setminus \Min_+(\ev)$, the expansion
\eqref{eq:saddle-exp} takes the form
\begin{equation}
u = e^{i\Phi(y)/x} \sum_k x^{a_{\beta'} - i\kappa_{k}}
(y')^{\beta'}w(Y'')v_{k}(Y''')
+u', 
\label{eq:saddle-exp}\end{equation}
with $\Phi$ smooth. Again it parameterizes the image of the zero section
under the canonical relation of $F^{-1}$.  In this case, the value of
$\Phi$ on the unstable manifold $\{ y'' = y''' = 0 \}$ is essential, but
only a finite number of terms in the Taylor series for $\Phi$ about this
unstable manifold are relevant.

These expansions were obtained directly in Part I (\ie without going via a
normal form) in the two dimensional case.
\end{rem}


\section{Effectively resonant operators}\label{sec:er}

If $P$ is effectively resonant, the simple expressions
\eqref{eq:homog-summand-1} and \eqref{eq:homog-summand-2} need to be
replaced by slightly more complicated ones in which positive integral
powers of $\log x$ also appear. Essentially, instead of powers, or Schwartz
functions, of $\frac{y_j}{x^{r_j}}$, factors of $\log x$ also arise in the
expressions for the $Y_l.$

First define a change of coordinates inductively that simplifies the
vector field
\begin{equation}\label{eq:model-Vf}
V=(xD_x)+\sum_{j=1}^{m-1} (r_j y_j+\cP_j(y_s,\ldots,y_{j-1}))D_{y_j}
\end{equation}
that appears in \eqref{eq:Pt-form} as the combinations of the linear
terms $\sum r_j y_j D_{y_j}$ and the effectively resonant vector fields
in $R_{\effr}$.
(Note that $r_j y_j$ and $\cP_j(y_s,\ldots,y_{j-1})$ are both homogeneous of
degree $r_j$.) We do this in two steps to clarify the argument, first
only dealing with the $y''$ terms, \ie $j=s,\ldots,m-1$.

The coordinates $Y_j$, $j=s,\ldots,m-1$, are
a modification of the coordinates
$\frac{y_j}{x^{r_j}}$
that appear in \eqref{eq:def-Y''}, so that $Y_j-\frac{y_j}{x^{r_j}}$
are polynomials
$\cP^\sharp_j$ in $Y_s,\ldots,Y_{j-1},t=\log x$.
Thus, we let
\begin{equation*}
Y_s=\frac{y_s}{x^{r_s}},\ \cP^\sharp_s=0,
\ \bar Y_s(Y_s,\log x)=Y_s+\cP^\sharp_s(\log x)
\end{equation*}
and provided that $Y_s,\ldots,Y_{j-1}$, $\cP^\sharp_s,\ldots,\cP^\sharp_{j-1}$
have been defined, we let
\begin{equation*}\begin{split}
&\cP^\sharp_j(Y_s,\ldots,Y_{j-1},t)
=\int_0^t \cP_j(\bar Y_s(Y_s,t'),\ldots,\bar Y_{j-1}
(Y_s,\ldots,Y_{j-1},t'))\,dt',\\
&Y_j=\frac{y_j}{x^{r_j}}-\cP^\sharp_j(Y_s,\ldots,Y_{j-1},\log x),\\
&\bar Y_j=Y_j+\cP^\sharp_j(Y_s,\ldots,Y_{j-1},\log x),
\ j=s,\ldots,m-1.
\end{split}\end{equation*}
The point of the construction is that $V$ annihilates $Y_j$ for all $j$.
This can be seen iteratively: for $Y_s$ this is straightforward, and
if $VY_s=\ldots=VY_{j-1}=0$ then (with $\pa_t \cP^\sharp_j$ denoting the
derivative with respect to the last variable, $t=\log x$)
\begin{equation*}\begin{split}
&VY_j\\
&=-r_j\frac{y_j}{x^{r_j}}+(r_j y_j+\cP_j(y_s,\ldots,y_{j-1}))x^{-r_j}
-(\pa_t \cP^\sharp_j)(Y_s,\ldots,Y_{j-1},\log x)\\
&=\cP_j(y_s x^{-r_s},\ldots,y_{j-1}x^{-r_{j-1}})-
\cP_j(\bar Y_s(Y_s,\log x),\ldots,\bar Y_{j-1}
(Y_s,\ldots,Y_{j-1},\log x))\\
&=0
\end{split}\end{equation*}
in view of the definition of $Y_s,\ldots,Y_{j-1}$ and
$\bar Y_s,\ldots,\bar Y_{j-1}$.

One can deal with the $j=1,\ldots,s-1$ terms similarly.
We define $\cP^\sharp_j$, $Y_j$ and $\bar Y_j$ inductively
as above, starting with $Y_{s-1}$.
Thus, we let
\begin{equation*}
Y_{s-1}=\frac{y_{s-1}}{x^{r_{s-1}}},\ \cP^\sharp_{s-1}=0,
\ \bar Y_{s-1}(Y_{s-1},\log x)=Y_{s-1}+\cP^\sharp_{s-1}(\log x)
\end{equation*}
and provided that $Y_{j+1},\ldots,Y_{s-1}$,
$\cP^\sharp_{j+1},\ldots,\cP^\sharp_{s-1}$
have been defined, we let
\begin{equation}\begin{split}\label{eq:saddle-vars}
&\cP^\sharp_j(Y_{j+1},\ldots,Y_{s-1},t)
=\int_0^t \cP_j(\bar Y_{j+1}(Y_{j+1},\ldots,Y_{s-1},t'),\ldots,\bar Y_{s-1}
(Y_{s-1},t'))\,dt',\\
&Y_j=\frac{y_j}{x^{r_j}}-\cP^\sharp_j(Y_{j+1},\ldots,Y_{s-1},\log x),\\
&\bar Y_j=Y_j+\cP^\sharp_j(Y_{j+1},\ldots,Y_{s-1},\log x),
\ j=1,\ldots,s-1.
\end{split}\end{equation}
With these definitions,
in the coordinates $X=x,Y_1,\ldots,Y_{m-1},
y_m,\ldots,y_{n-1}$, \ie $(X,Y',Y'',y''')$,
which correspond to a blow-up of
$x=y_s=\ldots=y_{m-1}=0$, $V=X^2D_X$.

The zeroth order term is
a polynomial $\cP_0$ in $y_s,\ldots,y_{m-1}$ which is homogeneous of
degree $1$ (where $y_j$ has degree $r_j$). Thus,
\begin{equation*}
x^{-1}\cP_0(y_s,\ldots,y_{m-1})=\cP_0(\bar Y_s(Y_s,\log x),\ldots,
\bar Y_{m-1}(Y_s,\ldots,Y_{m-1},\log x)).
\end{equation*}
Let
\begin{equation*}
\cP^\sharp_0(Y_s,\ldots,Y_{j-1},t)
=\int_0^t \cP_0(\bar Y_s(Y_s,t'),\ldots,\bar Y_{j-1}
(Y_s,\ldots,Y_{j-1},t'))\,dt',
\end{equation*}
which is thus a polynomial in $Y_s,\ldots,Y_{j-1},t$. Then
$e^{i\cP^\sharp_0(Y_s,\ldots,Y_{j-1},\log x)}$ can be used as an integrating
factor, conjugating $\tilde P$, to remove the zeroth order term in $R_{\effr}$.

Finally, to put the quadratic terms in a convenient form, we let
\begin{equation*}
Y_j=\frac{y_j}{x^{1/2}},\ j=m,\ldots,n-1
\end{equation*}
as before.

Suppose first that $\cP_0=0.$
With our definition of the $Y_j$, \eqref{eq:P_0-enr-min},
resp.\ \eqref{eq:P_0-enr-saddle}, holds if $\Rp$ is a source/sink,
resp.\ saddle point, of $V_0$.
Thus, the statement and the proof of Proposition~\ref{prop:enr-min}
holds without any changes, while the statement
and the proof of Proposition~\ref{prop:enr-saddle} carry over provided
$x^{a_{\beta'}}(y')^{\beta'}$ is replaced by
$x^{-i\tilde b}(Y')^{\beta'}$. A minor difference is that slightly
more effort is required to show that $|y'|$ decreases on the integral
curves of the vector field \eqref{eq:model-Vf} inside $|y'|<\delta_1$
for $\delta_1>0$ small. Namely we need to use that, as $\cP_j$, $j=1,\ldots,
s-1$ have no linear or constant terms by Lemma~\ref{effnonres-form},
$V|y'|^2=\sum_{j=1}^{s-1}r_j y_j^2+\cO(|y'|^3)
\leq r_{s-1}|y'|^2+\cO(|y'|^3)$, $r_{s-1}<0$, to conclude that
$V|y'|^2\leq 0$ for $|y'|<\delta_1$, $\delta_1>0$ small.

In general, with
$\tilde b = b - i\frac{n-m}{4}$ as in \eqref{eq:def-exps},
\eqref{eq:P_0-enr-min},
resp.\ \eqref{eq:P_0-enr-saddle}, are replaced by
\begin{equation}\label{eq:P_0-er-min}
P_0 = x D_x |_{Y} + \sum_j \tilde Q_j(Y'''_j, D_{Y'''_j}) +\cP_0+ \tilde b,
\end{equation}
respectively
\begin{equation}\label{eq:P_0-er-saddle}
P_0 = x D_x |_{y',Y'',Y'''} + \sum_{j=1}^{s-1} (r'_j y'_j+\cP_j) D_{y'_j}
+ \sum_j \tilde Q_j(Y'''_j, D_{Y'''_j}) + \cP_0+\tilde b.
\end{equation}
Thus,
\begin{equation}\label{eq:P_0-c-er-min}
e^{i\cP^\sharp_0}P_0e^{-i\cP^\sharp_0} = x D_x |_{Y} + \sum_j \tilde Q_j(Y'''_j, D_{Y'''_j}) + \tilde b,
\end{equation}
respectively
\begin{equation}\label{eq:P_0-c-er-saddle}
e^{i\cP^\sharp_0}P_0e^{-i\cP^\sharp_0} = x D_x |_{y',Y'',Y'''}
+ \sum_{j=1}^{s-1} (r'_j y'_j+\cP_j) D_{y'_j} + \sum_j \tilde Q_j(Y'''_j, D_{Y'''_j}) +\tilde b.
\end{equation}
Since multiplication by $e^{\pm i\cP^\sharp_0}$ preserves $\Isc^{(s)}(O,\cM)$,
the rest of the proof of the propositions is applicable with $u$ replaced
by $e^{i\cP^\sharp_0}u$, $f=P_0u$ replaced by $e^{i\cP^\sharp_0}f$.
We thus deduce the following analogues of
Propositions~\ref{prop:enr-min} -- \ref{prop:enr-saddle} in the effectively
resonant case.

\begin{prop}\label{prop:er-min} Suppose that the radial point $\Rp$
is a source/sink of $W$, and \eqref{eq:im-b-form-sa} holds, that $u \in
I^{(s)}(O,\cM)$, and $P_0 u \in I^{(s')}(O,\cM)$ where  $s < -1/2 <
s'$. Then $u$ takes the form 
\begin{equation}\label{eq:homog-summand-1-c}
u = \sum_{k} x^{-i\tilde b - i\kappa_{k}} e^{-i\cP^\sharp_0}
w_{k}(Y'')v_{k}(Y''')
+u'
\end{equation}
where the sum is over $k \in \NN$, $v_{k}(Y)$ is an $L^2$-normalized eigenfunction of the harmonic oscillator
\begin{equation}
\sum_{j=m}^{n-1}\tilde Q_j(Y_j,D_{Y_j}),\ \tilde Q_j(Y_j,
D_{Y_j})=Q_j(Y_j,D_{Y_j})-\frac{1}{4}(Y_j D_{Y_j}+D_{Y_j}Y_j),
\ Y_j=\frac{y_j}{x^{1/2}},
\label{model-efn-c}\end{equation}
with eigenvalue $\kappa_{k}$,
$w_{k}$ are Schwartz functions with each seminorm rapidly decreasing in
$k$, and $u' \in I^{(s' - \epsilon)}(O,\cM)$ for every $\epsilon > 0$.

Conversely, given any sequence $w_{k}$ of Schwartz functions in $Y''$ with
each seminorm rapidly decreasing in $k$, and given any
$f\in\Isc^{(s')}(O,\cM)$, there exists $u\in\cap_{s<-1/2}\Isc^{(s)}(O,\cM)$
of the form \eqref{eq:homog-summand-1-c} with $\WFsc(P_0 u-f)\cap
O=\emptyset$.  \end{prop}

\begin{prop}\label{prop:er-saddle}
Suppose that $\Rp$ is a saddle point of $W$,
and \eqref{eq:im-b-form-sa} holds, that $u \in I^{(s)}(O,\cM)$, and
$P_0 u \in I^{(s')}(O,\cM)$ for some $s < s' < \infty$.
Then $u$ takes the form 
\begin{equation}\label{eq:homog-summand-2-c}
u = \sum_{\beta', k} x^{-i\tilde b - i\kappa_{k}}
(Y')^{\beta'} e^{-i\cP^\sharp_0}w_{\beta', k}(Y'')v_{k}(Y''')
+u'
\end{equation}
where the sum is over $k \in \NN$ and a finite set of multiindices $\beta'$, $v_{k}(Y)$ and $\kappa_k$ are as above,
$w_{\beta', k}$ are Schwartz functions with each seminorm rapidly decreasing in $k$, and $u' \in I^{(s'-\ep)}(O,\cM)$ for every $\ep>0$.   

Conversely, given any sequence of Schwartz functions
$w_{\beta', k}$, finite in $\beta'$
with each seminorm rapidly decreasing in $k$, and any
$f\in\Isc^{(s')}(O,\cM)$ there exists
$u\in\cap_{s<-1/2}\Isc^{(s)}(O,\cM)$ of the form
\eqref{eq:homog-summand-2-c} with $\WFsc(P_0 u-f)\cap O=\emptyset$.
\end{prop}

We thus deduce the following analogue of Theorem~\ref{thm:enr-structure},
with a similar proof.

\begin{thm}\label{thm:er-structure} Suppose that $P(\ev)$ is effectively
resonant at $\Rp$, with normal form $P_{\norm}$ near $\Rp$ as in
Lemma~\ref{effnonres-form}, and \eqref{eq:im-b-form-sa} holds.

\begin{enumerate}
\item
If in addition $\Rp$ is a source/sink of $W$, then any microlocal
solution $u$ of $P_{\norm}$ has the form 
\eqref{eq:homog-summand-1-c}, and conversely
given any rapidly Schwartz sequence of functions $w_k$ there is
a microlocally outgoing solution $u$ of $P_{\norm}$ which has the form 
\eqref{eq:homog-summand-1-c}. Thus,
microlocal eigenfunctions at a source/sink are parameterized by Schwartz
functions of the variables $(Y'', Y''').$

\item
If $\Rp$ is a saddle point of $W$, then all microlocal
solutions are in $x^{-1/2 + \epsilon} L^2$ for some $\epsilon > 0$.
For each monomial in the variables $Y'$, each $k \in \NN$ and each
Schwartz function $w(Y'')$ there is a microlocally outgoing solution of the form 
\begin{equation}
u = x^{-i\tilde b- i\kappa_{k}} e^{-i\cP^\sharp_0}
(Y')^{\beta'}w(Y'')v_{k}(Y''')
+u', 
\end{equation}
where $u'$ is in a strictly faster decaying weighted $L^2$ space than $u$,
and every microlocally outgoing solution is a sum of such solutions, with
the $w=w_{k,\beta'}$ rapidly decreasing as $k\to\infty$ in every seminorm.
\end{enumerate}
\end{thm}


\section{Fr{o}m microlocal to approximate eigenfunctions}\label{sec:micro}

We are interested in the structure of (global) eigenfunctions of
$\Delta+V.$ While in the first half of the paper a rather general element
$P\in\Psisc^{*,-1}(X)$ was considered, from now on attention is limited to
\begin{equation*}
H=\Delta+V\in\Psisc^{*,0}(X),\ H(\ev)=H-\ev,
\end{equation*}
in particular {\em the order of $H$ at $\pa X$ is $0.$}

In the next section we obtain an iterative description of the
`smooth' eigenfunctions in terms of the microlocal eigenspaces. As the
first step, we show that if $\Rp$ is a radial point for $H(\ev)=H-\ev,$
then elements of $E_{\mic,\out}(\Rp,\ev),$ which are the
microlocally outgoing eigenfunctions near $\Rp,$ have representatives
satisfying $(H-\ev)u\in\dCinf(X),$ \ie they extend to approximate
eigenfunctions, with $\WFsc(u)$ a subset of the forward flow-out of $\Rp.$
Stated explicitly this is

\begin{prop}\label{prop:min-ident-8}
If $\Rp \in \RP_+(\ev)$ then every element of $E_{\mic,\out}(\Rp,\ev)$
has a representative $u$ such that
$(H-\ev)u\in\dCI(X),$ and $\WFsc(u)\subset \Phi_+(\{\Rp \})$.
\end{prop}

\begin{rem} Fr{o}m this result, given $u$ as in Proposition~\ref{prop:min-ident-8} it is easy to produce an exact eigenfunction $v$ such that $\WFsc(v) \cap \{ \nu \geq 0 \} \subset \Phi_+(\{\Rp \})$: we simply take $v = u - R(\ev - i0)(H - \ev)u$.
\end{rem}

The key ingredient of the proof, as in the two-dimensional case studied in
\cite{HMV1}, is the microlocal solvability of the eigenequation through
radial points. To avoid a microlocal construction along the lines of
H\"ormander \cite{MR48:9458}, we introduce, as in \cite[Lemma~5.3]{HMV1},
an operator
$\tilde H$ which arises from $H$ by altering $V$ appropriately. This is
chosen to be equal to $H$ near the radial point in question but to have no
other radial points in $\RP_+(\ev)$ at which $\nu$ takes a smaller value.
One may then assume, in any argument concerning $q \in \RP_+(\ev),$ that
there is no $q'\in \RP_+(\ev)$ with $\nu(q') < \nu(q).$

As in \cite[Definition~11.3]{HMV1},
we introduce a partial order on $\RP_+(\ev)$
corresponding to the flow-out under $W.$

\begin{Def}\label{Def:partial-order} If $\Rp,$ $\Rp'\in\RP_+(\ev)$ we say
that $\Rp\leq \Rp'$ if $\Rp'\in\Phi_+(\{\Rp\})$ and $\Rp<\Rp'$ if
$\Rp\leq\Rp'$ but $\Rp'\neq\Rp.$ A subset $\Gamma\subset\RP_+(\ev)$ is
closed under $\leq$ if, for all $\Rp\in\Gamma,$ $\{ q' \in \RP_+(\ev); q
\leq q' \} \subset \Gamma.$ We call the set $\{ q' \in \RP_+(\ev); q \leq
q' \}$ the string generated by $q.$
\end{Def}

\begin{rem}\label{rem:bich} This partial order relation between two radial
points in $\RP_+(\ev)$ corresponds to the existence of a sequence
$\Rp_j\in\RP_+(\ev),$ $j=0,\ldots,k,$ $k\geq 1$, with $\Rp_0=\Rp,$
$\Rp_k=\Rp'$ and such that for every $j=0,\ldots,k-1$, there is a
bicharacteristic $\gamma_j$ with $\lim_{t\to-\infty}\gamma_j=\Rp_j$ and
$\lim_{t\to+\infty}\gamma_j=\Rp_{j+1}.$
\end{rem}

\begin{lemma}\label{lemma:Vt} Given $\ev > \min \Vy$ and $\tilde\nu>0,$ set
$K=\Vy^{-1}((-\infty,\ev-\tilde\nu^2])\subset \pa X$ then there exists a
potential function $\tilde V\in\Cinf(X)$ with $\tilde\Vy$ Morse such that
\begin{enumerate}
\item \label{bigger}
$\tilde \Vy\geq \Vy,$
\item
$\tilde V_0 = V_0$  on a neighbourhood of $K$,
\item
no critical value
of $\tilde V$ lies in the interval $(\ev-\tilde\nu^2,\ev]$, 
\item\label{same}
if $\tilde\Sigma(\ev)$ is the characteristic variety at energy $\ev$ of
$\tilde H=\Lap+\tilde V$ then
\begin{equation}
\Sigma(\ev)\cap\{\nu\geq\tilde\nu\}
=\tilde\Sigma(\ev)\cap\{\nu\geq\tilde\nu\}\Mand
\label{eq:Sigma=Sigmat}\end{equation}
\item\label{no-ev}
$\tilde H - \ev$ has no $L^2$ null space.
\end{enumerate}
\end{lemma}

\begin{proof} Choose a smooth function $f$ on the real line so that $f'>0,$
$f(t)=t$ if $t\leq \ev-\tilde\nu^2$ and $f(t)>\ev$ for
$t\geq\min\{V(\Cp);dV(\Cp)=0 \Mand V(\Cp)>\ev-\tilde\nu^2\}>
\ev-\tilde\nu^2.$ Then let $\tilde V=f\circ V,$ so the critical points of
$\Vy$ and $\tilde \Vy$ are the same and are non-degenerate.

On $\Sigma(\ev)\cap\{\nu\geq\tilde\nu\},$ $\nu^2+|\mu|_y^2+\Vy=\ev$, hence
$\Vy\leq \ev-\tilde\nu^2,$ so $\Vy=\tilde \Vy$, and therefore
$\Sigma(\ev)\cap\{\nu\geq\tilde\nu\}\subset\tilde\Sigma(\ev).$ With the
converse direction proved similarly, \eqref{bigger}  -- \eqref{same}
follow. Property \eqref{no-ev} can be arranged by a suitable perturbation of
$\tilde V$ with compact support in the interior.
\end{proof}

These properties of $\tilde H$ are exploited in the proof of the following
continuation result.

\begin{lemma}[Lemma 5.5 of \cite{HMV1}]\label{HMV.136} Suppose $u\in\CmI(X)$
satisfies 
\begin{equation*}
\WFsc(u)\subset\{\nu\geq\nu_1\}\Mand\WFsc((H-\ev)u)\subset\{\nu\geq\nu_2\},
\end{equation*}
for some $0<\nu_1<\nu_2,$ then there exists $\tilde u\in\CmI(X)$ with
$\WFsc(u-\tilde u)\subset\{\nu\geq\nu_2\}$ and $(H-\ev)\tilde u\in\dCI(X).$
\end{lemma}

\begin{proof} We just sketch the proof here; for full details, see
\cite{HMV1}. The obvious idea of subtracting $R(\ev+i0)((H-\ev)u)$ from $u$
does not quite work, since the forward flowout of other critical points
$q' \in \RP_+(\ev)$ with $\nu(q')$ less than $\nu(\Rp)$ may strike $\Rp.$ To avoid
this problem, choose $\tilde \nu$ with $\nu_1 < \tilde \nu < \nu_2,$
sufficiently close to $\nu_2$ so that there are no radial points $\Rp$ with
$\nu(\Rp) \in [\tilde \nu, \nu_2)$, and a corresponding $\tilde V$ as in
Lemma~\ref{lemma:Vt}. Then consider the function $\tilde R(\ev + i0) (H -
\ev) Au$, where $A$ is equal to the identity
microlocally on $\{ \nu \leq \tilde \nu \}\cap\Sigma(\ev)$
and vanishes microlocally in $\{ \nu \geq \nu_2 \}.$ Since
$\tilde V_0$ has no critical points $\Rp$ with $0 < \nu(\Rp) < \nu_2$ it
follows readily $\tilde u = Au - \tilde R(\ev + i0) (H - \ev) Au$ satisfies
the desired conditions.
\end{proof}

{F}rom this we can readily deduce 

\begin{lemma}\label{lemma:min-ident-8} If $\Rp \in \RP_+(\ev)$ then every
element of $E_{\mic,\out}(\Rp,\ev)$ has a representative $\tilde u$ such that
$(H-\ev)\tilde u\in\dCI(X)$ and $\WFsc(\tilde u)$ is contained in the
union of $\Phi_+(\{ \Rp \})$ and the $\Phi_+(\{\Rp'\})$ for those $\Rp'
\in \RP_+(\ev)$ with $\nu(\Rp') > \nu(\Rp).$
\end{lemma}

\begin{proof} Let $O$ be a $W$-balanced neighbourhood of $\Rp$ (see
Definition~\ref{HMV.92}). Let $A\in\Psisc^{-\infty,0}(X)$ be microlocally
equal to the identity on $\Phi_+(\{ \Rp \})\cap \overline{O}$ and supported
in a small neighbourhood of $\Phi_+(\{ \Rp \})\cap \overline{O}.$ Then
there exists $\nu_2 > \nu(\Rp)$ such that $\nu > \nu_2$ on $\Phi_+(\{ \Rp
\})\setminus O$, and $\WFscp(A)\setminus O \subset \{ \nu \geq \nu_2 \}$. (Here $\WFscp(A)$ is the operator wavefront set of $A$, i.e. the complement in ${}^{\scl} T^*_{\partial X} X$ of the set where $A$ is microlocally trivial; see \cite{Melrose43}.)
Now let $u$ be any representative. Since $\WFsc(u)\cap O\subset \Phi_+(\{
\Rp \}),$ $\WFsc(Au-u)\cap O=\emptyset.$ In addition,
$\WFsc(Au)\subset\WFscp(A)\cap\WFsc(u),$ hence $\nu\geq\nu(\Rp)$ on
$\WFsc(Au).$ Moreover, $\WFsc(Au-u)\cap O=\emptyset$ implies that
\begin{equation*}
\WFsc((H-\ev)Au)\cap O=
\WFsc((H-\ev)Au-(H-\ev)u)\cap O=\emptyset,
\end{equation*}
so $\WFsc((H-\ev)Au)\subset\WFscp(A)\setminus O,$ hence is contained in $\{
\nu \geq \nu_2 \}.$ Then, by Lemma~\ref{HMV.136}, there exists $\tilde
u\in\dist(X)$ such that $\nu\geq\nu_2$ on $\WFsc(\tilde u-Au)$ and
$(H-\ev)\tilde u\in\dCI(X)$.  In particular, $\nu\geq\nu(\Rp)$ in
$\WFsc(\tilde u).$ Moreover, $\nu\geq\nu_2$ on $\WFsc(\tilde u-u)\cap O,$
hence by Lemma~\ref{lemma:E-mic-structure-8}, $\WFsc(\tilde u-u) \cap
O=\emptyset,$ so $\tilde u$ and $u$ have the same image in
$E_{\mic,\out}(O,\ev).$
\end{proof}

Finally, we can show that each microlocally outgoing eigenfunction is
represented by an approximate eigenfunction.

\begin{proof}[Proof of Proposition~\ref{prop:min-ident-8}]
Let $\tilde u$ be a representative as in Lemma~\ref{lemma:min-ident-8}. If
we choose $\Rp'$ from the set
\begin{equation}
\{ \Rp' \in \RP_+(\ev) \cap \WFsc(\tilde u); \nu(\Rp') > \nu(\Rp), \
\Rp' \notin \Phi_+(\{ \Rp \}) \},
\label{eq:removing}\end{equation}
with $\nu(\Rp')$ minimal, then, localizing $\tilde u$ near $\Rp',$ gives
an element $v$ of $E_{\mic,+}(\Rp').$ By subtracting from $\tilde u$ a
representative of $v$ given by Lemma~\ref{lemma:min-ident-8},
we remove the wavefront set near $\Rp'.$ Inductively choosing radial points
from \eqref{eq:removing} and performing this procedure repeatedly, all
wavefront set may be removed from $\tilde u$ except that contained in
$\Phi_+(\{ \Rp \}).$
\end{proof}


\section{Microlocal Morse decomposition}\label{sec:mmd}

Next we show that global smooth eigenfunctions can, in an appropriate
sense, be decomposed into components originating, in the sense of the Introduction, at a single radial
point. We do this by defining subspaces of $E^\infty_{\ess}(\ev)$
corresponding to the location of scattering wavefront set in $\{ \nu > 0
\}$ and showing that suitable quotients of these spaces are isomorphic to
the spaces of microlocal eigenfunctions $E^\infty_{\mic, +}(q,\ev),$ $q
\in \RP_+(\ev),$ analyzed in sections~\ref{sec:enr} and
\ref{sec:er}. Since each of the spaces $E^\infty_{\mic, +}(q, \ev),$ $q
\in \RP_+(\ev),$ is non-trivial this shows that each such radial point
gives rise to eigenfunctions. However, as noted previously in \cite{Herbst-Skibsted1}, \cite{Herbst-Skibsted2}, \cite{Herbst-Skibsted3} and
\cite{HMV1} in some special cases, there is a qualitative difference
between the radial points corresponding to local minima of $V_0$ and the
others. This is expressed by Proposition~\ref{density} where we show that
the eigenfunctions $u\in E^\infty_{\Min,+}(\ev)$ originating only at
minimum radial points are dense in $E^0_{\ess}(\ev)$ (definitions of these
spaces are given below).

Recall from \cite[Equation~(3.14)]{HMV1}
the spaces of eigenfunctions of fixed growth
\begin{equation}
E_{\ess}^s(\ev) = \{ u \in E_{\ess}^{-\infty}(\ev); \WFsc^{0,s-1/2}(u) \cap
\{ \nu = 0 \} = \emptyset  \}.
\label{Edefn}\end{equation}
This condition is equivalent to requiring that 
\begin{equation}
Bu \in x^{s - 1/2} L^2(X)
\label{B-op}\end{equation} 
for some pseudodifferential operator $B \in \Psisc^{0,0}(X)$ with boundary
symbol which is elliptic on $\Sigma(\ev) \cap \{ \nu = 0 \}$ and
microsupported in $\{ |\nu| < a(\ev) \},$ where
\begin{equation*}
a(\ev) = \min \{ |\nu(q)|; q \in \RP(\ev) \}.
\end{equation*}
The space $E_{\ess}^0(\ev)$ is of particular interest.  Choose an
operator $A \in \Psisc^{0,0}(X)$ whose boundary symbol is $0$ for $\nu \leq
-a(\ev)$ and $1$ for $\nu \geq a(\ev)$. The space
$E_{\ess}^0(\ev)$ is a Hilbert space with norm
\begin{equation}
\| u \|_{E_{\ess}^0(\ev)}^2 = \langle i[H, A]u, u \rangle.
\label{E0-norm}\end{equation}
The positive-definiteness of this form, and its independence of the choice
of operator $A,$ was shown in \cite{HMV1}, Section 12. An equivalent norm
is  
\begin{equation}
\| Bu \|_{x^{-1/2} L^2} + \| u \|_{x^{-1/2 - \epsilon} L^2}
\label{equiv-norm}\end{equation}
where $\epsilon > 0$ and $B$ is as in \eqref{B-op}; see \cite{HMV1}, section 3. 

We now define subspaces of $E_{\ess}^s(\ev)$ depending on the location
of the scattering wavefront set inside $\{ \nu > 0 \}.$ Given any
$\leq$-closed subset $\Gamma$ of $\RP_+(\ev)$, we define
\begin{equation}
E^s_{\ess}(\ev, \Gamma) = \{ u \in E_{\ess}^s(\ev); \WFsc(u)
\cap \RP_+(\ev) \subset \Gamma \}.
\label{EGamma-defn}\end{equation}
The set of radial points $q \in \RP_+(\ev)$ lying above local minima of
$V$ is an example of a $\leq$-closed subspace and  will be denoted
$\Min_+(\ev)$. In this case we use the notation
$$
E^s_{\Min,+}(\ev) \equiv E^s_{\ess}(\ev, \Min_+(\ev)) = \{ u
\in E_{\ess}^s(\ev); \WFsc(u) \cap \RP_+(\ev) \subset
\Min_+(\ev) \}
$$
to be consistent with \cite{HMV1}.

\begin{prop}\label{prop:Gamma}
Suppose that $\Gamma\subset\RP_+(\ev)$ is $\leq$-closed and $\Rp$ is a
$\leq$-minimal element of $\Gamma$. Then with $\Gamma '=\Gamma
\setminus\{\Rp\}$ 
\begin{equation*}\xymatrix{ 0 \ar[r]^{}&
E_{\ess}^\infty(\ev,\Gamma')\ar[r]^{\iota}&
E_{\ess}^\infty(\ev,\Gamma)\ar[r]^{r_\Rp}&
E_{\mic,\out}(\ev,\Rp) \ar[r]^{}& 
0}
\end{equation*}
is a short exact sequence, where $\iota$ is the inclusion map and $r_\Rp$ is
the microlocal restriction map.
\end{prop}

\begin{proof} The injectivity of $\iota$ follows from the definitions. The
null space of the microlocal restriction map $r_{\Rp},$ which can be viewed
as restriction to a $W$-balanced neighbourhood of $\Rp,$ is precisely the
subset of $E_{\ess}^\infty(\ev,\Gamma)$ with wave front set disjoint
from $\{\Rp\},$ and this subset is $E_{\ess}^\infty(\ev,\Gamma').$ Thus
it only remains to check the surjectivity of $r_{\Rp}.$

We do so first for the strings generated by $q \in \RP_+(\ev)$. For $q \in
\Min_+(\ev)$, the string just consists of $q$ itself and the result follows
trivially. So consider the string $S(q)$ generated by $q \in \RP_+(\ev) \setminus \Min_+(\ev)$. 
By Proposition~\ref{prop:min-ident-8} any element of
$E_{\mic,+}(q,\ev)$ has a representative $\tilde u$ satisfying
$(H-\ev)\tilde u\in\dCI(X)$ with $\WFsc(\tilde
u)\subset\Phi_+( \{ \Rp \} ).$ Then $u = \tilde u - R(\ev - i0)(H -
\ev) \tilde u \in E^\infty_{\ess}(\ev, \Gamma)$, which gives
surjectivity in this case.

For any $\leq$-closed set $\Gamma$ and $\leq$-minimal element $\Rp$,
the string $S(q)$ is contained in $\Gamma$, so the surjectivity of $r_{\Rp}$
follows in general.
\end{proof}

Notice that we can always find a sequence 
$\emptyset=\Gamma_0\subset\Gamma_1\subset\ldots\subset\Gamma_n=\RP_+(\ev),$
of $\leq$-closed sets with
$\Gamma_j\setminus\Gamma_{j-1}$ consisting of a single point $\Rp_j$
which is $\leq$-minimal in $\Gamma_j$: we simply order the $q_i \in
\RP_+(\ev)$ so that $\nu(q_1) \geq \nu(q_2) \geq \dots$, and set $\Gamma_i =
\{ q_1, \dots, q_i \}$. Then Proposition~\ref{prop:Gamma} implies the
following

\begin{thm}[Microlocal Morse Decomposition]\label{thm:mmd}
Suppose that
$\emptyset=\Gamma_0\subset\Gamma_1\subset\ldots\subset\Gamma_n=\RP_+(\ev),$
is as described in the previous paragraph. Then
\begin{equation}
\{0\}\longrightarrow E^\infty_{\ess}(\ev,\Gamma_1)\hookrightarrow \ldots
\hookrightarrow
E^\infty_{\ess}(\ev,\Gamma_{n-1})\hookrightarrow E^\infty_{\ess}(\ev),
\end{equation}
with
\begin{equation}
E^\infty_{\ess}(\ev,\Gamma_{j}) / E^\infty_{\ess}(\ev,\Gamma_{j-1})
 \simeq E_{\mic,+}(\Rp_j,\ev),\ j=1,2,\ldots,n.
\end{equation}
\end{thm}


\section{$L^2$-parameterization of the generalized eigenspaces}\label{sec:para}

Recall from Theorem~\ref{thm:enr-structure}, or
Theorem~\ref{thm:er-structure} in the effectively resonant case, that there
is a surjective map
\begin{multline}
M_+(\ev) : E^{\infty}_{\Min,+}(\ev) \to \oplus_{q \in \Min_+(\ev)}
\mathcal{S}(\bbR^{n-1}), \\ \ev \in (\min V_0, \infty) \setminus \big( \Cv(V) \cup \cup_{z \in \Cv(V)} \cR_{\hesst, z} \big),
\label{Mplus}\end{multline}
given by taking $u \in E^{\infty}_{\Min,+}(\ev)$, microlocally restricting
$u$ to a neighbourhood of each $q$ giving $u_q \in E^\infty_{\mic,+}(\ev,
q)$ and sending $u$ to the sum of the leading coefficients $\sum_k
w_k(Y'')v_k(Y''')$, $(Y'',Y''') \in \bbR^{n-1},$ of each of the
$u_q$. Since the $v_k$ are normalized eigenfunctions of a harmonic
oscillator and the $w_k$ are Schwartz functions of $Y''$ with seminorms
rapidly decreasing in $k$, the sum is a Schwartz function of $(Y'', Y''')$.

Let us regard $\oplus_q \mathcal{S}(\bbR^{n-1})$ as a subspace of $\oplus_q
L^2(\bbR^{n-1})$, endowed with the norm
\begin{equation}
\| (w_q)_{q \in \Min_+(\ev)} \|^2 = \sum_q
\int_{\bbR^{n-1}} |w_q(Y)|^2 d\omega_{q,\sigma},\ d\omega_{q,\sigma}
= 2\sqrt{\ev - V(\pi(q))}\,d\omega_q,
\label{Hilbert-norm}\end{equation}
where $\omega_q$ is the measure induced by Riemannian measure, namely the
measure $$x^{n - (n-m)/2 - \sum_j r''_j} dg$$ divided by $dx/x$ and
restricted to $x=0$. (It takes the form $dY'' \, dY'''$ provided that the
$y$ are normal coordinates, centred at the critical point, for the metric
$h(0, y, dy)$.)

The next result is the main content of this section.

\begin{thm}\label{thm:iso} The map $M_+(\ev)$ in \eqref{Mplus} has a unique
  extension to an unitary isomorphism
\begin{equation*}
M_+(\ev) : E^0_{\ess}(\ev) \to \oplus_{q \in \Min_+(\ev)} L^2(\bbR^{n-1}).
\end{equation*}
\end{thm}

\begin{rem} Here, and throughout this section, we take $\ev \in (\min V_0,
\infty) \setminus \Cv(V)$.
\end{rem}

To prove the theorem, we establish several intermediate results. First we show

\begin{prop}\label{density} The space $E^\infty_{\Min,+}(\ev)$ is dense in
  $E^\infty_{\ess}(\ev)$ in the topology of $E^0_{\ess}(\ev).$
\end{prop}

\begin{proof} The proof is by induction. We consider a sequence $\Gamma_0
\subset \Gamma_1 \subset \dots \subset \Gamma_n = \RP_+(\ev)$ as in the
previous section, but with the additional condition that the radial points
are ordered so that, among the points with equal values of $\nu$, those
corresponding to local minima of $V_0$ are placed last. We shall prove by
induction that
\begin{equation}
E^\infty_{\ess}(\ev, \Gamma_i \cap \Min_+(\ev)) 
\text{ is dense in } E^\infty_{\ess}(\ev, \Gamma_i) \text{ in the topology
  of } E^0_{\ess}(\ev).
\label{ind}\end{equation}
For $i=1$ there is nothing to prove. Assume that \eqref{ind} is true for
$i=k.$ Let $\Gamma_{k+1} \setminus \Gamma_k = \{ q \}.$ If $q$ arises from
a local minimum of $V_0,$ then using a microlocal decomposition, any $u \in
E^\infty_{\ess}(\ev, \Gamma_{k+1})$ can be written as the sum of $u_1 \in
E^\infty_{\ess}(\ev, \{ q \})$ and $u_2 \in E^\infty_{\ess}(\ev,
\Gamma_{k})$. A similar statement is true for $u \in E^\infty_{\ess}(\ev,
\Gamma_{k+1} \cap \Min_+(\ev))$, which proves \eqref{ind} for $i = k+1$.

Next suppose that $q$ does not arise from a local minimum of $V_0$. Then we
adapt the argument of Proposition 11.6 of \cite{HMV1} to prove \eqref{ind}
for $i=k+1$. We first make the assumption that $\ev$ is not in the point
spectrum of $H$. Using our inductive assumption, it is enough to show that
$E^\infty_{\ess}(\ev, \Gamma_{k})$ is dense in $E^\infty_{\ess}(\ev,
\Gamma_{k+1})$. Let $u \in E^\infty_{\ess}(\ev, \Gamma_{k+1})$.  Let $Q \in
\Psisc^{0,0}(X)$ be microlocally equal to the identity near $\Gamma_k \cap
\Min_+(\ev)$, and microsupported sufficiently close to $\Gamma_k \cap
\Min_+(\ev)$. Then away from $\Min_+(\ev)$, $u \in x^{-1/2 + \epsilon} L^2$
by (ii) of Theorem~\ref{thm:enr-structure} and thus $(H- \ev) Q u = [H, Q]
u \in x^{1/2 + \epsilon} L^2$ for some $\ep>0$. This is also true near
$\Min_+(\ev)$ since $Q$ is microlocally the identity there, so we have $(H-
\ev) Q u \in x^{1/2 + \epsilon} L^2$ everywhere. This implies that
\begin{equation}
u = Qu - R(\ev - i0)(H - \ev) Qu,
\label{u-eqn}\end{equation}
since $v = u - (Qu - R(\ev - i0)(H - \ev) Qu)$ satisfies $(H- \ev) v = 0$
and $v \in x^{-1/2 + \epsilon} L^2$ microlocally for $\nu > 0.$ 

Now choose a modified potential function $\tilde V$ as in
Lemma~\ref{lemma:Vt}, where we choose $\tilde \nu$ larger than $\nu(q)$ but
smaller than $\nu(q')$ for every $q' \in \Gamma_k \cap \Min_+(\ev)$. (This
is possible because of the way we ordered the $q_i$.) Since $\WFsc(Qu)$
lies in $ \{ \nu > \tilde \nu \}$, we have
\begin{equation}
Qu = \tilde R(\ev + i0) (\tilde H - \ev) Qu.
\label{uid2}\end{equation}

Now take $u'_j = \phi(x/r_j) u$, where $\phi \in \CI(\bbR)$, $\phi(t) = 1$
for $t \geq 2$, $\phi(t) = 0$ for $t \leq 1$ and $r_j \to 0$ as $j \to
\infty$. Then $u'_j \in \CIdot(X)$, and $w_j$ defined by
$$
w_j = \tilde R(\ev + i0) (\tilde H - \ev) Qu'_j
$$
converge to $Qu$ in $x^{-1/2 - \epsilon} L^2$. Our choice of $\tilde V$
ensures that
$$
\WFsc(w_j) \cap \RP_+(\ev) \subset \Gamma_k.
$$
Moreover, 
\begin{equation}
(H - \ev) w_j \text{ converges to } (H - \ev) Qu \text{ in } x^{1/2 +
    \epsilon} L^2.
\label{wj-conv}\end{equation}
Now define
$$
u_j = w_j - R(\ev - i0)(H - \ev) w_j.
$$
Then $u_j \in E^\infty_{\ess}(\ev, \Gamma_k)$. We claim that $u_j \to u$ in
the topology of $E^0_{\ess}(\ev)$. Certainly, $u_j \to u$ in $x^{-1/2 -
\epsilon} L^2$. We must also show that $Bu_j \to Bu$ in $x^{-1/2} L^2$,
where $B$ is as in \eqref{B-op}.  To do this we write
\begin{equation*}\begin{gathered}
B u_j - B u = B \big( w_j - R(\ev - i0)(H - \ev) w_j \big) - B \big( 
(\Id - Q) u + Qu \big) \\
= B \Big( \tilde R(\ev + i0) (\tilde H - \ev) Qu'_j - R(\ev - i0)(H - \ev)
w_j \\ + R(\ev - i0)(H - \ev) Qu - \tilde R(\ev + i0) (\tilde H - \ev) Qu
\Big),
\end{gathered}\end{equation*}
using \eqref{u-eqn} and \eqref{uid2}, and this goes to zero in $x^{-1/2}
L^2$ by \eqref{wj-conv} and propagation of singularities, Theorem 3.1 of
\cite{HMV1}, as in the proof of \cite[Proposition~11.6]{HMV1}.

If $\ev$ is in the point spectrum of $H$, then equation \eqref{u-eqn} must
be replaced by
\begin{equation}
u = \Pi \Big( Qu - R(\ev - i0)(H - \ev) Qu \Big),
\label{u-eqn-2}\end{equation}
where $\Pi$ is projection off the $L^2$ $\ev$-eigenspace. Consequently we
must define $w_j$ by $\Pi \tilde R(\ev + i0) (\tilde H - \ev) Qu'_j$, and
then the rest of the proof goes through.
\end{proof}

The second intermediate result we need is
\begin{prop}\label{prop:bdy-pairing} The Hilbert norm \eqref{E0-norm} on
  the subspace $E^\infty_{\Min,+}(\ev) \subset E^0_{\ess}(\ev)$ is given by
the formula
\begin{equation}
\| u \|_{E_{\ess}^0(\ev)}^2 = \sum_{q \in \Min_+(\ev)} 2\sqrt{\ev - V(\pi(q))} \int_{\bbR^{n-1}} \big|
M^+(q, \ev) u \big|^2 \, d\omega_q .
\label{norm-form}\end{equation}
\end{prop}

\begin{proof} The proof is the same as the one dimensional case, which is
  proved in Proposition 12.6 of \cite{HMV1}, so we just give a sketch here.

Let $\phi$ be as in the proof of Proposition~\ref{density}. Then we can
write the natural norm \eqref{E0-norm} on $E^0_{\ess}(\ev)$ as a limit
\begin{equation}
\lim_{r \to 0} i \langle (H - \ev) A u, \phi(x/r) u \rangle = \lim_{r \to
  0} i \langle A u, [H, \phi(x/r)] u \rangle.
\label{0-pairing}\end{equation}

Since $u \in x^{-1/2 - \epsilon} L^2$, the only term in $[H, \phi(x/r)]$
contributing in the limit is $2 (x^2 D_x) \phi(x/r) (x^2 D_x)$. The
cutoff operator $A$ restricts attention to $\{ \nu > 0 \}$, and the limit
vanishes when localized to any region where $u \in x^{-1/2 + \epsilon}
L^2$, so we can substitute for $u$ a sum of expressions $u_q$ as in
\eqref{eq:homog-summand-111} in the effectively non-resonant case,
or its analogue in the effectively resonant setting arising from
\eqref{eq:homog-summand-1-c} (namely $e^{i\Phi(y)/x}$ times an expression
as in \eqref{eq:homog-summand-1-c}, cf.\ Remark~\ref{original-exp}),
one for each $q \in \Min_+(\ev)$. A
straightforward computation then gives \eqref{norm-form}.  \end{proof}

\begin{proof}[Proof of Theorem~\ref{thm:iso}]
Proposition~\ref{prop:bdy-pairing} shows that $M_+(\ev)$ maps
$E^\infty_{\Min,+}(\ev)$ into a dense subspace of $\oplus_q
L^2(\bbR^{n-1})$, with the Hilbert norm of $M_+(\ev) u$, $u \in
E^\infty_{\Min,+}(\ev)$, equal to that of $u$. By
Proposition~\ref{density}, $E^\infty_{\Min,+}(\ev)$ is dense in
$E^\infty_{\ess}(\ev)$, and by Corollary 3.13 of \cite{HMV1},
$E^\infty_{\ess}(\ev)$ is dense in $E^0_{\ess}(\ev)$. The result follows.
\end{proof}

So far we have only considered the microlocal restriction of eigenfunctions near
radial points $q$ satisfying $\nu(q)>0.$ For each critical point of
$V_0,$ there are two corresponding radial points with opposite signs of
$\nu$, and we can equally well consider microlocal restriction near radial
points with $\nu(q)<0.$ This leads to an operator
$$
M_-(\ev) : E^0_{\ess}(\ev) \to \oplus_{q \in \Min_-(\ev)} L^2(\bbR^{n-1})
$$ 
and the analogue of Theorem~\ref{thm:iso} holds also for $M_-(\ev)$. 

\begin{Def} The inverses of $M_{\pm}(\ev),$ $P_\pm(\ev) : \oplus_{q
  \in\Min_\pm(\ev)} L^2(\bbR^{n-1}) \to E^0_{\ess}(\ev)$ of $M_\pm(\ev)$
  are called the \emph{Poisson operators at energy} $\ev.$ 
\end{Def}

We can identify $\oplus_{q \in \Min_+(\ev)} L^2(\bbR^{n-1})$ and $\oplus_{q
  \in \Min_-(\ev)} L^2(\bbR^{n-1})$ in the obvious way, and may therefore
assume that the $M_\pm(\ev)$ have the same range, identified with the
domain of $P_\pm(\ev).$

\begin{cor} For $\ev \notin \Cv(V)$, 
the S-matrix may be identified as the unitary operator
$S(\ev)=M_+(\ev)P_-(\ev)$ on $\oplus_{z\in\Min}L^2(\bbR^{n-1}).$
\end{cor}

\begin{rem}
For $n=2,$ the structure of $S(\ev)$ was described rather precisely in
\cite{HMV1EDP} as an anisotropic Fourier integral operator.
\end{rem}

Theorem~\ref{thm:iso} is essentially a pointwise version of asymptotic
completeness in $\ev.$ Integrating gives a version of the usual statement,
but some uniformity in $\ev$ is required for this. So we proceed to discuss
an extension of part (i) of Theorem~\ref{thm:enr-structure} that is valid
in an interval rather than just at one value. For this purpose, let $I
\subset (\min V_0, \infty)$ be a compact interval disjoint from the set of
effectively resonant energies, the set of Hessian thresholds and
$\Cv(V)$. Then for each $\ev \in I$, the sets $\Min_+(\ev) \subset
\RP_+(\ev)$ can be identified; we write $\Min_+(I)$ for this set. Each
element of $\Min_+(I)$ is thus a continuous family $q(\ev)$ of minimal
radial points, with $q(\ev) \in \Min_+(\ev)$.

\begin{prop}\label{prop:I-asymp} Let $I \subset (\min V_0, \infty)$ be as
above, and let the $\Rp(\ev) \in \Min_+(I)$ be an outgoing radial point
associated to a minimum point $z$ of $V_0$, with $Y'', Y'''$ the associated
coordinates given by \eqref{eq:def-Phi}. For any $h(\ev, \cdot)\in \CI(I;
\Schwartz(\overline{\bbR^{n-1}}))$ there is $\phi\in\dCinf(X)$ orthogonal to $E_{\pp}(I)$ such that
for every $\ev \in I$,
\begin{equation}\begin{split}\label{eq:I-asymp}
F(\ev)^{-1}R(\ev+i0) \phi= \sum_j x^{-i \tilde b - i\kappa_{j}}
w_j(Y'',\ev)v_{j}(Y''', \ev)
+u', \\
h(\ev,Y'',Y''')=\sum_j w_j(Y'',\ev)v_{j}(Y''', \ev),
\end{split}\end{equation}
where $w_j$,
$v_{j}$, $\kappa_j$ and $\tilde b$ are as in Proposition~\ref{prop:enr-min}, $F(\sigma)$ is as in Theorem~\ref{thm:model-smooth}, and
where $u'\in \Cinf(I;\Isc^{(l)}(X,\cM))$ for some $l>-\frac{1}{2}.$
\end{prop}

\begin{rem}
The statement $u'\in \Cinf(I;\Isc^{(l)}(X,\cM))$ is meant to underline
that this is a global claim, namely $u'\in\Cinf(I;\Isc^{(l)}(O,\cM))$
and that it is $\Cinf$ with values in $\dCinf(X)$ microlocally away from
$\{\Rp(\ev); \ev\in I\},$ \ie for all $A\in\Psisc(X)$ with
$\WFsc'(A)\cap \{\Rp(\ev); \ev\in I\}=\emptyset$, $Au'\in\Cinf(I;\dCinf(X)).$
\end{rem}

\begin{proof} By the construction of Section~\ref{sec:enr}, for each
$\ev\in I$ there is an approximate microlocally outgoing solution $u_\ev$
with $f_\ev=(H-\ev)u_\ev\in\dCI(X)$ and $F(\ev)^{-1}u_\ev$ of the same form
as the right hand side of \eqref{eq:I-asymp}.  Indeed, the construction is
smooth in $\ev$, in the sense that $(d/d\ev)^k u \in I^{s}(O, \mathcal{M})$
for each $k$ and each $s < -1/2,$ so that with $f(\ev,.)=f_\ev(.)$, we have
$f\in\CI(I; \dCI(X)).$ Notice that there is no need to `globalize' using
Proposition~\ref{prop:min-ident-8}, since microlocally outgoing solutions
over sources/sinks (\ie minima of $V_0$) are localized at $\Rp(\ev)$.

Let $\tilde f\in\dCI_c(\Cx\times X)$ be an almost analytic extension of $f$
with compact support, so $\overline{\pa}_{\ev} f$ vanishes to infinite
order at $\bbR\times X$, and let
\begin{equation*}
\phi=\frac{-1}{2\pi i}\int_{\Cx} R(\ev) \overline{\pa}_{\ev} \tilde f\,d\ev
\wedge d\bar\ev.
\end{equation*}
Thus, $\phi\in\dCI(X)$ since $\overline{\pa}_{\ev} \tilde f$ vanishes to
infinite order on the real axis.

We also claim that \eqref{eq:I-asymp} holds. Indeed, let $\ev_0\in \bbR$,
$\chi\in\Cinf_c(\bbR)$, $\chi$ identically $1$
near $\ev_0$, let $\tilde\chi$ be an almost analytic extension of $\chi$
of compact support.
Thus,
\begin{equation*}
f(\ev,.)=f(\ev_0,.)\chi(\ev)+(\ev-\ev_0)g(\ev,.),
\ \tilde f(\ev,.)=f(\ev_0,.)\tilde \chi(\ev)+(\ev-\ev_0)\tilde g(\ev,.)
\end{equation*}
with $g\in\dCI_c(\bbR\times X)$, $\tilde g\in\dCI_c(\Cx\times X)$.
Then, writing $\ev-\ev_0=(H-\ev_0)-(H-\ev)$,
\begin{equation*}\begin{split}
\phi=&\frac{-1}{2\pi i}\left(\int_{\Cx} R(\ev) \overline{\pa}_{\ev} \tilde
\chi\,d\ev\wedge d\bar\ev\right) f(\ev_0,.)\\
&-\frac{1}{2\pi i}(H-\ev_0)\int_{\Cx} R(\ev) \overline{\pa}_{\ev} \tilde
g\,d\ev\wedge d\bar\ev
+\frac{1}{2\pi i}\int_{\Cx} \overline{\pa}_{\ev} \tilde
g\,d\ev\wedge d\bar\ev,
\end{split}\end{equation*}
where in the last term the identity $(H-\ev)R(\ev)=\Id$ is used. Since the
last term vanishes (as $\tilde g$ is smooth), and
the integral in the second term is in $\dCI(X)$, while the integral in the
first term is $\chi(H)$, we deduce that
$$
\phi=f_{\ev_0}+(H-\ev_0)f'_{\ev_0} = (H - \ev_0)(u_{\sigma_0} + f'_{\sigma_0})
$$
for some $f'\in\dCI(I\times X)$. Then if $v \in E_{\pp}(I)$, $(H - \sigma_0) v = 0$, we have $v \in \dCI(X)$, so $\langle \phi, v \rangle = \langle u_{\sigma_0} + f'_{\sigma_0}, (H - \sigma_0)v \rangle = 0$. Also 
$R(\ev_0+i0)\phi-R(\ev_0+i0)f_{\ev_0} = f'_{\ev_0}\in\dCI(X)$, so $R(\ev_0+i0)\phi$ and
$R(\ev_0+i0)f_{\ev_0}$
indeed have the same asymptotics. In particular, \eqref{eq:I-asymp} holds
for every $\ev_0\in\bbR$.
\end{proof}

Now we state asymptotic completeness in a more standard form.

\begin{thm}[Asymptotic completeness]\label{thm:AC}
Let $I \subset (\min V_0, \infty)$ be a compact interval as
above. Then 
\begin{equation*}
M_+(\cdot) \circ \Sp(\cdot):{\operatorname{Ran}}( \Pi_I)\ominus E_{\pp}(I)\to
\oplus_{\Rp\in\Min_+(I)}L^2(I\times \bbR^{n-1}_q;
\,2\pi\,d\sigma\,d\omega_{q,\sigma})
\end{equation*}
is unitary. Here, as before, $d\omega_{q,\sigma}=2\sqrt{\ev-V(\pi(\Rp))}
\,d\omega_q.$
\end{thm}

\begin{proof}
For $f\in\dCI(X)$ orthogonal to $E_{\pp}(I),$ let
\begin{equation*}\begin{split}
u=u(\ev)=(2\pi i)^{-1} (R(\ev+i0)f - R(\ev - i0) f) =  \Sp(\ev) f,\\
\qquad \Mwhere
\Sp(\ev) = (2\pi i)^{-1} (R(\ev+i0) - R(\ev - i0))
\end{split}\end{equation*}
is the spectral measure. 
The norm of $u$ in $E^0_{\ess}(\ev)$ is given by $\langle i (H - \ev) A u, u \rangle$, where $A$ is as in \eqref{E0-norm}. Notice that
\begin{equation*}\begin{gathered}
2\pi i(H - \ev) Au - f = (H - \ev)A \big(R(\ev + i0) - R(\ev - i0)\big) f - (H - \ev)R(\ev + i0) f  \\
= (H - \ev) \Big( (A - \Id) R(\ev + i0) f - A R(\ev - i0) f \Big) = (H - \ev) v, \quad v \in \CIdot(X),
\end{gathered}\end{equation*}
since
\begin{equation*}
\WFsc'(A)\cap\WFsc (R(\ev - i0) f)=\emptyset\Mand
\WFsc'(A - \Id)\cap\WFsc (R(\ev + i0) f)=\emptyset.
\end{equation*}
Hence
$$
2\pi \| u \|_{E^0_{\ess}(\ev)}^2 = 2\pi i \langle  (H - \ev) A u, u \rangle =  \langle f + (H - \ev) v, u \rangle = \langle f, \Sp(\ev) f \rangle.
$$
The right hand side is continuous, hence so is the left hand side.

Integrating over $\ev$ in $I,$ denoting
the spectral projection of $H$ to $I$ by $\Pi_I,$ and using
Proposition~\ref{prop:bdy-pairing}, we deduce that
$M_+(\ev) \Sp(\ev) f$ is continuous with values in $L^2$ and
\begin{equation}
\|\Pi_I f\|^2= 
2\pi \,\int_I \| M_+(\ev) \Sp(\ev) f \|^2 \,d\ev,
\label{eq:Pi_I-isometry}\end{equation}
so $M_+(\cdot) \circ \Sp(\cdot)$ is an isometry on the orthocomplement of the
finite dimensional space $E_{\pp}(I)$ in the range of $\Pi_I.$

It remains to prove that the range is dense in  $\oplus_{\Rp\in\Min}L^2(I\times \bbR^{n-1})$. It suffices to show that if $h\in\oplus_{\Rp\in\Min}\dCI(I\times\overline{\bbR^{n-1}})$,
then there is a
$f\in\dCI(X)$ with $M_+(\ev)\Sp(\ev)f=h(\ev,.)$. But this was proved in
Proposition~\ref{prop:I-asymp}, so the proof of the theorem is complete.
\end{proof}

\begin{rem} The results of this section can be related more closely with
Theorem~\ref{thm:mmd} by considering the closure of $E^{\infty}_{\Min,
+}(\ev)$ as a subset of $E^{\infty}_{\ess}(\ev)$ in the topology of
$E^s_{\ess}(\ev)$ for varying values of $s$. We have seen in
Proposition~\ref{density} that $E^{\infty}_{\Min, +}(\ev)$ is dense, in the
topology of $E^0_{\ess}(\ev)$. In fact the proof of
Proposition~\ref{density} shows that this is true in the topology of
$E^s_{\ess}(\ev)$ for $0 \leq s < s_0$, where $s_0$ is the smallest number
such that every $u \in E^\infty_{\mic}(q)$, for every $q \in \RP_+(\ev)
\setminus \Min_+(\ev)$, is in $x^{-1/2 + s_0}L^2$ locally near $\pi(q)$;
that $s_0$ is strictly positive follows from (ii) of
Theorem~\ref{thm:enr-structure}. By contrast, $E^\infty_{\Min, +}(\ev)$ is
closed in the $E^\infty_{\ess}(\ev)$ topology. What happens as $s$
increases is that the closure of $E^\infty_{\Min, +}(\ev)$ in the
$E^s_{\ess}(\ev)$ topology changes discretely, as $s$ crosses certain
values determined by the structure of eigenfunctions at the non-minimal
critical points.

One way to understand this is in terms of microlocally \emph{incoming}
eigenfunctions at the outgoing radial points, \ie microlocal eigenfunctions
$u$ with scattering wavefront set near $q$ is contained in $\Phi_-(q)$ as
opposed to $\Phi_+(q)$. In Part I we showed (in all dimensions) that there
are nondegenerate pairings
\begin{equation*}\begin{gathered}
E_{\mic, +}(q, \ev) \times E_{\mic, -}(q, \ev) \to \mathbb{C}, \\
E^s_{\ess}(\ev) \times E^{-s}_{\ess}(\ev)  \to \mathbb{C}
\end{gathered}\end{equation*}
(Lemma 12.2 and Proposition 12.3 of \cite{HMV1}). The closure of
$E^\infty_{\Min, +}(\ev),$ in the topology of $E^s_{\ess}(\ev),$ may be
identified with the annihilator, in $E^\infty_{\ess}(\ev)$, of the
eigenfunctions which are in $E^{-s}_{\ess}(\ev)$ and have scattering
wavefront set contained in
$$
\bigcup_{q \in \RP_+(\ev) \setminus \Min_+(\ev)} \Phi_-(q) \cup \{ \nu < 0 \}.
$$
This set is trivial for $s < s_0$, and nontrivial for $s > s_0.$  The fact
that this set of eigenfunctions jumps discretely with $s$ in shown in the
two dimensional case in Section 10 of Part I.
\end{rem}


\section{Time-dependent Schr\"odinger equation}\label{sec:time}

\subsection{Long-time asymptotics}
In this final section we apply the earlier results to deduce the long-time
asymptotics for solutions of the initial value problem
\begin{equation}
(D_t+H)u=0,\ u|_{t=0}=u_0,\ u_0\in\dCinf(X),
\label{time-dept-eqn}\end{equation}
for a dense set (in $L^2 \ominus E_{\pp}(H)$) of initial data. 

Our approach is to use the spectral resolution of $u_0$ and the functional
calculus. In this way, we deduce the long-time asymptotics of $u$ from the
asymptotics of generalized eigenfunctions of $H$ using the stationary phase
lemma.

We first define the space $\Xsch$ on which the asymptotics of the solution
$u$ of \eqref{time-dept-eqn} will be described. Let us first choose a
global boundary defining function $x$ satisfying
\eqref{sc-metric}; we can specify, for example, that $x \equiv 1$ outside a
collar neighbourhood of $\partial X$. We then introduce the variable $\tau
= tx$, where $t$ is time. Let us compactify the real $\tau$-line $\RR$ to
an interval $\Rbar$ using $\tau^{-1}$ as a boundary defining function near
$\tau = \infty$, and $-\tau^{-1}$ as a boundary defining function near
$\tau = -\infty$. Then we define
\begin{equation}
\Xsch = X \times \Rbar_{\tau} 
\label{Xsch}\end{equation}
Thus $\Xsch$ is a compact manifold with corners, with boundary
hypersurfaces if (the `infinity face') at $\tau = \pm\infty$ (or $t = \pm
\infty$), naturally diffeomorphic to two copies of $X$ (one at $t =
+\infty$, one at $t=-\infty$), and a boundary hypersurface af (the
`asymptotic face') diffeomorphic to $\partial X \times \Rbar_\tau$. At af,
every point with $\tau > 0$ corresponds to $t = +\infty$ and 
every point with $\tau < 0$ corresponds to $t = -\infty$, so this is the
place to look for long-time (and large-distance) asymptotics of the
Schr\"odinger wave $u$. The variable $\tau$ has an interpretation of
inverse speed; a particle travelling asymptotically radially at speed
$\tau_0^{-1}$ will end up at af after infinite time at $\tau = \tau_0$.

We now specify a good subset of $L^2$ initial data $u_0$, for which the
asymptotics as $t \to +\infty$ of the solution, $u,$ to
\eqref{time-dept-eqn} are particularly simple. Let $I \subset (\min V_0,
\infty)$ be a compact interval disjoint from $\Cv(V)$ and from the set of
effectively resonant energies and Hessian thresholds. Let $(h(\ev,
\cdot))q \in \CI(I;\Schwartz(\RR^{n-1}))$ be a collection of smooth
functions from $I$ into Schwartz functions of $n-1$ variables, one for each
$q \in \Min_+(I)$, and let $\phi = \phi(I,h) = \sum_q \phi(I, h_q) \in
\CIdot(X)$ be the function constructed in Proposition~\ref{prop:I-asymp}. Let
\begin{equation*}
\mathcal{A}_I = \{ \phi(I,h); h(\ev, \cdot) \in \CI(I;
\Schwartz(\RR^{n-1})) \} \text{ and } \mathcal{A} = \sum_I \mathcal{A}_I
\end{equation*}
be the (algebraic) vector space sum of $\mathcal{A}_I$ over all such $I$ as
above. It is clear from Theorem~\ref{thm:AC} that $\mathcal{A}_I$ is dense
in $\Ran \Pi_I(H) \ominus E_{\pp}(I)$, and hence that $\mathcal{A}$ is
dense in $L^2 \ominus E_{\pp}(H) = H_{\ac}(H)$.  To give the asymptotics of
\eqref{time-dept-eqn} with initial data from $\mathcal{A}$ it suffices to
give the asymptotics starting from $u_0 = \phi(I,h)$ for some $h$ as above.

\begin{thm}\label{thm:enr-structure-t} Suppose that $I$ is as above and that $\phi = \phi(I, h) \in \mathcal{A}_I$.
Let $u(\cdot, t)$ be the solution of \eqref{time-dept-eqn} with initial
data $u_0 = \phi$, regarded as a function on $\Xsch$. Then $u$ has trivial
asymptotics (i.e. $u$ and all its derivatives are $O(t^{-\infty})$) at if. Also, if $w \in \partial X$ is not a local minimum of $V_0$, and $\tau > 0$, then $u$ has trivial asymptotics in a neighbourhood of $(w, \tau) \in \, $af. 

Let $z$ be a local minimum of $V_0$, and let $(Y'', Y''')$ be the coordinates given by \eqref{eq:def-Phi}, where $\sigma$ is determined in terms of $\tau$ by \eqref{energy-eqn}. Then, in a neighbourhood of $(z, \tau) \in \, $af,  $u$ takes the
form
\begin{equation}\begin{gathered}\label{eq:t-asymp}
u(x, Y'',Y''', \tau) = c \tau^{-3/2} 
\sum_j x^{-i \tilde b- i\kappa_{j} + 1/2} e^{i\Psi(y, \tau)/x} w_j(Y'',\ev(\tau))v_{j}(Y''', \ev(\tau))
+u', \\
\text{ where } \ h(\ev(\tau),Y'',Y''')=\sum_j w_j(Y'',\ev(\tau))v_{j}(Y''', \ev(\tau)), \quad
c =  \frac1{2\sqrt{\pi}}  e^{-3i\pi/4},  
\end{gathered}\end{equation} 
\begin{equation}
\ev(\tau)=V_0(z)+\frac{1}{4\tau^2}, 
\label{energy-eqn}\end{equation}
$\tilde b$ is as in \eqref{eq:def-exps} and \eqref{eq:im-b-form-sa}, $\kappa_j$ is as in \eqref{eq:homog-summand-1}, 
$\Psi$ is a smooth function of $y$ and $\tau$, $h$ is decomposed as in
Proposition~\ref{prop:I-asymp}, and $u'$ decays faster than the leading term.
\end{thm}

\begin{proof}
Let $v(\ev) = \Sp(\ev) \phi = (2\pi i)^{-1}(R(\ev + i0) - R(\ev - i0))
\phi$. Then
\begin{equation*}
u(t, \cdot) = \frac1{2\pi i} \int_I e^{-it\ev} (R(\ev + i0) - R(\ev - i0))
\phi \, d\ev. 
\end{equation*}
Shifting the contour of integration shows that, as $t \to \infty$, $R(\ev -
i0) \phi$ has trivial asymptotics. Hence it is enough to consider
\begin{equation}
u(t, \cdot) = \frac1{2\pi i} \int_I e^{-it\ev} R(\ev + i0) \phi \, d\ev.
\label{spec-int}\end{equation}

Let $F(\ev)$ be the FIO constructed in Theorem~\ref{thm:model-smooth},
which conjugates $x^{-1} (H - \sigma)$ to normal form microlocally near the
point $q \in \RP_+(\sigma)$ with $\pi(q) = z.$ 
By construction, $F({\ev})^{-1} R(\ev + i0) \phi$ has asymptotics
\eqref{eq:I-asymp} for every $\ev$. Since $F(\ev)$ is a smooth family of
FIOs, and \eqref{eq:I-asymp} is a Legendre distribution associated to the
zero section (i.e. it is conormal at $x=0$ with no oscillatory factor),
 it follows that $R(\ev + i0) \phi$ itself has asymptotics
\begin{equation}
R(\ev + i0) \phi = x^{-i\tilde b - i \kappa_j} e^{i\Phi(y, \ev)/x} a(Y'', Y''', x, \ev) + v',
\label{res-exp}\end{equation}
where $\Phi(\cdot, \sigma)$ is a smooth function, parametrizing the image of the zero section under the canonical relation of $F(\sigma)$ (as in \eqref{eq:homog-summand-111}). 
By assumption, 
$a$  is smooth in $\ev$, conormal in $x$ and Schwartz in $(Y'', Y''')$.  At the critical point $z$ we have
\begin{equation}
\Phi(z, \ev) = \sqrt{\ev - V_0(z)},\quad  z = \pi(q), \quad q \in \Min_+(\ev). 
\label{Phi-value}\end{equation}
We may substitute
\eqref{res-exp} into \eqref{spec-int} and compute
\begin{equation}
u(t, \cdot) = \frac1{2\pi i} \int_I e^{-it\ev} \Big( e^{i\Phi(y, \ev)/x}
a(Y'', Y''', x, \ev) + v' \Big) \, d\ev ,
\label{spec-int-2}\end{equation}
exploiting the smoothness of $\Phi$ and $a$ in $\ev$.

Let $p \in X$ be an interior point. Then $(R(\ev \pm i0) \phi)(p)$ is a
smooth function of $\ev$ by Proposition~\ref{prop:I-asymp}.
It follows that for a fixed interior point $p$ the integral
\eqref{spec-int-2} is rapidly decreasing as $t \to \infty$, being the
Fourier transform of a smooth, compactly supported function. Hence the asymptotics of $u$ are trivial at if.

To investigate asymptotics at af, where $x \to 0$, we rewrite
\eqref{spec-int-2} as
\begin{equation}
u(\tau, x, Y'', Y''') = \frac1{2\pi i} \int e^{i(-\tau\ev + \Phi(y,\ev))/x}
\Big( a(Y'', Y''', x, \ev) + v'(Y'',Y''',x, \ev) \Big) \, d\ev ,
\label{spec-int-3}\end{equation}
and apply stationary phase to the integral. We first note that for any $w \in \partial X$ which is not a local minimum of $V_0$, the integrand is rapidly decreasing as $x \to 0$ in a neighbourhood of $(w, \ev) \in \partial X \times I$, uniformly in $\ev$, so $u$ is rapidly decreasing as $x \to 0$ in a neighbourhood of $(w, \tau) \in$ af. So we may restrict attention to a neighbourhood of $(z, \tau) \in$ af, where $z$ is a local minimum of $V_0$. 

To do this we apply stationary phase to \eqref{spec-int-3}. 
The phase is critical when $\tau = d_{\ev} \Phi(y, \ev)$. Since $\Phi$ is smooth in $y$, this gives 
$$
\tau_{\operatorname{critical}} = d_{\ev}\Phi(z, \ev) + O(Y_i x^{r_i})
$$
and, since $a$ is Schwartz in $Y$, to compute the expansion of $u$ to leading order we may drop the $O(Y_i x^{r_i})$ terms when we substitute $\tau = \tau_{\operatorname{critical}}$ into $a$ in \eqref{spec-int-2}. Since $\Phi(z, \ev)$ is given by \eqref{Phi-value}, we may therefore take $\tau$ in the argument of $a$ to be given by $\tau = d_{\ev} \Phi(z, \ev)$ which implies \eqref{energy-eqn}. 
Moreover, the Hessian of the phase function at the critical point is $(4x)^{-1} (\ev - V_0(z))^{-3/2} = 2 x^{-1} \tau^3$. The stationary phase lemma then gives \eqref{eq:t-asymp},
with $\Psi(y,\tau) = -\tau \ev(\tau) + \Phi(y, \ev(\tau)).$
\end{proof}

\begin{rem} Equation \eqref{energy-eqn} is just the energy equation `total
energy = potential energy+ kinetic energy' at infinity, since $1/\tau$ is
the asymptotic speed. The factor $1/4$ comes from the fact that in writing
our Hamiltonian as $\Delta + V$, we have taken the value of mass to be
$1/2$ in our units.  \end{rem}

\begin{rem}\label{HS} We may \emph{not} replace $\Psi(y,\tau)$ with $\Psi(z, \tau)$ in \eqref{eq:t-asymp}, due to the singular factor $1/x$ in the phase. In fact, if we expand $\Psi(y, \tau)$ in a Taylor series about $y=z$, written in terms of the variables $Y_i = y_i/x^{r_i}$, then we get an asymptotic expansion involving polynomials $y^\alpha$ in the variables $y_i$ multiplied by nonnegative powers $x^r$, where $r = \sum \alpha_i r_i$. We may discard all terms in the Taylor series of $\Psi$ with $r > 1$ since these will only contribute to the term $u'$ decaying faster than the leading term, but we must keep all terms with $r \leq 1$. The number of such terms is always finite, but depends on $\ev - V_0(z)$ and the eigenvalues of the Hessian of $V_0$ at $z$. 
\end{rem}

\subsection{Asymptotic completeness: time dependent formulation}
We see that solutions of the time dependent Schr\"odinger equation (at
least those with initial data in $\mathcal{A}$) have expansions at af which
are equivalent to first spectrally resolving the initial data and looking
at the expansion of the corresponding family of generalized
eigenfunctions; the variable $\ev$ in the time-independent setting, and $\tau$ in the time-dependent setting, play equivalent roles and are linked by \eqref{energy-eqn}. In view of this, we can recast Theorem~\ref{thm:AC} in time-dependent terms as follows:

\begin{thm}\label{thm:time-AC} Let $I$ and $\mathcal{A}_I$ be as in Theorem~\ref{thm:enr-structure-t}, let $u_0 \in \mathcal{A}_I$ and let  $u$ be the solution of the time
dependent Schr\"odinger equation \eqref{time-dept-eqn} with initial data $u_0$. For a given local minimum $z$ of $V_0$, let $\Min_+(I)$ be the associated family of outgoing radial points, and let $\tilde Q = \sum_j \tilde Q_j$, $\tilde b$ and $\kappa_j$ be as in
Proposition~\ref{prop:enr-min}. The map
\begin{equation}
\mathcal{A}_I \ni u_0 \mapsto \oplus_{q \in \Min_+(I)} \Big(e^{i \log x
  \tilde Q}x^{i \tilde b - 1/2} e^{-i\Psi(y,\tau)/x} u(x, \tau, Y'', Y''') \Big)\big|_{x = 0}
\label{ttt}\end{equation} 
whose existence is guaranteed by Theorem~\ref{thm:enr-structure-t}
extends uniquely by linearity and continuity to a unitary isomorphism 
\begin{equation}
L^2 \ominus E_{\pp}(H) \to \oplus_{q} L^2\big(\RR^+_\tau \times
\RR^{n-1}_q; \frac{ d\tau}{2\tau^4} \otimes \omega_{q,\tau} \big).
\label{unit-iso}\end{equation}
Here $\omega_{q, \tau}$ is the measure in \eqref{Hilbert-norm} and $\tau
= \tau(q, \ev)$ is given by \eqref{energy-eqn}. \end{thm}

\begin{rem} The operator $e^{i \log x \tilde Q}$ simply removes the factors
of $x^{-i \sigma_k}$ in the expansion \eqref{eq:t-asymp}, so that we can
take a limit as $x \to 0$.  \end{rem}

\begin{rem} The measure in \eqref{unit-iso} should be thought of as the product of $\tau^{-1} \omega_{q,\tau}$, which is the measure in Proposition~\ref{prop:bdy-pairing}, tensored with the measure $d\ev =  \tau^{-3} d\tau/2$. 
\end{rem} 

\begin{proof} Let $u_0 = \phi(I, h)$ be as in Theorem~\ref{thm:enr-structure-t}. We may take the $L^2$ norm of \eqref{eq:t-asymp} for a fixed $t$, and take the limit as $t \to \infty$. To do this, we write $x = \tau/t$ and integrate with respect to the measure on $X$ which is given by $t \tau^{-2} x^{-1 -2 \Im \tilde b} dY d\tau$. If we just look at the principal term in \eqref{eq:t-asymp} then the powers of $t$ cancel exactly and we get
\begin{equation}
\sum_q \int \frac1{2\pi} \big| \sum_j w_j(Y'', \ev(\tau)) v_j(Y''', \ev(\tau)) \big|^2 \, dY \, \frac{d\tau}{2\tau^4}.
\label{leading-term}\end{equation}
Since $\omega_{q, \ev} = \tau^{-1} dY$, and $d\sigma =  \tau^{-3} d\tau/2$ using \eqref{energy-eqn}, this is given by
$$\begin{gathered}
\frac1{2\pi} \sum_q \int_I  2 \sqrt{\ev - V_0(z)} \big| \sum_j w_j(Y'', \ev) v_j(Y''', \ev) \big|^2 \, d\omega_{q,\ev} \, d\sigma .
\end{gathered}$$
The expression $\sum_j w_j v_j$ is equal to $M_+(q, \ev)((R(\ev + i0) - R(\ev - i0)) u_0$, or equivalently $2\pi i M_+(q, \ev) \Sp(\ev) u_0$. Also, the norm on $\oplus_q L^2(\RR^{n-1})$ is given by \eqref{Hilbert-norm}. So we get, using Theorems \ref{thm:iso} and \ref{thm:AC},
$$
\eqref{leading-term} \ = 2\pi \sum_q \int_I  \| M_+(q,\ev) \Sp(\ev) u_0 \|^2  \, d\sigma 
= \| u_0 \|_{L^2}^2 .
$$
But \eqref{leading-term} is precisely the square of the norm of the right hand side of \eqref{ttt}. 
So we have established the conclusion of the theorem for the principal term in the asymptotic expansion in \eqref{eq:t-asymp}. 
Since the remainder term $u'$ decays faster than the principal term, the $L^2$ norm of $u'(\cdot, t)$ goes to zero as $t \to \infty$, so the proof is complete. 
\end{proof}

{F}rom Theorem~\ref{thm:time-AC}  we can deduce the following result first proved by Herbst and Skibsted, using a direct method involving
the uncertainty principle, rather than proceeding via the structure of
generalized eigenfunctions as here.

\begin{cor}[Absence of $L^2$ channels at non-minimal critical
    points]\label{cor:channels} Let $\chi \in \CI(X)$ vanish in a neighbourhood of the local minima of $V_0$ on $\partial X$. Let $u$ be
  the solution of \eqref{time-dept-eqn} on $X \times \RR$ with initial
  value $u_0 \in L^2(X) \ominus E_{\pp}(H)$. Then 
\begin{equation}
\lim_{t \to \infty} \| \chi u(t, \cdot) \|_{L^2(X)} \to 0.
\label{channel-absence}\end{equation}
\end{cor}

\begin{proof} We may assume that $0 \leq \chi \leq 1$ without loss of generality. 

Let $\epsilon > 0$ be given. Then by
density of $\mathcal{A}$ in $L^2 \ominus E_{\pp}(H)$, we can find $\phi \in
\mathcal{A}$, with $\phi$ equal to a sum of a finite number of $\phi_j(I_j,
h_j) \in \mathcal{A}_{I_j}$, such that $\| u_0 - \phi \|_{L^2} <
\epsilon$. Without loss of generality we may assume that all the $I_j$ are
disjoint. Let $u'$ be the solution with initial condition $\phi$. By
direct calculation from \eqref{eq:t-asymp} we find that

$$
 \lim_{t \to \infty}  \| (1 - \chi) u'(t, \cdot) \|_{L^2}^2 =
\sum_j \int_I    \sqrt{\ev - V_0(\pi(q))} \| h_j \|_{L^2(R^{n-1})}^2 \, d\ev,
$$
which by Theorem~\ref{thm:AC} is equal to $\| \phi \|_{L^2}^2$. But by
 unitarity of $e^{-itH}$, we have
$$
\| u'(t, \cdot) \|_{L^2}^2 = \| \phi \|_{L^2}^2 \text{ for each } t. 
$$
Since $0 \leq \chi \leq 1$, an elementary calculation shows that 
$$
\| (1 - \chi) u' \|_{L^2}^2 + \| \chi u' \|_{L^2}^2 \leq \| u' \|_{L^2}^2,
$$
which implies that 
$$
\lim_{t \to \infty}  \| \chi u'(t, \cdot) \|_{L^2}^2 = 0. 
$$

So \eqref{channel-absence} is true for $u'.$ On the other hand,
$$
\lim\sup_{t \to \infty} \| \chi (u(t, \cdot) - u'(t, \cdot)) \|_{L^2} \leq
\epsilon,
$$
so $\lim \sup_{t \to \infty} \| \chi
u(t, \cdot) \|_{L^2} \leq \epsilon $. Since this is
true for every $\epsilon > 0$, the result follows.
\end{proof}

\subsection{Comparison with results of Herbst-Skibsted}
We first show that our results on the asymptotics of the solutions to the time-dependent equation \eqref{time-dept-eqn} are consistent with the comparison dynamics of Herbst-Skibsted \cite{Herbst-Skibsted3}. Herbst and Skibsted define comparison dynamics, i.e. a family of unitary operators $U_0(t)$ for a given local minimum of $V_0$ and for either a `low energy' range or a `high energy' range which depends on the behaviour of the $r_i(\sigma)$ from Lemma~\ref{HMV2r.191}. It has the property that the strong limit
$$
\lim_{t \to \infty} e^{itH} U_0(t)
$$
exists in $L^2(X)$ and defines a unitary wave operator. 

Let us compare their results on long-time asymptotics with ours. For simplicity, we consider the `very low energy' energy interval in which all of the exponents $r_i$ are complex, with real part $1/2$ (this is `below the Hessian threshold', in our terminology). For simplicity we also assume, as do Herbst and Skibsted, that $V_0(z) = 0$. In this case, the exponent $-i\tilde b$ in \eqref{eq:t-asymp} is equal to $(n-1)/4$, and there are no $Y''$ variables. Moreover, the function $\Phi(y, \sigma)$ is equal to $\sqrt{\sigma}(1 - |y|^2/4)$ (see \cite{HMV1}, section 7, particularly (7.22) and (7.23) for the case $n=2$), which implies that $\Psi(y,\tau) = (1 - |y|^2/4)/\tau$. If we substitute $x = \tau/t$ into \eqref{eq:t-asymp} then we get
\begin{equation}
c\sum_j t^{-(n-1)/4 - 1/2 + i \kappa_j} \tau^{(n-1)/4 +1 - i \kappa_j} e^{it(1 - |y|^2/4)/\tau^2} w_j(\tau) v_j(Y''', \tau).
\end{equation}
To compare this with Herbst and Skibsted's comparison dynamics, we adopt their notation: we decompose the variable $\underline{x} \in \RR^n$ as $\underline{x} = (x_1, x^\perp)$ where $(1, 0, \dots, 0)$ is the point on the sphere at infinity where $V_0$ has a local minimum, and $x^\perp$ are  $n-1$ orthogonal linear coordinates. We can identify our boundary defining function $x$ with $1/x_1$. Thus $\tau = t/x_1$ and $y = x^\perp/x_1$, and we can write the expression above as
\begin{equation}
c\sum_j t^{-(n-1)/4 - 1/2 + i \kappa_j} \tau^{(n-1)/4 +1 - i \kappa_j} e^{it/\tau^2} e^{i|x^\perp|^2/4t} w_j(\tau) v_j(x^\perp/\sqrt{x_1}, \tau).
\label{lta-us}\end{equation}

In this very low-energy case the Herbst-Skibsted comparison dynamics is given explicitly by
\begin{equation}
U_0(t) = S_{t^{-1/2}} e^{i|x^\perp|^2/4} e^{-it p_1^2/2} e^{-i (\log t) H_2} \hat U_0,
\label{HS-low}\end{equation}
where the operator $p_i$ stands for $D_{x_i} = -i \partial_{x_i}$,    $S_{t^{-1/2}}$ is the scaling
$$
S_{t^{-1/2}} f(x_1, x^\perp) = t^{-(n-1)/4} f(x_1, t^{-1/2} x^\perp),
$$
the operator $H_2$ is given by
$$H_2 = \frac{|p^\perp|^2}{2} + \frac1{2} \langle x^\perp, \big( p_1^{-2} V^{(2)} - \Id/4 \big) x^\perp \rangle
$$ (where $V^{(2)}$ is the Hessian of $V_0$ at the critical point), and finally $\hat U_0$ is an arbitrary unitary operator.

To compare this to our long-time asymptotic expansion \eqref{lta-us}, it is convenient to take $\hat U_0$ to be inverse Fourier transform mapping functions of $p_1$ to functions of  $x_1$. Then $H_2$ is a family of harmonic oscillators parametrized by $p_1$. The operator $e^{-itp_1^2/2}$ acting on $W(p_1, x^\perp)$ then takes the form
$$
(2\pi)^{-1/2} \int e^{ix_1 p_1} e^{-it p_1^2/2} W(p_1, x^\perp) \, dp_1
$$
and by stationary phase we see that the large $t$ asymptotics of this operation is given by
$$
W(p_1, x^\perp) \mapsto t^{-1/2} e^{ix_1^2/2t} W(\frac{x_1}{t}, x^\perp).
$$
Let us expand $W(p_1, x^\perp)$ in eigenfunctions of the operator $H_2 = H_2(p_1)$ as
$$
W(p_1, x^\perp) = \sum_j \omega_j(p_1) \chi_j(x^\perp, p_1),
$$
and write $\tau$ for $t/x_1$. A computation shows that $S_{t^{-1/2}} H_2 S_{t^{1/2}}$ is equal to $p_1^{-1} \tilde Q$ where $\tilde Q$ is the operator from \eqref{model-efn}. The comparison dynamics 
 therefore maps $W$ to
\begin{equation}
 t^{-(n-1)/4 - 1/2} e^{it/2\tau^2} e^{-i|x^\perp|^2/4t} \sum_j t^{i \kappa_j} \omega_j(\tau^{-1}) \chi_j(x^\perp/\sqrt{t}, \tau^{-1}). 
\label{lta-them}\end{equation} 
This agrees with \eqref{lta-us}, if we identify $w_j(\tau)$ with $\tau^{-1} \omega_j(\tau^{-1})$, and $v_j(Y''', \tau)$ with $\chi_j(Y''' \tau^{-1/2}, \tau^{-1})$. (The imaginary powers of $\tau$ simply amount to a different choice of normalized eigenfunction $v_j$. Also there are some discrepancies of factors of $1/2$ since Herbst-Skibsted's operator is \emph{half} the Laplacian plus $V$.) Thus, the two expansions are consistent. 

In the high energy regime, it is easier to check the agreement of the two expansions. In this case, there are no $Y'''$ variables. The Herbst-Skibsted comparison dynamics takes the form
$$
\tilde U_0(t)f(\underline{x}) = e^{iS(t, \underline{x})} J(t, \underline{x})^{1/2} f(k(t, \underline{x}), w(t, \underline{x})),
$$
where $S(t,\underline{x})$ is a solution to the eikonal equation
\begin{equation}
\partial_t S(t, \underline{x}) + \frac1{2} |\nabla_{\underline{x}} S(t, \underline{x})|^2 + V(\underline{x}) = 0
\label{eik}\end{equation}
and $k(t, \underline{x})$ is the energy function 
$$
\frac{k^2}{2} = \frac1{2} |\nabla_{\underline{x}} S(t, \underline{x})|^2 + V(\underline{x}).
$$

To make the link with our long time expansion \eqref{eq:t-asymp}, we begin by showing that $S$ corresponds to our phase function $\Psi/x$. Indeed, $\Psi$ is obtained from $\Phi$, the phase function in \eqref{res-exp} by performing stationary phase as in \eqref{spec-int-3}. The phase function $\Phi/x$ parametrizes a Legendrian submanifold which is the image of the zero section under the FIO $F$ in Remark~\ref{original-exp}. Since the zero section is the flowout from the critical point in the eigendirections (of the linearized flow) with eigenvalues $\lambda r_i$ (as opposed to $\lambda (1-r_i)$), as can be computed easily from \eqref{HMV2r.189}, the same is true of the Legendrian submanifold parametrized by $\Phi/x$; in particular, it corresponds precisely to Herbst-Skibsted's Lagrangian submanifold $\mathcal{M}_k$ parametrized by $\overline{S}$ (using the correspondence between conic Lagrangian submanifolds and Legendre manifolds `at infinity'); see \cite{Herbst-Skibsted3}, Theorem 2.1. Then the way $\Psi/x$ is obtained from $\Phi/x$ is exactly the same as the Legendre transform by which Herbst-Skibsted obtain $S$ from $\overline{S}$ (see \cite{Herbst-Skibsted3}, p559), with $k^2$ corresponding to our $\sigma$ and $S$ corresponds to our $\Psi$. Moreover, from \cite{Herbst-Skibsted3}, p561, we have 
$$
w_j = t^{-\beta_j(k)} k^{1-\beta_j(k)} (1 + 2\beta_j(k)) u_j + O(u_j |u|).
$$
In our notation, $\beta_j(k) = -r_j$, $u_j = y_j$ and $t = \tau/x$. Setting $Y''_j = y_j x^{-r_j}$ as above, we get
$$
w_j = g_j(k) Y''_j + O(x^{\min r_j}).
$$
Thus, up to an energy-dependent factor $g_j(k)$, the coordinate $w$ in Herbst-Skibsted is equivalent to our $Y''$. 
The asymptotics \eqref{eq:t-asymp} in this regime 
 (where now there are no $Y'''$ variables, hence no sum over $j$) thus take the form
 $$
c x^{-i \tilde b + 1/2} \tau^{-3/2} e^{i\Psi(y, \sigma(\tau))/x} w(Y'', \sigma)
$$
which is consistent with the Herbst-Skibsted comparison dynamics.

It is a little more difficult to make the link between our asymptotics and
Herbst-Skibsted's in the low-energy regime where not all the $r_i$ have
real part equal to $1/2$, but the real parts are all at least $1/3$.   We can,
however, offer some explanation as to why the low-energy comparison
dynamics fails to work \emph{above} this energy level. Referring to
Remark~\ref{HS}, above this energy we cannot approximate the function
$\Psi(y, \sigma(\tau))$ by its quadratic approximation; we need to include
at least cubic terms in the Taylor series of $\Psi$ at $y=z$. These in turn
depend on the cubic terms in the Taylor series for $V_0$ at $z$. The
Herbst-Skibsted low energy comparison dynamics neglects these terms. It
cannot therefore be expected to provide an accurate approximation to
the long-time asymptotics of solutions to \eqref{time-dept-eqn}, since we
have seen that in \eqref{eq:t-asymp} that one \emph{cannot} replace $\Psi$
by its quadratic approximation.

We emphasize that our long time asymptotic formula \eqref{eq:t-asymp} works
for all energies (except for the discrete set of eigenvalues, effectively resonant energies and Hessian thresholds), whereas in Herbst and
Skibsted's results there is a gap of `intermediate energies' in which they
do not give any comparison dynamics. The formula \eqref{eq:t-asymp}
correctly interpolates between low energies, below the Hessian threshold,
and high energies, where all the exponents $r_i$ are real.

\appendix

\section{Errata for \cite{HMV1}}
\subsection{Correction to the proof of \cite[Proposition~6.7]{HMV1}}

There is an error in the proof of
\cite[Proposition~6.7]{HMV1}:
with the stated assumptions, the proof of Proposition 6.7 needs to be 
two-step, and the conclusion is slightly modified, although this does not
affect any of its applications, in particular
\cite[Proposition~6.9]{HMV1}, which is its
only use in \cite{HMV1}.
Below equation numbers of the form (6.xx) refer to
\cite{HMV1}, while equation numbers
of the form (A.xx) refer to this appendix.

The error in the proof arises from the microlocalizers
$Q\in\Psisc^{-\infty,0}(X)$ considered there,
in (6.27), so we recall the assumptions on it. With $O_m$ a neighborhood
of $\Rp$ as (6.24) or (6.25),
we assume that $\WFscp(Q)\subset O_m$
$\Rp\nin\WFscp(\Id-Q),$
\begin{equation}\begin{split}
&i[Q^*Q,P-\ev]= x^{1/2} (\tilde B^*\tilde B + \tilde G)x^{1/2}
+x^{1/2}\tilde Fx^{1/2}, \text{ where }\\
&\tilde B,\ \tilde F\in\Psisc^{0,0}(O), \ \tilde G \in \Psisc^{0,1}(X),
\ \Rp\nin\WFscp(\tilde F),
\end{split}
\label{eq:[Q,P]-0}\end{equation}
and in addition, $\tilde F$ satisfies
$\WFscp(\tilde F) \subset \{ \nu < \nu(\Rp) \}$. (This condition on $\tilde F$
ensures that $\WFscp(\tilde F) \cap \WFsc(u) = \emptyset$ for the
application in \cite[Section~9]{HMV1}.)

In fact, due to the two step nature of the proof below, we also need
another microlocalizer $Q'\in\Psisc^{-\infty,0}(X)$ satisfying analogous
assumptions with $\tilde B$, etc., replaced by $\tilde B'$, etc.,
\begin{equation}
i[(Q')^*Q',P-\ev]= x^{1/2} ((\tilde B')^*\tilde B' + \tilde G')x^{1/2}
+x^{1/2}\tilde F'x^{1/2},
\label{eq:[Q,P]-p}\end{equation}
with properties analogous to \eqref{eq:[Q,P]-0},
except that $\WFscp(Q')\subset O'_m$, etc., where $O'_m$ is the elliptic
set of $Q$.

The following is a slightly modified version of
Proposition~\ref{prop:HMV1-mod}, in that we need to assume the
existence of $Q'$ as above, and that the conclusion is on
the elliptic set of $Q'$ rather than that of $Q$.

\begin{prop}\label{prop:HMV1-mod}
(modified version of \cite[Proposition~6.7]{HMV1})
Suppose that $m>0,$ $s<-1/2,$
$\Rp\in\RP_+(\ev),$ $\ev\nin\Cv(V),$ either
(6.14) or (6.15) hold, and let
$O_m$ be as in (6.24) (or (6.25)).
Suppose that $u\in\Isc^{(s),m-1}(O_m,\cM),$
$\WFsc((P-\ev)u)\cap O_m=\emptyset$ and that there exists 
$Q,Q'\in\Psisc^{-\infty,0}(O_m)$ elliptic at $\Rp$
that satisfies \eqref{eq:[Q,P]-0}-\eqref{eq:[Q,P]-p} with
$\WFscp(\tilde F)\cap \WFsc(u)=\emptyset,$
$\WFscp(\tilde F')\cap \WFsc(u)=\emptyset.$
Then $u\in\Isc^{(s),m}(O'',\cM)$
where $O''$ is the elliptic set of $Q'.$
\end{prop}

The issue with the argument presented in \cite[Proof of
Proposition~6.7]{HMV1}: is that
it gains a whole extra factor in the module at once:
$u\in\Isc^{(s),m-1}(O_m,\cM),$ is assumed, and $u\in\Isc^{(s),m}(O',\cM)$
is concluded. Now, the novel part of such a statement, corresponding to the
terms arising from factors from the module $\cM\subset\Psisc^{-\infty,-1}(X)$,
is properly dealt with in the (erroneous) proof presented in
\cite{HMV1}.
However, there 
is a problem with the microlocalizer $Q$ unless (6.27) is 
strengthened to make the error term $\tilde G$ have two orders higher 
decay than the main term, i.e. to make it order $(0,2)$.
This is of 
course the same issue as what makes one gain $1/2$ order at a time 
usually in positive commutator proofs for the propagation of singularities
for operators of real principal type.
Factors from the module
$\cM$ are fine because they essentially get 
reproduced by the commutator with $P-\ev$.
The problem is that $\tilde G$ cannot be 
written as a multiple of $Q$, in general. Technically, this shows up in 
(6.29) where 
$\ep\|A_{\alpha,s}u'\|^2$ cannot be absorbed in the left hand side
for it does not 
have a factor of $Q$. (One needs to remember that $Au'$ is the vector of 
$QA_{\alpha,s}u'$, so all terms arising by commutators with the module 
generators are OK, the only issue is the microlocalizer $Q$.)

This error is easily remedied by a two-step argument. The cost of
this is that the open set on which we conclude regularity is shrunk
slightly from the elliptic set of $Q$ to that of $Q'$, although
in relevant situations one can usually recover the original
statement of \cite[Proposition~6.7]{HMV1}
easily as in \cite[Proposition~6.9]{HMV1}.
First, the argument
given in \cite[Proof of Proposition~6.7]{HMV1}
proves the following lemma.

\begin{lemma}\label{lemma:HMV1-mod}
Suppose that $m>0,$ $r<-1/2,$
$\Rp\in\RP_+(\ev),$ $\ev\nin\Cv(V),$ either
(6.14) or (6.15) hold, and let
$O_m$ be as in (6.24) (or (6.25)).
Suppose that $u\in\Isc^{(r),m-1}(O_m,\cM),$
$\WFsc((P-\ev)u)\cap O_m=\emptyset$ and that there exists 
$Q\in\Psisc^{-\infty,0}(O_m)$ elliptic at $\Rp$
that satisfies \eqref{eq:[Q,P]-0} with
$\WFscp(\tilde F)\cap \WFsc(u)=\emptyset.$ Then
$u\in\Isc^{(r-1/2),m}(O',\cM)$
where $O'$ is the elliptic set of $Q.$
\end{lemma}

Notice that under the same hypothesis as Proposition~\ref{prop:HMV1-mod},
this lemma proves regularity under $\cM^m$
(as Proposition~\ref{prop:HMV1-mod}), but does so at the cost of losing
half an order of decay: $u\in\Isc^{(r-1/2),m}(O',\cM)$ rather than
$u\in\Isc^{(r),m}(O',\cM)$.

\begin{proof}[Proof of Lemma~\ref{lemma:HMV1-mod}]
With the notation of
\cite[Proof of Proposition~6.7]{HMV1},
let $s=r-1/2$ (so in particular $s<-1/2$), let $A_{\alpha,s}$, etc., be
as there. Then
the pairing $\langle A_{\alpha,s}u',\tilde GA_{\alpha,s}u'\rangle$ (where 
$u'$ will be regularizations of $u$) is controlled by the a priori control 
of $u'$ in $\Isc^{(s+1/2),m-1}(O_m,\cM)=\Isc^{(r),m-1}(O_m,\cM)$.
Indeed, $x^{1/2}A_{\alpha,s}$ and 
$x^{-1/2}\tilde GA_{\alpha,s}$ are both the product of an element
of $\Psisc^{(0,-s+1/2)}(O_m)$
and $m$ factors in the module $\cM\subset\Psisc^{0,-1}(O_m)$, hence in 
particular can be thought of (by combining the factor from
$\Psisc^{(0,-s+1/2)}(O_m)$ with a factor from $\cM$)
as the product of an element of  $\Psisc^{(0,-s-1/2)}(O_m)$
with $m-1$ factors in $\cM$. So this gives
$u\in 
\Isc^{(s),m}(O',\cM)=\Isc^{(r-1/2),m}(O',\cM)$, proving the lemma.
\end{proof}

We can now prove Proposition~\ref{prop:HMV1-mod}.

\begin{proof}[Proof of Proposition~\ref{prop:HMV1-mod}.]
  Lemma~\ref{lemma:HMV1-mod} shows that
$u\in\Isc^{(s-1/2),m}(O',\cM)$ with $O'$ as in Lemma~\ref{lemma:HMV1-mod}.
With this additional knowledge,
the argument stated in
\cite[Proof of Proposition~6.7]{HMV1},
applied with the same $s$, goes through. (But now we apply it with $Q$
replaced by $Q'$, etc!)
Indeed, the pairing $\langle 
A_{\alpha,s}u',\tilde G'A_{\alpha,s}u'\rangle$
is controlled by the a priori information, as
$x^{1/2}A_{\alpha,s}u'=A_{\alpha,s-1/2}u'$, so it is controled in $L^2$ if 
$u'$ is a priori controlled in $\Isc^{(s-1/2),m}(O',\cM)$ (which 
we just have proved), and a similar conclusion holds for $x^{-1/2}\tilde 
G'A_{\alpha,s}u'$ as $x^{-1/2}\tilde G'\in\Psisc^{0,1/2}(X)$ just like 
$x^{1/2}$ is. Thus, $u\in\Isc^{(s),m}(O'',\cM)$, with $O''$
the elliptic set of $Q'$, as desired. This finishes the proof.
\end{proof}

\subsection{Correction to Proposition 9.4 of \cite{HMV1}}

The proof of Proposition 9.4 contains the statement  `` Since $r_1 < 0$, the vector field
$x \partial_x + r_1 y \partial_y$ is nonresonant'', which is false. To correct the proof,
that statement should be deleted and the sentence following it replaced by: 
``By a change of coordinates $x'=a(y)x,$ $y'=b(y)y$, where $a, b \in C^\infty$ near $y=0$ satisfy the ODEs
$$
a'(y) =-\frac{a(y) F(y)}{r_1 + y G(y)}, \quad b'(y) = -\frac{b(y)}{r_1 + y G(y)}, \quad a(0) = b(0) = 1
$$
the $F$ and $G$ terms are eliminated and 
the vector field becomes
$$
-\frac{2\tilde\nu}{a(y)}(((x')^2D_{x'})+r_1 y(x'D_{y'})),
$$
modulo terms in $x^2\tcM^2$ and subprincipal
terms.'' 
This proves the proposition apart from the prefactor of $a(y)^{-1}$ in front of $\tilde P_0$
which is irrelevant for the application of this proposition. 

Of course, Proposition 9.4 also follows by applying the results of the present paper, noting that the case considered there  is effectively nonresonant.

\end{document}